%% file: TRPCDProbAspect.tex
\documentclass[10pt]{article}

\usepackage[dvips]{epsfig} 
\usepackage{amssymb,amsmath,amsthm,mathrsfs} 
\usepackage{graphicx}
\usepackage{lineno}

\usepackage[authoryear,sort&compress]{natbib}
\usepackage{url}
\usepackage{enumerate}
\usepackage{colordvi,color,pspicture}
\usepackage{psfrag}

\DeclareMathOperator{\argmax}{argmax}
\DeclareMathOperator{\argmin}{argmin}
\DeclareMathOperator{\arginf}{arginf}

\newcommand{\E}[1]{\mathbf{E}\,#1}

\newcommand{\I}{\mathbf{I}}

\newcommand{\C}{\mathcal{C}}

\newcommand{\U}{\mathcal{U}}

\newcommand{\g}{\gamma}
\newcommand{\G}{\Gamma}
\newcommand{\V}{\mathcal{V}}
\newcommand{\A}{\mathcal{A}}
\newcommand{\Y}{\mathcal{Y}}
\newcommand{\y}{\mathsf{y}}
\newcommand{\X}{\mathcal{X}}
\newcommand{\M}{\mathcal{M}}

\newcommand{\mS}{\mathcal{S}}

\newcommand{\NY}{N_{\Y}}
\newcommand{\NS}{N_S}
\newcommand{\NAS}{N_{AS}}
\newcommand{\NCS}{N_{CS}}
\newcommand{\NCSt}{N_{CS}^{\tau}}
\newcommand{\NPE}{N_{PE}}
\newcommand{\TY}{T\left(\Y_3\right)}

\newcommand{\R}{\mathbb{R}}
\newcommand{\T}{\mathcal{T}}

\newcommand{\RS}{\mathscr{R}_S}
\newcommand{\UT}{\U(\TY)}
\newcommand{\Tr}{\mathscr T^r}
\newcommand{\ve}{\varepsilon}

\theoremstyle{plain}
\newtheorem{theorem}{Theorem}[section]
\newtheorem{lemma}[theorem]{Lemma}
\newtheorem{proposition}[theorem]{Proposition}
\newtheorem{corollary}[theorem]{Corollary}

\theoremstyle{definition}
\newtheorem{definition}[theorem]{Definition}

\theoremstyle{remark}

\newtheorem{remark}[theorem]{Remark}

\begin{document}


\title{Technical Report \# KU-EC-09-3:\\
A Probabilistic Characterization of Random Proximity Catch Digraphs and the Associated Tools}
\author{
Elvan Ceyhan
\thanks{Address:
Department of Mathematics, Ko\c{c} University, 34450 Sar{\i}yer, Istanbul, Turkey.
e-mail: elceyhan@ku.edu.tr, tel:+90 (212) 338-1845, fax: +90 (212) 338-1559.
}
}

\date{\today}
\maketitle


\begin{abstract}
\noindent
Proximity catch digraphs (PCDs) are based on proximity maps which yield proximity regions
and are special types of proximity graphs.
PCDs are based on the relative allocation of points from two or
more classes in a region of interest and have applications in various fields.
In this article, we provide auxiliary tools for and various characterizations of PCDs
based on their probabilistic behavior.
We consider the cases in which the vertices of the PCDs come from uniform
and non-uniform distributions in the region of interest.
We also provide some of the newly defined proximity maps as illustrative examples.
\end{abstract}

\noindent
{\it Keywords:}
class cover catch digraph (CCCD); central similarity PCD; Delaunay triangulation;
domination number; proportional-edge PCD; proximity graph; random graph; relative arc density




\newpage


\section{Introduction}
\label{sec:intro}
The proximity catch digraphs (PCDs) are a special type of proximity graphs
which are based on proximity maps and are used in disciplines where
shape and structure are crucial.
Examples include computer vision
(dot patterns), image analysis, pattern recognition (prototype selection),
geography and cartography, visual perception, biology, etc.
\emph{Proximity graphs} were first introduced by \cite{toussaint:1980},
who called them \emph{relative neighborhood graphs}.
The notion of relative neighborhood graph has been
generalized in several directions
and all of these graphs are now called proximity graphs.
From a mathematical and algorithmic point of view,
proximity graphs fall under the category of \emph{computational geometry}.

In recent years, a new classification and spatial pattern analysis approach
which is based on the relative positions
of the data points from various classes has been developed.
\cite{priebe:2001} introduced the class cover catch digraphs (CCCDs)
and gave the exact and the asymptotic distribution of the domination number of the CCCD
based on two data sets $\X_n$ and $\Y_m$ both of which are
random samples from uniform distribution on a compact interval in $\R$.
\cite{devinney:2002a}, \cite{marchette:2003}, \cite{priebe:2003a},
\cite{priebe:2003b}, and \cite{devinney:2006}
applied the concept in higher dimensions and demonstrated relatively
good performance of CCCD in classification.
The employed methods involve data reduction (condensing)
by using approximate minimum dominating sets as prototype sets,
since finding the exact minimum dominating set is
in general an NP-hard problem
--- in particular, for CCCD --- (see \cite{devinney:Phd-thesis}).
Furthermore, the exact and
the asymptotic distribution of the domination number of the CCCDs are
not analytically tractable in higher dimensions.
\cite{ceyhan:Phd-thesis} extended the concept of CCCDs
by introducing PCDs, which do not suffer from
some of the shortcomings of CCCDs in higher dimensions.
In particular,
two new types of PCDs (namely, \emph{proportional-edge} and
\emph{central similarity PCDs}) are introduced;
distribution of the domination number of proportional-edge PCDs is calculated,
and is applied in testing spatial patterns of segregation and association
(\cite{ceyhan:dom-num-NPE-SPL,ceyhan:dom-num-NPE-MASA}).
The distributions of the relative arc density of these PCD families
are also derived and used for the same purpose
(\cite{ceyhan:arc-density-PE} and \cite{ceyhan:arc-density-CS}).

A general definition of proximity graphs is as follows:
Let $V$ be any finite or infinite set of points in $\R^d$. Each (unordered) pair of points
$(p,q) \in V \times V$ is associated with a neighborhood
$\mathfrak N(p,q) \subseteq \R^d$.  Let $\mathfrak P$ be a property
defined on $\mathfrak N=\{\mathfrak N(p,q):\; (p,q) \in V \times V\}$.
A \emph{proximity} (or \emph{neighborhood}) \emph{graph} $G_{\mathfrak N,\mathfrak P}(V,E)$
defined by the property $\mathfrak P$ is a graph with the set of vertices
$V$ and the set of edges $E$ such that $(p,q) \in E$ iff $\mathfrak N(p,q)$
satisfies property $\mathfrak P$.
Examples of most commonly used proximity graphs are the Delaunay tessellation,
the boundary of the convex hull, the Gabriel graph, relative neighborhood
graph, Euclidean minimum spanning tree, and sphere of influence graph
of a finite data set. See, e.g., \cite{jaromczyk:1992}.

The relative allocation of the data points are used to construct a proximity digraph.
A \emph{digraph} is a directed graph, i.e., a graph with directed edges from
one vertex to another based on a binary relation. Then the pair
$(p,q) \in V \times V$ is an ordered pair and $(p,q)$ is an arc (directed edge)
denoted $pq$ to reflect its difference from an edge.
For example, the nearest neighbor (di)graph in \cite{paterson:1992} is a proximity digraph.
 The nearest neighbor digraph, denoted $NND(V)$, has the vertex set $V$
and $pq$ an arc iff $d(p,q) = \min_{v \in V \setminus\{p\}} d(p,v)$.
That is, $pq$ is an arc of $NND(V)$ iff $q$ is a nearest neighbor of $p$.
Note that if $pq$ is an arc in $NND(V)$, then $(p,q)$ is an edge in $RNG(V)$.
Our PCDs are based on the property $\mathfrak P$ that is
determined by the following mapping which is defined in a more general
space than $\R^d$.
Let $(\Omega,\M)$ be a measurable space. The \emph{proximity map} $N(\cdot)$
is given by  $N:\Omega \rightarrow \wp(\Omega)$,
where $\wp(\cdot)$ is the power set functional, and the
\emph{proximity region} of $x \in \Omega$, denoted $N(x)$, is the image
of $ x \in \Omega$ under $N(\cdot)$.
 The points in $N(x)$ are thought of as being ``closer" to $x \in \Omega$
than are the points in $\Omega \setminus N(x)$. Proximity maps are the
building blocks of the \emph{proximity graphs} of \cite{toussaint:1980};
an extensive survey is available in \cite{jaromczyk:1992}.

The \emph{PCD} $D=(\V,\A)$ has the vertex set
$\V=\bigl\{ p_1,p_2,\ldots,p_n \bigr\}$ and
the arc set $\A$ is defined by $p_ip_j \in \A$ iff $p_j \in N(p_i)$ for $i\not=j$.
Notice that $D$ depends on the \emph{proximity} map $N(\cdot)$,
and if $p_j \in N(p_i)$, then $N(p_i)$ is said to \emph{catch} $p_j$.
Hence the name \emph{proximity catch digraph}.
If arcs of the form $p_ip_i$ (i.e., loops) were allowed,
$D$ would have been called a \emph{pseudodigraph} according to some
authors (see, e.g.,  \cite{chartrand:1996}).

In this article,
we provide a probabilistic characterization of the proximity maps,
and the associated regions and PCDs,
and introduce auxiliary tools for the PCDs.
We define the proximity maps and
data-random PCDs in Section \ref{sec:PCD},
describe the auxiliary tools (such as edge and vertex regions)
for the construction of PCDs in Section \ref{sec:prox-maps-prelim},
provide $\G_1$-regions and the related concepts for
proximity maps in Section \ref{sec:gamma1-regions},
discuss the examples of proximity maps in
Delaunay triangles in Section \ref{sec:example-PCDs},
and provide the transformations preserving uniformity
on triangles in $\R^2$ in Section \ref{sec:transformations}.
We investigate the characterization of proximity regions and
the associated PCDs in Section \ref{sec:characterize-PCDs},
introduce $\G_k$-regions for proximity maps in $\TY$ in Section \ref{sec:Gammak-regions},
$\kappa$-values for the proximity maps in $\TY$ in Section \ref{sec:kappa-value},
and provide discussion and conclusions in Section \ref{sec:disc-conc}.

\section{Proximity Maps and Data-Random PCDs}
\label{sec:PCD}
Let 
$\X_n=\bigl\{ X_1,X_2,\ldots,X_n \bigr\}$ and
$\Y_m=\bigl\{ Y_1,Y_2,\ldots,Y_m \bigr\}$ be two data sets from classes $\X$ and $\Y$ of
$\Omega$-valued random variables whose joint pdf is $f_{X,Y}$.
Let $d(\cdot,\cdot):\Omega \times \Omega \rightarrow [0,\infty)$
be a distance function.  The class cover problem for a target class,
say $\X_n$, refers to finding a collection of neighborhoods,
$N(X_i)$ around $X_i \in \X_n$ such that
(i) $\X_n \subset \bigl( \bigcup_i N(X_i) \bigr)$ and
(ii) $\Y_m \cap \bigl( \bigcup_i N(X_i)\bigr) =\emptyset$.
A collection of neighborhoods satisfying both conditions is
called a {\em class cover}.  A cover satisfying condition
(i) is a {\em proper cover} of class $\X$ while a collection satisfying
condition (ii) is a {\em pure cover} relative to class $\Y$.
From a practical point of view, for example for classification,
of particular interest are the class covers satisfying
both (i) and (ii) with the smallest collection of neighborhoods,
i.e., minimum cardinality cover.
This class cover problem is a generalization of the set cover
problem in \cite{garfinkel:1972} that emerged in statistical pattern
recognition and machine learning, where an edited or condensed
set (prototype set) is selected from $\X_n$ (see, e.g., \cite{devroye:1996}).

In particular, we construct the proximity regions using data sets from two classes.
Given $\Y_m \subseteq \Omega$,
the {\em proximity map} $\NY(\cdot)$ 
associates a {\em proximity region} $\NY(x) \subseteq \Omega$ with each point
$x \in \Omega$. The region $\NY(x)$ is defined in terms of the distance
between $x$ and $\Y_m$.  More specifically, our proximity maps will be
based on the relative position of points from class $\X$ with respect to
the Delaunay tessellation of the class $\Y$.
See \cite{okabe:2000} and \cite{ECarXivPCDGeo:2009} for more on Delaunay tessellations.

If $\X_n$ is a set of $\Omega$-valued
random variables then $\NY(X_i)$ are random sets.
If $X_i$ are independent identically distributed then so
are the random sets $\NY(X_i)$.
We define the data-random PCD $D=(\V,\A)$ ---
associated with $\NY(\cdot)$ ---
with vertex set $\V=\X_n$ and arc set $\A$ by
$X_iX_j \in \A \iff X_j \in \NY(X_i)$.
Since this relationship is not symmetric,
a digraph is needed rather than a graph.
The random digraph $D$ depends on
the (joint) distribution of the $X_i$ and on the map $\NY(\cdot)$.

The PCDs are closely related to the
{\em proximity graphs} of \cite{jaromczyk:1992}
and might be considered as a special case of {\em covering sets}
of \cite{tuza:1994} and {\em intersection digraphs} of
\cite{sen:1989}.  This data random proximity digraph is a
{\em vertex-random proximity digraph} which is not of standard type.
The randomness of the PCDs lies in the fact that the vertices
are random with joint pdf $f_{X,Y}$, but arcs $X_iX_j$ are
deterministic functions of the random variable $X_j$ and the set $\NY(X_i)$.

For example, the CCCD of \cite{priebe:2001}
can be viewed as an example of PCD with $\NY(x)=B(x,r(x))$,
where $r(x):=\min_{\y \in \Y_m}d(x,\y)$.
The CCCD is the digraph of order $n$ with vertex set $\X_n$ and
an arc from $X_i$ to $X_j$ iff $X_j \in B(X_i,r(X_i))$.
That is, there is an arc from $X_i$ to $X_j$ iff there exists an open ball
centered at $X_i$ which is ``pure" (or contains no elements) of $\Y_m$,
and simultaneously contains (or ``catches") point $X_j$.


\section{Auxiliary Tools for the Construction of PCDs in $\R^d$}
\label{sec:prox-maps-prelim}
Recall the proximity map (associated with CCCD) in $\R$ is defined as $B(x,r(x))$
where $r(x)=\min_{y\in \Y_m}d(x,y)$ with $d(x,y)$
being the Euclidean distance between $x$ and $y$ (\cite{priebe:2001}).
Our goal is to extend this idea to higher dimensions and investigate the associated digraph.
Now let $\Y_m=\left \{\y_1,\y_2,\ldots,\y_m \right\} \subset \R^d$.
For $d=1$ the proximity map associated with CCCD is defined
as the open ball $\NS(x):=B(x,r(x))$ for all $x \in \R \setminus \Y_m$
and for $x \in \Y_m$, define $\NS(x)=\{x\}$.
Furthermore, dependence on $\Y_m$ is through $r(x)$.
Hence $\NS(x)$ is based on the intervals
$I_{i-1}=\left(\y_{(i-1):m},\y_{i:m}\right)$ for $i=1,2,\ldots,(m+1)$ with
$\y_{0:m}=-\infty$ and $\y_{(m+1):m}=\infty$
where $\y_{i:m}$ is the $i^{th}$ order statistic in $\Y_m$.
This intervalization can be viewed as a tessellation
since it partitions $C_H\left(\Y_m\right)$, the convex hull of $\Y_m$.
For $d>1$, a natural tessellation that
partitions $C_H\left(\Y_m\right)$ is the Delaunay tessellation
(see \cite{okabe:2000} and \cite{ECarXivPCDGeo:2009}).
Let $\T_i$ for $i=1,2,\ldots,J$ be the $i^{th}$
Delaunay cell in the Delaunay tessellation of $\Y_m$.
In $\R$, we implicitly use the cell that contains $x$ to
define the proximity map.

A natural extension of the proximity region $\NS(x)$
to multiple dimensions (i.e., to $\R^d$ with $d>1$)
is obtained by the same definition as above;
that is, $\NS(x):=B(x,r(x))$ where $r(x):=\min_{\y \in \Y_m} d(x,\y)$.
Notice that a ball is a sphere in higher dimensions,
hence the name \emph{spherical proximity map} and the notation $\NS$.
The spherical proximity map $\NS(x)$ is well-defined for all
$x \in \R^d$ provided that $\Y_m \not= \emptyset$.
Extensions to $\R^2$ and higher dimensions with the spherical proximity
map --- with applications in classification --- are
investigated in  \cite{devinney:2002a}, \cite{devinney:2002b}, \cite{marchette:2003},
\cite{priebe:2003b}, \cite{priebe:2003a}, and \cite{devinney:2006}.
However, finding the minimum
dominating set of the PCD associated with $\NS(\cdot)$ is an NP-hard problem and
the distribution of the domination number
is not analytically tractable for $d>1$ (\cite{ceyhan:Phd-thesis}).
This drawback has motivated us to define new types of proximity maps
in higher dimensions.
Note that for $d=1$, such problems do not occur.
\cite{ECarXivPCDGeo:2009} states some appealing properties of the proximity map
$\NS(x)=B(x,r(x))$ in $\R$ and uses them as guidelines for
extending proximity maps to higher dimensions and defining new proximity maps.
After a slight modification,
the spherical proximity maps gives rise to arc-slice proximity maps which is defined as
$\NAS(x):=B(x,r(x)) \cap \TY$ for $x \in \TY$ (i.e., when the arc-slice proximity region
is the spherical proximity region restricted to the Delaunay triangle $x$ lies in).
However, for $x \not\in C_H(\Y_m)$
(i.e., $x$ is not in any of the Delaunay triangles based on $\Y_m$),
$\NAS(x)$ is not defined.

For $x \in I_i$, $\NS(x)=I_i$ iff
$x=\left(\y_{(i-1):m}+\y_{i:m}\right)/2$.
We define an associated region for such points in the general context.
\begin{definition}
The \emph{superset region} for any proximity map $N(\cdot)$
in $\Omega$ is defined to be $\RS(N):=\bigl\{ x \in \Omega: N(x) = \Omega \bigr\}$.
When $X$ is $\Omega$-valued random variable, then we assume $X \in \RS(N)$
if $N(X) = \Omega$ a.s.
$\square$
\end{definition}

For example, for $\Omega=I_i \subsetneq \R$ with $i=1,2,\ldots,(m-1)$,
$\RS(\NS):=\{x \in I_i: \NS(x) = I_i\}=\left\{ \left(\y_{(i-1):m}+\y_{i:m}\right)/2 \right\}$,
and for $i=0,m$ (i.e., $\Omega=I_0$ or $\Omega=I_m$,
then $\RS(\NS)=\emptyset$ since $\NS(x) \subsetneq I_i$ for all $x \in I_i$ for $i=0,m$.
More generally for $\Omega=\T_i \subsetneq \R^d$ (i.e., $i^{th}$ Delaunay cell),
$\RS(\NS):=\{x \in \T_i: \NS(x) = \T_i\}$.
Note that for $x \in I_i$, $\lambda(\NS(x)) \le \lambda(I_i)$
and $\lambda(\NS(x)) = \lambda(I_i)$ iff $x \in \RS(\NS)$
where $\lambda(\cdot)$ is the Lebesgue measure on $\R$
(also called as $\R$-Lebesgue measure).
So the proximity region of a point in $\RS(\NS)$ has the largest
$\R$-Lebesgue measure.
Note that for $\Y_m=\left\{ \y_1,\y_2,\ldots,\y_m \right\} \subset \R$ (i.e. $\Omega = \R$),
$\RS(\NS)=\emptyset$, since $\NS(x)\subseteq I_i$ for all $x \in I_i$
so $\NS(x)\subsetneq \R$ for all $x \in \R$.
Note also that given $\Y_m$, $\RS(\NS)$ is not a random set,
but $\I(X\in \RS(\NS))$ is a random variable.

\subsection{Vertex and Edge Regions}
\label{sec:vertex-edge-regions}
In $\R$, the spherical proximity maps are defined as open intervals
where one of the endpoints is in $\Y_m$.
In particular, for $x \in I_i=(\y_{i:m},\y_{(i+1):m})$ for $i=1,2,\ldots,m$,
$\NS(x)=(\y_{i:m},\y_{i:m}+2r(x))$ for all $x \in (\y_{i:m},\y_{(i-1):m}+\y_{i:m})/2)$
where $r(x)=d(x,\y_{i:m})$
and
$\NS(x)=(\y_{(i+1):m},\y_{(i+1):m}-2r(x))$ for all $x \in (\y_{(i-1):m}+\y_{i:m})/2,\y_{(i+1):m})$
where $r(x)=d(x,\y_{(i+1):m})$.
Hence there are two subinterval in $I_i$ each touching an edge and the midpoint of the interval,
and $\NS(x)$ depends on which of these regions $x$ lies in.

In $\R^d$ with $d>1$, intervals become Delaunay tessellations,
and our proximity maps are based on the Delaunay cell
$\T_i$ that contains $x$.
The region $\NY(x)$ will also depend on the location of $x$
in $\T_i$ with respect to
the vertices or faces (edges in $\R^2$) of $\T_i$.
Hence for $\NY(x)$ to be well-defined, the vertex or face of
$\T_i$ associated with $x$
should be uniquely determined.
This will give rise to two new concepts:
\emph{vertex regions} and face regions (\emph{edge regions} in $\R^2$).

\begin{figure}[ht]
\begin{center}
\psfrag{A}{\scriptsize{$\y_1$}}
\psfrag{B}{\scriptsize{$\y_2$}}
\psfrag{C}{\scriptsize{$\y_3$}}
\psfrag{R(A)}{\scriptsize{$R_{CM}(\y_1)$}}
\psfrag{R(B)}{\scriptsize{$R_{CM}(\y_2)$}}
\psfrag{R(C)}{\scriptsize{$R_{CM}(\y_3)$}}
\psfrag{CM}{\scriptsize{$M_{C}$}}
\psfrag{P1}{\scriptsize{$M_3$}}
\psfrag{P2}{\scriptsize{$M_1$}}
\psfrag{P3}{\scriptsize{$M_2$}}
\psfrag{x}{}
\epsfig{figure=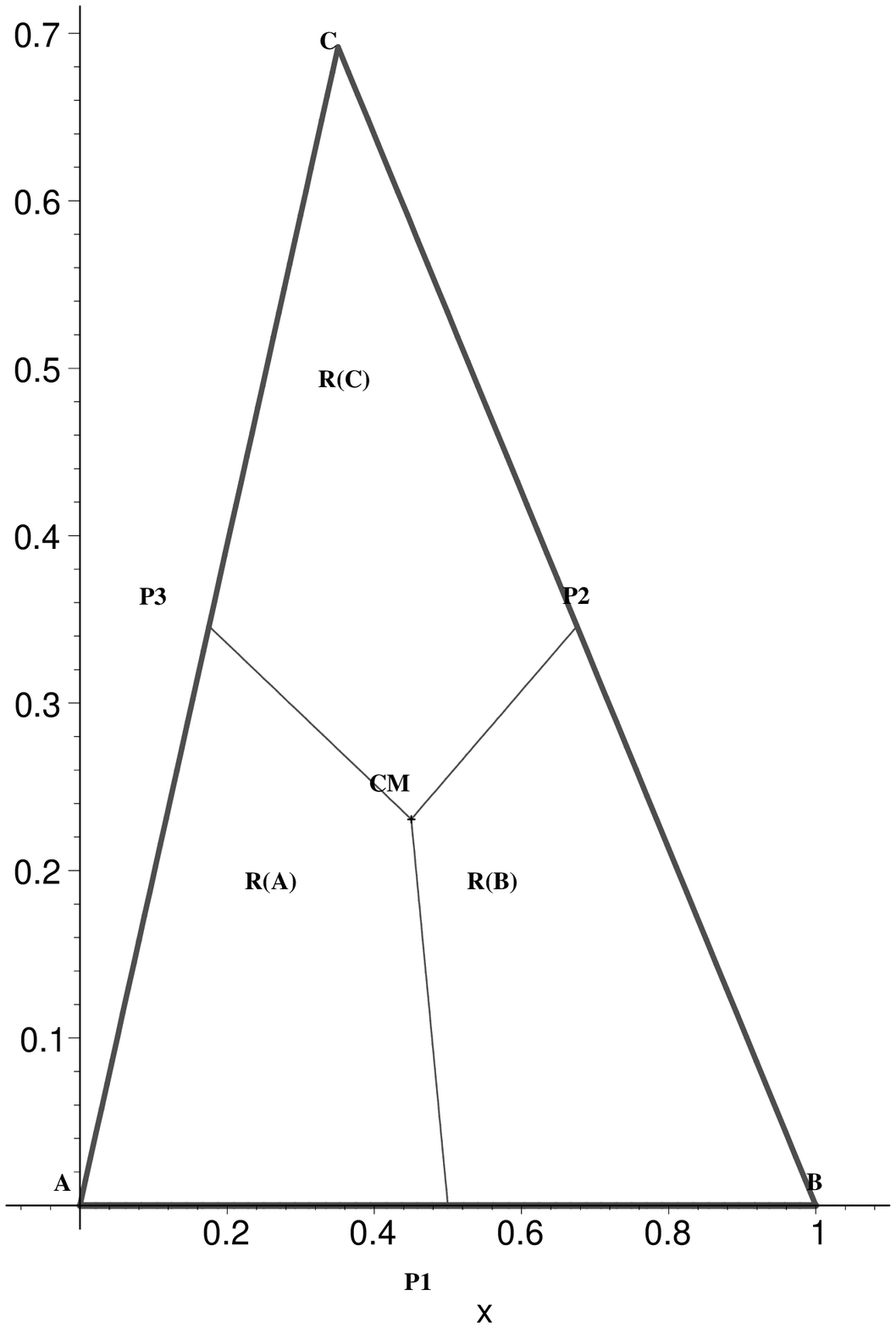, height=140pt , width=200pt}
\psfrag{R(AB)}{\scriptsize{$R_{CM}(e_3)$}}
\psfrag{R(BC)}{\scriptsize{$R_{CM}(e_1)$}}
\psfrag{R(AC)}{\scriptsize{$R_{CM}(e_2)$}}
\epsfig{figure=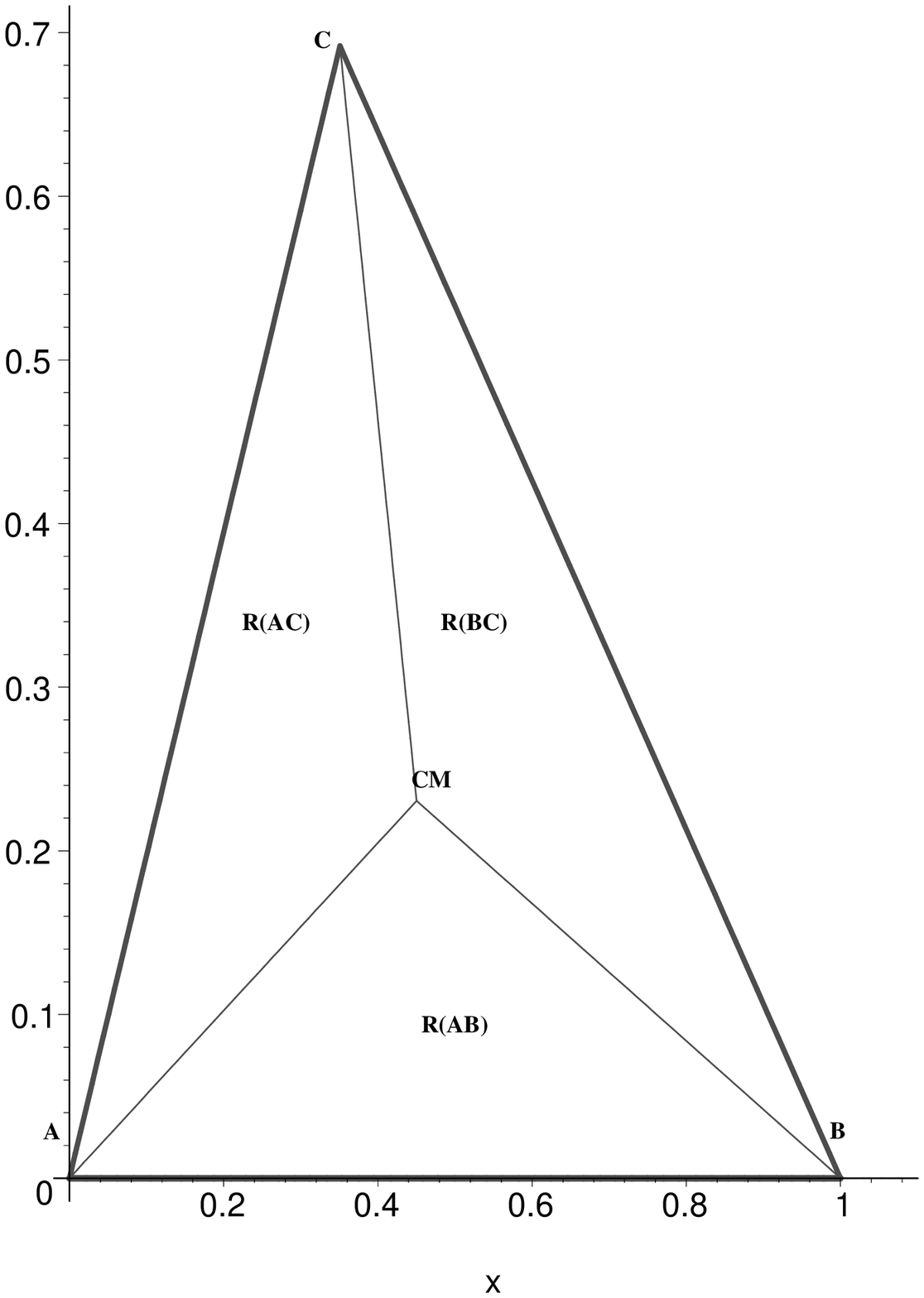, height=140pt , width=200pt}
\end{center}
\caption{
\label{fig:vertex-edge-regions}
The $CM$-vertex regions (left)
and $CM$-edge regions (right) with median lines.}
\end{figure}

Let $\Y_3=\{\y_1,\y_2,\y_3\}$ be three non-collinear points
in $\R^2$ and $\TY=T(\y_1,\y_2,\y_3)$ be the triangle
with vertices $\Y_3$.
To define new proximity regions based on some sort of distance
or dissimilarity relative to the vertices $\Y_3$, we
associate each point in $\TY$ to a vertex of $\TY$.
This gives rise to the concept of \emph{vertex regions}.
\begin{definition}
The connected regions that partition the triangle, $\TY$,
(in the sense that the pairwise intersections of the regions have zero $\R^2$-Lebesgue measure)
 such that each region has one and only one vertex of $\TY$ on its boundary
are called \emph{vertex regions}.
$\square$
\end{definition}
This definition implies that we have three vertex regions.
In fact, we can describe the vertex regions
starting with a point $M \in \R^2 \setminus \Y_3$ as follows.
Join the point $M$ to a point on each edge by a curve such
that the resultant regions satisfy the above definition.
We call such regions \emph{$M$-vertex regions} and
denote the vertex region associated with
vertex $\y$ as $R_M(\y)$ for $\y \in \Y_3$.
Vertex regions can be defined using any point $M \in \R^2 \setminus \Y_3$
by joining $M$ to a point on each edge.
In particular, we use a \emph{center} of $\TY$ as
the starting point $M$ for vertex regions.
See the discussion of triangle centers in (\cite{ECarXivPCDGeo:2009})
with relevant references.
We think of the points in $R_M(\y)$ as
being ``closer" to $\y$ than to the other vertices.
It is reasonable to require that the area of the region $R_M(\y)$ gets larger
as $d(M,\y)$ increases.
Unless stated otherwise, $M$-vertex regions will refer to
regions constructed by joining $M$ to the edges
with \emph{straight line segments}.
Vertex regions with circumcenter, incenter, and center of mass
are investigated in \cite{ECarXivPCDGeo:2009}.
For example $M$-vertex regions can be constructed with $M \in \TY^o$ by
using the extensions of the line segments joining $\y$ to $M$ for all $\y \in \Y_3$.
See Figure \ref{fig:vertex-edge-regions} (left) with $M=M_C$.

We can also view the endpoints of the interval $I_i$ as the edges of the interval
which suggests the concept of \emph {edge regions}.
\begin{definition}
The connected regions that partition the triangle, $\TY$, in such a way that
each region has one and only one edge of $\TY$ on its boundary,
are called \emph{edge regions}.
$\square$
\end{definition}
This definition implies that we have exactly three edge regions.
In fact, we can describe the edge regions starting with $M$ in $\TY^o$, the interior of $\TY$.
Join the point $M$ to the vertices by curves such that
the resultant regions satisfy the above definition.
We call such regions \emph{$M$-edge regions} and denote
the region for edge $e$ as $R_M(e)$ for $e \in \{e_1,e_2,e_3\}$.
Unless stated otherwise, $M$-edge regions will refer to the regions
constructed by joining $M$ to the vertices by straight lines.
In particular, we use a \emph{center} of $\TY$
for the starting point $M$ as the edge regions.
See Figure \ref{fig:vertex-edge-regions} (right) with $M=M_C$.
We can also consider the points in $R_M(e)$ to be
``closer" to $e$ than to the other edges.
Furthermore, it is reasonable to require that the area of
the region $R_M(e)$ gets larger as $d(M,e)$ increases.
Moreover, in higher dimensions, the corresponding regions are called ``face regions".

%
%
Edge regions for incenter, center of mass, and orthocenter are
investigated in (\cite{ECarXivPCDGeo:2009}).

\section{$\G_1$-Regions for Proximity Maps and the Related Concepts}
\label{sec:gamma1-regions}
For any set $B \subseteq \Omega$,
the {\em $\G_1$-region} of $B$ associated with $N(\cdot)$,
is defined to be the region $\G_1(B,N):=\{z \in \Omega: B \subseteq  N(z)\}$.
For $x \in \Omega$, we denote $\G_1\bigl( \{x\},N \bigr)$ as $\G_1(x,N)$.
Note that $\G_1$-region is based on the proximity region $N(\cdot)$.
If $\X_n=\bigl\{ X_1,X_2,\ldots,X_n \bigr\}$ is a set
of $\Omega$-valued random variables,
then $\G_1\left(X_i,N\right)$, $i=1,2,\ldots,n$ are random sets.
If the $X_i$ are independent and identically distributed,
then so are the random sets $\G_1\left(X_i,N\right)$.
Additionally, $\G_1\left(\X_n,N\right)$ is also a random set.

In a digraph $D=(\V,\A)$, a vertex $v \in \V$ \emph{dominates}
itself and all vertices of the form $\{u: (v,u) \in \A\}$.
A \emph{dominating set} $S_D$ for the digraph $D$ is a subset of $\V$
such that each vertex $v \in \V$ is dominated by a vertex in $S_D$.
A \emph{minimum dominating set} $S^*_{D}$ is a dominating set of
minimum cardinality and the \emph{domination number} $\g(D)$ is
defined as $\g(D):=|S^*_{D}|$ (see, e.g., \cite{lee:1998}) where
$|\cdot|$ denotes the set cardinality functional.

For $X_1,X_2,\ldots,X_n \stackrel{iid}{\sim} F$
the domination number of the associated data-random PCD,
denoted $\g_n(N)$, is the minimum number of points
that dominate all points in $\X_n$.
 Note that, $\g_n(N)=1$ iff $\X_n \cap \G_1\left(\X_n,N\right)\not= \emptyset$.
Hence the name \emph{$\G_1$-region}.
Suppose $\mu$ is a measure on $\Omega$.
Following are some general results about $\G_1(\cdot,N)$.

\begin{proposition}
\label{prop:RSsubsetG1}
For any proximity map $N$ and set $B \subseteq \Omega$, $\RS(N) \subseteq \G_1(B,N)$.
\end{proposition}
\noindent \textbf{Proof:}
For $x \in \RS(N)$,
$N(x) = \Omega$, so $B \subseteq N(x)$ since $B \subseteq \Omega$.
Then $x \in \G_1(B,N)$,
hence $\RS(N) \subseteq \G_1(B,N)$. $\blacksquare$

\begin{lemma}
\label{lem:gamma1-intersect}
For any proximity map $N$ and $B \subseteq \Omega$, $\G_1(B,N)= \bigcap_{x \in B}\G_1(x,N)$.
\end{lemma}
\noindent \textbf{Proof:}
Given a proximity map $N$ and subset $B \subseteq \Omega$,
$y \in \G_1(B,N)$ iff $B \subseteq N(y)$ iff $x \in N(y)$ for all $x \in B$
iff $y \in \G_1(x,N)$ for all $x \in B$ iff $y \in \bigcap_{x \in B}\G_1(x,N)$.
Hence the result follows. $\blacksquare$

\begin{corollary}
\label{cor:gamma1-intersect}
For any proximity map $N$ and a realization $\X_n=\bigl\{ x_1,x_2,\ldots,x_n \bigr\}$ from $F$
with support $\mS(F) \subseteq \Omega$, $\G_1\left(\X_n,N\right)= \bigcap_{i=1}^{n}\G_1\left(x_i,N\right)$.
\end{corollary}


A problem of interest is finding, if possible, a subset of $B$,
say $G \subseteq B$, such that $\G_1(B,N)=\bigcap_{x \in G}\G_1(x,N)$.
This implies that only the points in $G$ are used in determining $\G_1(B,N)$.
\begin{definition}
\label{def:active-set}
An \emph{active set} of points $S_A(B) \subseteq \Omega$
for determining $\G_1(B,N)$
is defined to be a subset of $B$ such that
$\G_1(B,N)=\bigcap_{x \in S_A(B)}\G_1(x,N)$.
$\square$
\end{definition}
This definition allows $B$ to be an active set,
which always holds by Lemma \ref{lem:gamma1-intersect}.
If $B$ is a set of finitely many points, so is the associated active set.
Among the active sets, we seek an active set of minimum cardinality.
\begin{definition}
\label{def:min-active-set}
Let $B$ be a set of finitely many points.
An active subset of $B \subset \Omega$ is called a \emph{minimal active subset},
denoted $S_{\mu}(B)$, if there is no other active subset $S_A$
of $B$ such that $S_A(B) \subsetneq S_{\mu}(B)$.
The minimum cardinality among the active subsets of $B$  is
called the  \emph{$\eta$-value} and denoted as $\eta(B,N)$.
An active subset of cardinality $\eta(B,N)$ is called a
\emph{minimum active subset} denoted as $S_M(B)$; that is,
$\eta(B,N):=|S_M(B)|$. $\square$
\end{definition}
Note that Definitions \ref{def:active-set} and \ref{def:min-active-set}
can be extended for any subset $B \subseteq \Omega$, in a similar fashion.
Moreover, a minimal active set of minimum cardinality  is a minimum active set.
We will suppress the dependence on $B$ for $S_A(B)$, $S_{\mu}(B)$, and $S_M(B)$
if there is no ambiguity.
In particular, if $B=\X_n$ is a set of $\Omega$-valued random variables,
then $S_A$ and $S_M$ are random sets and $\eta_n(N)$ is a random quantity.

For example, in $\R$ with $\Y_2=\{0,1\}$, and $\X_n$ a random sample
(i.e., set of iid random variables) of size $n>1$ from
$F$ whose support is in $(0,1)$,
$\displaystyle \G_1\left(\X_n,\NS\right)=\left(\frac{X_{n:n}}{2},\frac{1+X_{1:n}}{2} \right)$,
where $X_{i:n}$ is the $i^{th}$ largest value in $\X_n$.
So the extrema (minimum and maximum) of the set
$\X_n$ are sufficient to determine the $\G_1$-region;
i.e., $S_M=\bigl\{ X_{1:n},X_{n:n} \bigr\}$.
Then $\eta_n\left(\NS\right)=1+\I(n>1)$ a.s.
for $\X_n$ being a random sample from a continuous distribution
with support in $(0,1)$.

In the multidimensional case there is no natural extension of
ordering that yields natural extrema such as minimum or maximum.
Some extensions of ordering are proposed under the title of
``statistical depth" (see for example \cite{liu:1999})
which is not pursued here.
To get the minimum active sets associated with our proximity maps,
we will resort to some other type of extrema, such as,
the closest points to edges or vertices in $\TY$.

For any proximity map $N$ and $\X_n$, $\eta_n(N)\le n$ follows trivially,
since
$$ \eta_n(N):=\min_{A \subseteq \mathcal X_n}
\left\{|A|: \G_1(\mathcal X_n,N)=\bigcap_{Z \in A}\G_1(Z)\right\}.$$

\begin{lemma}
\label{lem:gamma1-nonincreasing}
Given a sequence of $\Omega$-valued random variables
$X_1,X_2,\ldots$ from distribution $F$,
let $\X(n):=\X(n-1) \cup \{X_n\}$ for $n=0,1,2,\ldots$ with $\X(0):= \emptyset$.
Then $\G_1\left(\X(n),N\right)$ is non-increasing in $n$ in the sense that
 $\G_1\left(\X(n+1),N\right) \subseteq \G_1\left(\X(n),N\right)$.
\end{lemma}
\noindent \textbf{Proof:}
Given a particular type of proximity map $N$ and
a data set $\X(n)=\bigl\{ X_1,X_2,\ldots,X_n \bigr\}$,
by Lemma \ref{lem:gamma1-intersect},
$\G_1\left(\X(n),N\right)= \bigcap_{i=1}^{n}\G_1\left(X_i,N\right)$
and by definition, $\X(n+1)=\X(n) \cup \{X_{n+1}\}$.
So,
\begin{multline*}
\G_1\left(\X(n+1),N\right)= \bigcap_{i=1}^{n+1}\G_1\left(X_i,N\right)=\\
\left[ \bigcap_{i=1}^{n}\G_1\left(X_i,N\right) \right]\bigcap \G_1\left(X_{n+1},N\right)=
\G_1\left(\X(n),N\right)\cap \G_1\left(X_{n+1},N\right) \subseteq \G_1\left(\X(n),N\right).
\end{multline*}

Thus we have shown that $\G_1\left(\X(n),N\right)$ is non-increasing in $n$; i.e.,
$\G_1\left(\X(n+1),N\right) \subseteq \G_1\left(\X(n),N\right)$.
$\blacksquare$
\begin{remark}
By monotone sequential continuity from above (\cite{billingsley:1995}),
the sequence $\bigl\{ \G_1\left(\X(n),N\right) \bigr\}_{n=1}^{\infty}$ has a limit
\begin{eqnarray}
\label{eqn:Limit-of-G1}
G_1&:=&\bigcap_{j=1}^{\infty} \G_1\left(\X(j),N\right)=
\lim_{m \rightarrow \infty}\bigcap_{j=1}^{m} \G_1\left(\X(j),N\right)=
\lim_{m \rightarrow \infty}\G_1\left(\X(m),N\right)\\
&=&\lim_{m \rightarrow \infty} \bigcap_{i=1}^m \G_1\left(X_i,N\right)=
\bigcap_{i=1}^{\infty} \G_1\left(X_i,N\right).\nonumber\;\;
\square
\end{eqnarray}
\end{remark}

\begin{theorem}
\label{thm:G1-as-conv1}
Given a sequence of random variables
$X_1,X_2,\ldots$ which are identically distributed as $F$ on $\Omega$,
let $\X(n):=\X(n-1) \cup \{X_n\}$ with $\X(0):= \emptyset$.
Then $\G_1\left(\X(n),N\right) \downarrow \RS(N)$, as $n \rightarrow \infty$ a.s.
in the sense that $\G_1\left(\X(n+1),N\right) \subseteq \G_1\left(\X(n),N\right)$ and
$\mu(\G_1\left(\X(n),N\right)\setminus \RS(N)) \downarrow 0$ a.s.
\end{theorem}
\noindent \textbf{Proof:}
By Lemma \ref{lem:gamma1-nonincreasing}, and by monotone
sequential continuity from above,
$\{\G_1\left(\X(n),N\right)\}_{n=1}^{\infty}$ has a limit, namely,
$G_1$ in Equation \eqref{eqn:Limit-of-G1}.
We claim that $G_1=\RS(N)$ a.s.

Suppose $\RS(N) \subsetneq \Omega$, since if $\RS(N) = \Omega$ then
$N(x)= \Omega$ for all $x \in \Omega$, so $\G_1\left(\X(n),N\right)=\Omega$ for
all $x$ hence the result would follow trivially.
Since $\RS(N) \subseteq \G_1\left(\X(n),N\right)$ for all $n$, $\RS(N) \subseteq G_1$.
From Proposition 5.6. in \cite{karr:1992},
$Y_n \xrightarrow{a.s.} Y \text{ iff } \forall\, \ve >0, \;\;
\lim_{n\rightarrow \infty}P\left(\sup_{k \ge n}|Y_k-Y| >\ve\right)=0.$
Let $\ve > 0$. Then
$$P\left(\sup_{k \ge n}\mu\left(\bigcap_{i=1}^{k}\G_1\left(X_i,N\right)\setminus \RS(N)\right) \le \ve\right)=
P\left(\mu\left(\bigcap_{i=1}^{n}\G_1\left(X_i,N\right)\setminus \RS(N)\right) \le  \ve\right) \rightarrow 1$$
as $n \rightarrow \infty$, because if $\G_1\left(\X(n),N\right) \setminus \RS(N)$
had positive measure, then for each $y \in \G_1\left(\X(n),N\right) \setminus \RS(N)$,
$\Omega \setminus N(y)$ will contain data points from $\X(k)$
with positive probability for sufficiently large $k \ge n$.
So $y$ can not be in $\G_1\left(\X(n),N\right)$, which is a contradiction.
Hence the desired result follows.
$\blacksquare$

Note however that $\G_1\left(\X_n,N\right)$ is neither strictly decreasing
nor non-increasing provided that $\RS(N) \not= \Omega$ for all
$\X_n$, because we might have $\G_1\left(\X_n,N\right) \subsetneq
\G_1\left(\X_m,N\right)$ for some $m > n$.
Nevertheless, the following two results hold.

\begin{proposition}
\label{prop:stoch-order-G1}
Suppose $\Omega \setminus \RS(N)$ has positive measure.
For positive integers $m > n$,
let $\X_n\}$ and $\X_m$ be two samples from $F$ on $\Omega$.
Then $\mu(\G_1\left(\X_m,N\right)) \le^{ST}\mu(\G_1\left(\X_n,N\right))$.
\end{proposition}
\noindent \textbf{Proof:}
Recall that for $X \sim F$ and $Y \sim G$, $X \le^{ST} Y$
if $F(x) \ge G(x)$ for all $x$ with strict inequality
holding for at least one $x$.

Let $m>n$ and $\X_n$ and $\X_m$ be two samples from $F$.  Then
$\mu(\G_1\left(\X_m,N\right))$ is more often smaller than $\mu(\G_1\left(\X_n,N\right))$.
Hence $P[\mu(\G_1\left(\X_m,N\right)) \le \mu(\G_1\left(\X_n,N\right))] \ge 1/2$ which only
shows stochastic precedence (\cite{boland:2004}).

Now, let $\displaystyle t \in \left(\mu(\RS(N)),\mu(\Omega) \right)$,
then $\displaystyle \mu(\G_1\left(\X_m,N\right))\le t$
happens more often than $\displaystyle \mu(\G_1\left(\X_n,N\right))\le t$,
hence $\displaystyle  P\left(\mu(\G_1\left(\X_m,N\right)) \le t\right) \ge
P\left(\mu(\G_1\left(\X_n,N\right)) \le t \right)$;
that is, $F_m(t) \ge F_n(t)$, where $F_i(\cdot)$ is the
distribution function for $\mu(\G_1\left(\X_i,N\right))$ for $i=m,n$.
For $t < \mu(\RS(N))$ or $t> \mu(\Omega)$, $F_i(t)=0$ for $i=m,n$.
Letting $N_n:=|\X_n \setminus \RS(N)|$ and $N_m:=|\X_m \setminus \RS(N)|$,
then $P(N_n \not=N_m)>0$ since $\mu(\Omega\setminus \RS(N))>0$.
In fact, $P(N_m>N_n) \ge 1/2$.
But if $F_m(t) = F_n(t)$ for all $t$
were the case, then $P(N_n=N_m)=1$ would hold, which is a contradiction.
$\blacksquare$

\begin{theorem}
\label{thm:AG1-converge-in-p}
Let $\{\X_n\}_{n=1}^{\infty}$ be a
sequence of samples of size $n$ from distribution $F$
with support on $\Omega$.
Then $\G_1\left(\X_n,N\right) \stackrel{p}{\longrightarrow} \RS(N)$ in
the sense that $\mu(\G_1\left(\X_n,N\right) \setminus \RS(N))
\stackrel{p}{\longrightarrow} 0$ as $n \rightarrow \infty$.
\end{theorem}
\noindent \textbf{Proof:}
Given a sequence $\{\X_n\}_{n=1}^{\infty}$ as in the theorem.
By Proposition \ref{prop:RSsubsetG1},
$\RS(N) \subseteq \G_1\left(\X_n,N\right) $ for each $n$.
If $\G_1\left(\X_n,N\right) \setminus \RS(N)$ has zero measure as $n \rightarrow
\infty$ then result follows trivially.
Otherwise, if $\G_1\left(\X_n,N\right) \setminus \RS(N)$ had positive measure in the limit,
for each $ y \in \lim_{n \rightarrow \infty}\G_1\left(\X_n,N\right) \setminus \RS(N)$,
$\Omega \setminus N(y)$ would have positive measure with
positive probability, then $\X_n \cap [\Omega \setminus N(y)] \not=\emptyset$
with positive probability for sufficiently large $n$,
then $y \notin \G_1\left(\X_n,N\right)$, which is a contradiction.
$\blacksquare$

Since $\G_1\left(\X_n,N\right)=\bigcap_{i=1}^n\,\G_1\left(X_i,N\right)$ for
a given realization of the data set $\X_n$,
first we describe the region $\G_1(X,N)$ for $X \in \X_n$,
and then describe the region $\G_1\left(\X_n,N\right)$.


\begin{theorem}
If the superset region for any type of proximity map $N$ has positive measure
(i.e., $\mu(\RS(N)) > 0$), then $P(\g_n(N)=1)
\rightarrow 1$ as $n \rightarrow \infty$.
\end{theorem}
\noindent {\bf Proof:}
Notice that if there is at least one data point in $\RS(N)$ then
$\g_n(N)=1$, because any point $x \in \RS(N)$
will have $N(x) = \Omega$,
so $P(\text{there is at least 1 point in $\RS(N)$}) \le P(\g_n(N)=1)$.
Now, $\displaystyle P(\text{there is at least 1 point in $\RS(N)$})=1-P(\mathcal X_n \cap
\RS=\emptyset)=1-\left(\frac{\mu(\Omega)-\mu(\RS(N))}{\mu(\Omega)}\right)^n$, which goes to 1 as
$n \rightarrow \infty$.
Hence $P(\g_n(N)=1) \rightarrow
1$ as $n \rightarrow \infty$. $\blacksquare$

The \emph{relative arc density} of a digraph $D=(\V,\A)$
of order $|\V| = n$,
denoted as $\rho(D)$,
is defined as
$\displaystyle \rho(D) = \frac{|\A|}{n(n-1)}$
where $|\cdot|$ denotes the cardinality of sets (\cite{janson:2000}).
Thus $\rho(D)$ represents the ratio of the number of arcs
in the digraph $D$ to the number of arcs in the complete symmetric
digraph of order $n$, which is $n(n-1)$.

\begin{theorem}
If support of the joint distribution of $\X_n$ is subset of the  superset region
for any type of proximity map,
then relative arc density $\rho(\mathcal X_n)=1$ a.s.
\end{theorem}
\noindent {\bf Proof:}
Suppose the support of $\X_n$ is a subset of the superset region,
then the corresponding digraph is complete with $n(n-1)$ arcs.
Hence the relative arc density is 1 with probability 1.
$\blacksquare$

\section{Transformations Preserving Uniformity on Triangles in $\R^2$}
\label{sec:transformations}
The proximity regions and hence the corresponding PCDs
are based on the Delaunay tessellation of $\Y_m$,
which partitions $C_H(\Y_m)$.
So, suppose the set $\X_n$ is a set of iid uniform random variables on
the convex hull of $\Y_m$; i.e., a random sample from $\U(\C_H\left(\Y_m\right))$.
In particular, conditional on $|\X_n \cap T_i|>0$
being fixed, $\X_n \cap T_i$ will
also be a set of iid uniform random variables on $T_i$ for $i \in \{1,2,\ldots,J\}$,
where $T_i$ is the $i^{th}$ Delaunay cell
and $J$ is the total number of Delaunay cells.
Reducing the cell (triangle in $\R^2$) $T_i$ as much as possible
while preserving uniformity and the
probabilities related to PCDs will simplify the notation and calculations.
For simplicity we consider $\R^2$ only.

Let $\Y_3 =\{\y_1,\y_2,\y_3\} \subset \R^2$ be three non-collinear points
and $\TY$ be the triangle (including the interior)
with vertices $\y_1,\y_2,\y_3$. Let $X_i \stackrel{iid}{\sim} \; \UT$,
the uniform distribution on $\TY$, for $i=1,2,\ldots, n$.
The probability density function (pdf) of $\UT$ is
$$f(u)=\frac{1}{A(\TY)}\I(u \in \TY),$$
where $A(\cdot)$ is the area functional.

The triangle $\TY$ can be carried into the first quadrant by a
composition of transformations (scaling, translation, rotation, and reflection)
in such a way that the largest edge has unit length and lies on the $x$-axis, and the
$x$-coordinate of the vertex nonadjacent to largest edge is less than $1/2$.
We call the resultant triangle the \emph{basic triangle} and denote it as $T_b$
where $T_b=\bigl((0,0),(1,0),(c_1,c_2)\bigr)$ with $0< c_1 \le 1/2$,
and $c_2 > 0$ and $(1-c_1)^2+c_2^2 \le 1$.
The transformation from any triangle to $T_b$ is denoted by $\phi_b$.
See \cite{ECarXivPCDGeo:2009} for a detailed description of $\phi_b$.
Notice that $\TY$ is transformed into $T_b$, then $\TY$ is similar to $T_b$
and $\phi_b\left(\TY\right)=T_b$.
Thus the random variables $X_i \stackrel{iid}{\sim}\UT$ transformed
along with $\TY$ in the described fashion by $\phi_b$ satisfy $\phi_b(X_i) \stackrel{iid}{\sim}\U(T_b)$.
So, without loss of generality, we can assume $\TY$ to be $T_b$ for uniform data.
The functional form of $T_b$ is
$$T_b=\left\{ (x,y) \in \R^2 : y \ge 0;\; y \le (c_2\,x)/c_1;\; y \le c_2\,(1-x)/(1-c_1) \right\}.$$

There are other transformations that preserve uniformity
of the random variable, but not similarity of the triangles.
We only describe the transformation that maps $T_b$ to the
standard equilateral triangle,
$T_e=T\left((0,0),(1,0),\left(1/2,\sqrt{3}/2\right)\right)$
for exploiting the symmetry in calculations using $T_e$.

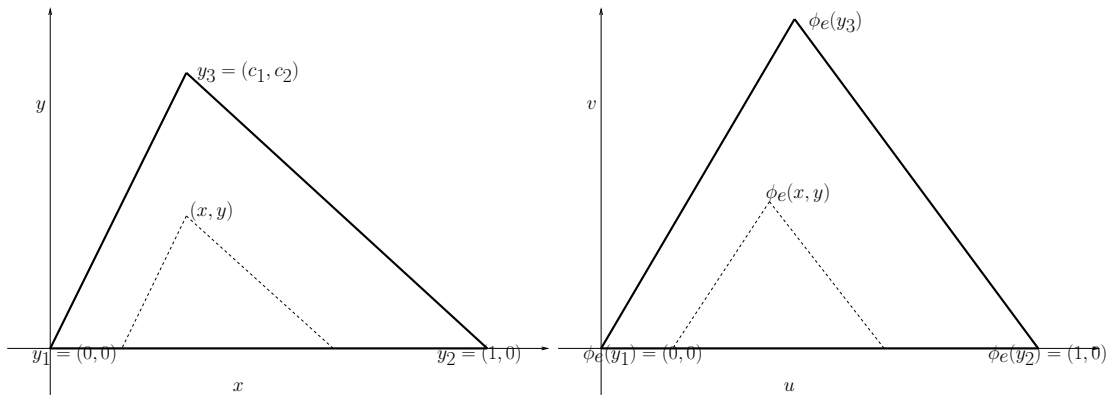
\begin{figure} [ht]
\centering
\scalebox{.3}{\input{Basic_Tri.pstex_t}}
\scalebox{.3}{\input{Equal_Tri_phi.pstex_t}}
\caption{The description of $\phi_e(x,y)$ for $(x,y) \in T_b$ (left) and
the equilateral triangle $\phi_e(T_b)=T_e$ (right).}
\label{fig:transform-equal}
\end{figure}

Let $\phi_e:\,(x,y) \rightarrow (u,v)$, where
$\displaystyle u(x,y)=x+\frac{1-2\,c_1}{\sqrt{3}}\,y$ and
$\displaystyle v(x,y)=\frac{\sqrt{3}}{2\,c_2}\,y$.
Then $\y_1$ is mapped to $(0,0)$, $\y_2$ is mapped to $(1,0)$,
and $\y_3$ is mapped to $\left(1/2,\sqrt{3}/2\right)$.
See also Figure \ref{fig:transform-equal}.

Note that the inverse transformation is
$\phi_e^{-1}(u,v)=\bigl(x(u,v),y(u,v)\bigr)$ where
$\displaystyle x(u,v)=u-\frac{(1-2\,c_1)}{\sqrt{3}}\,v$ and
$\displaystyle y(u,v)=\frac{2\,c_2}{\sqrt{3}}\,u$.
Then the Jacobian is given by
\begin{eqnarray*}
J(x,y)&=&\left|\begin{array}{cc}
\frac{\partial{x}}{\partial{u}} & \frac{\partial{x}}{\partial{v}}
\vspace{.1 in}\\
\frac{\partial{y}}{\partial{u}} & \frac{\partial{y}}{\partial{v}}
 \end{array}\right|=\left|\begin{array}{cc}
1 & \frac{2\,c_1-1}{\sqrt{3}}
\vspace{.1 in}\\
0 & \frac{2\,c_2}{\sqrt{3}}
 \end{array}\right|=\frac{2\,c_2}{\sqrt{3}}.
\end{eqnarray*}
So $\displaystyle f_{U,V}(u,v)=f_{X,Y}(\phi_e^{-1}(u,v))\,|J| =
\frac{4}{\sqrt{3}}\, \I\bigl((u,v)\in T_e\bigr)$.
Hence uniformity is preserved.

\begin{remark}
\label{rem:geo-dependence-NS-NAS}
The probabilities for uniform data in $\TY$
involve the ratio of the event region to $A(\TY)$.
If such ratios is not preserved under $\phi_e$,
then the probability content for $N$ depends on the geometry of $\TY$.
In particular, the probability content for uniform data for $\NS$ and $\NAS$
depend on the geometry of the triangle, hence is not geometry invariant (\cite{ECarXivPCDGeo:2009}).
For example, $P(X \in \NS(Y))$ and $P(X \in \NAS(Y,M))$ depends
on $(c_1,c_2)$, hence one has to do the computations
for all of (uncountably many) of these triangles.
Hence we do not investigate $\NS$ and $\NAS$ further in this article.
$\square$
\end{remark}

\section{Sample Proximity Maps in $\R^d$}
\label{sec:example-PCDs}
Let $\Y_m=\left \{\y_1,\y_2,\ldots,\y_m \right\}$ be $m$ points
in general position in $\R^d$ and $\T_i$ be the $i^{th}$
Delaunay cell for $i=1,2,\ldots,J$, where $J$ is the number
of Delaunay cells.
Let $\X_n$ be a random sample from a distribution $F$ in $\R^d$
with support $\mS(F) \subseteq \C_H\left(\Y_m\right)$.

In particular, for illustrative purposes, we focus on $\R^2$,
where a Delaunay tessellation is a triangulation,
provided that no more than three points in $\Y_m$ are cocircular.
Furthermore, for simplicity,
let $\Y_3=\{\y_1,\y_2,\y_3\}$ be three non-collinear points
in $\R^2$ and $\TY=T(\y_1,\y_2,\y_3)$ be the triangle
with vertices $\Y_3$.
Let $\X_n$ be a random sample from $F$ with
support $\mS(F) \subseteq \TY$.
In this section, we will describe two families of triangular proximity regions.

\subsection{Proportional-Edge Proximity Maps}
\label{sec:prop-edge}
The first type of triangular proximity map we introduce is the
proportional-edge proximity map.
For this proximity map,
the asymptotic distribution of domination number and
the relative density of the corresponding
PCD will have mathematical tractability.

\begin{figure} [ht]
\centering
\scalebox{.35}{\input{Nofnu2.pstex_t}}
\caption{Construction of proportional edge proximity region,
$\NPE^2(x)$ (shaded region).}
\label{fig:ProxMapDef1}
\end{figure}
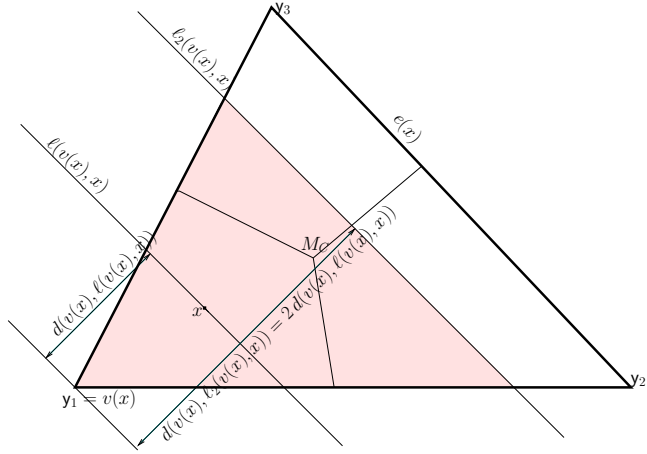

For $r \in [1,\infty]$,
define $\NPE^r(\cdot,M):=N(\cdot,M;r,\Y_3)$ to be the
\emph{proportional-edge proximity map}
with $M$-vertex regions as follows
(see also Figure  \ref{fig:ProxMapDef1} with $M=M_C$ and $r=2$).
For $x \in \TY \setminus \Y_3$, let $v(x) \in \Y_3$ be the
vertex whose region contains $x$; i.e., $x \in R_M(v(x))$.
If $x$ falls on the boundary of two $M$-vertex regions,
we assign $v(x)$ arbitrarily.
Let $e(x)$ be the edge of $\TY$ opposite $v(x)$.
Let $\ell(v(x),x)$ be the line parallel to $e(x)$ through $x$.
Let $d(v(x),\ell(v(x),x))$ be the Euclidean (perpendicular) distance
from $v(x)$ to $\ell(v(x),x)$.
For $r \in [1,\infty)$, let $\ell_r(v(x),x)$ be the line parallel to $e(x)$
such that
\begin{gather*}
d(v(x),\ell_r(v(x),x)) = r\,d(v(x),\ell(v(x),x)) \text{ and }
d(\ell(v(x),x),\ell_r(v(x),x)) < d(v(x),\ell_r(v(x),x)).
\end{gather*}
Let $T_r(x)$ be
the triangle similar to
and with the same orientation as $\TY$
having $v(x)$ as a vertex
and $\ell_r(v(x),x)$ as the opposite edge.
Then the {\emph proportional-edge proximity region}
$\NPE^r(x,M)$ is defined to be $T_r(x) \cap \TY$.
Notice that $\ell(v(x),x)$ divides the edges of $T_r(x)$
(other than $\ell_r(v(x),x)$) proportionally with the factor $r$.
Hence the name \emph{proportional edge proximity region}.

Notice that $r \ge 1$ implies $x \in \NPE^r(x,M)$.
Furthermore,
$\lim_{r \rightarrow \infty} \NPE^r(x,M) = \TY$ for all $x \in \TY \setminus \Y_3$,
so we define $\NPE^{\infty}(x,M) = \TY$ for all such $x$.
For $x \in \Y_3$, we define $\NPE^r(x,M) = \{x\}$ for all $r \in [1,\infty]$.

Notice that $X_i \stackrel{iid}{\sim} F$,
with the additional assumption
that the non-degenerate two-dimensional
probability density function $f$ exists
with support $\mS(F) \subseteq \TY$,
implies that the special case in the construction
of $\NPE^r$ ---
$X$ falls on the boundary of two vertex regions ---
occurs with probability zero.
Note that for such an $F$, $\NPE^r(X)$ is a triangle a.s.

The functional form of $\NPE^r(x,M)$ for $x=(x_0,y_0) \in T_b$ is given in \cite{ECarXivPCDGeo:2009}.
Of particular interest is $\NPE^r(x,M)$ with any $M$ and $r \in \{\sqrt{2},\,3/2,2\}$.
For $r=\sqrt{2}$, $\ell(v(x),x)$ divides $T_{\sqrt{2}}(x)$ into
two regions of equal area, hence $\NPE^{\sqrt{2}}(x,M)$ is also
referred to as \emph{double-area proximity region}.
For $r=2$, $\ell(v(x),x)$ divides the edges of $T_2(x)$ ---other than
$\ell_r(v(x),x)$ --- into two segments of equal length, hence
$\NPE^2(x,M)$ is also referred to as \emph{double-edge proximity region}.
For $r < 3/2$, $\RS(\NPE^r,M_C)=\emptyset$,
and for $r > 3/2$, $\RS(\NPE^r,M_C)$
has positive area; for $r=3/2$, $\RS(\NPE^r,M_C)=\{M_C\}$.
Therefore, $r=3/2$ is the threshold for the superset region
to be nonempty.
Furthermore, $r=3/2$ will be the value at which the asymptotic distribution of the
domination number of the PCD based on $\NPE^r(\cdot,M_C)$ is nondegenerate
(\cite{ceyhan:dom-num-NPE-SPL}).

\begin{figure}
\begin{center}
\psfrag{A}{\scriptsize{$\y_1$}}
\psfrag{B}{\scriptsize{$\y_2$}}
\psfrag{C}{\scriptsize{$\y_3$}}
\psfrag{CC}{\scriptsize{$M_{CC}$}}
\psfrag{x}{\scriptsize{$x$}}
\psfrag{M1}{\scriptsize{$M_3$}}
\psfrag{M2}{\scriptsize{$M_1$}}
\psfrag{M3}{\scriptsize{$M_2$}}
\epsfig{figure=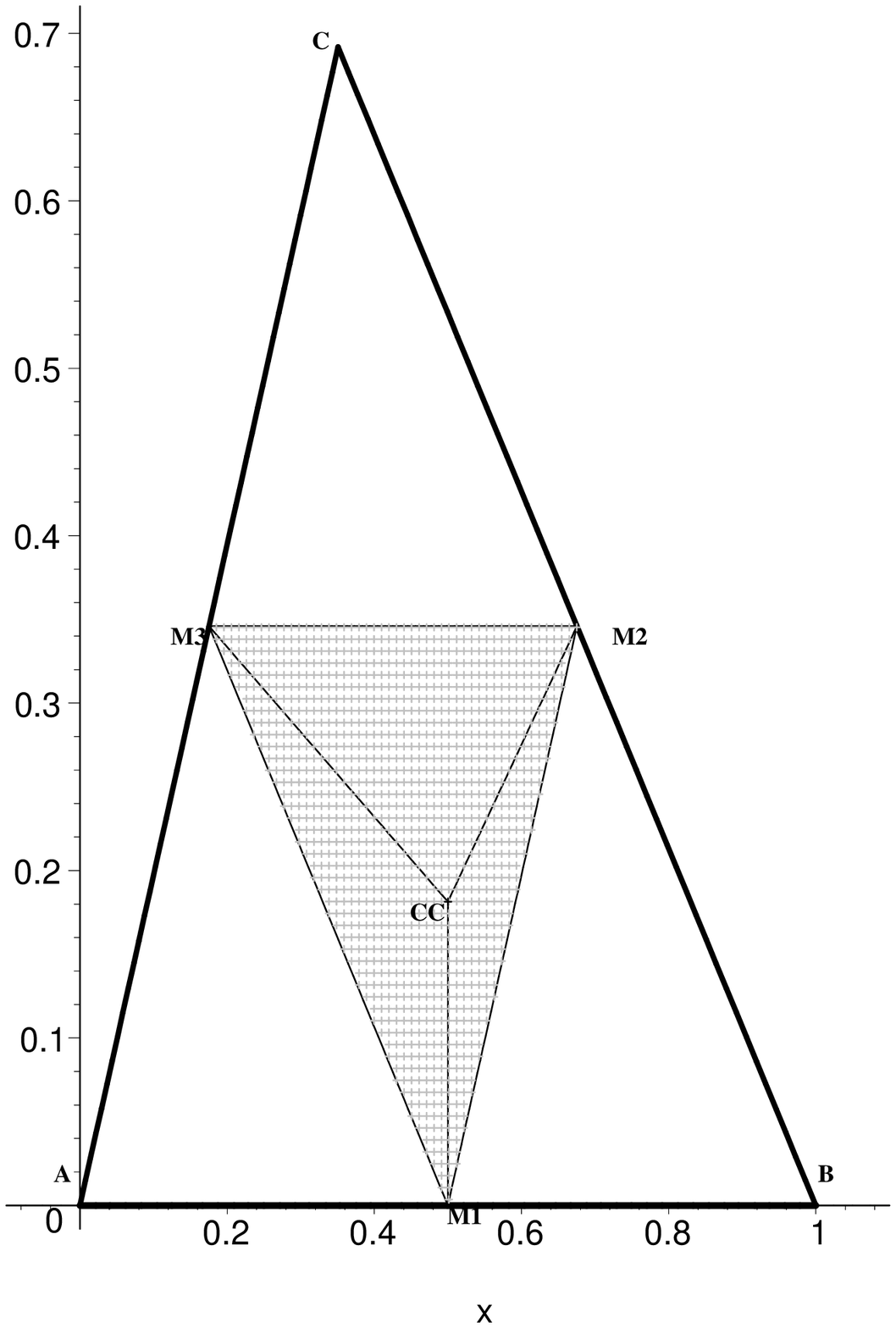,  height=140pt , width=200pt}
\psfrag{IC}{\scriptsize{$M_I$}}
\psfrag{P1}{\scriptsize{$P^{IC}_3$}}
\psfrag{P2}{\scriptsize{$P^{IC}_1$}}
\psfrag{P3}{\scriptsize{$P^{IC}_2$}}
\psfrag{Q2}{}
\psfrag{Q3}{}
\epsfig{figure=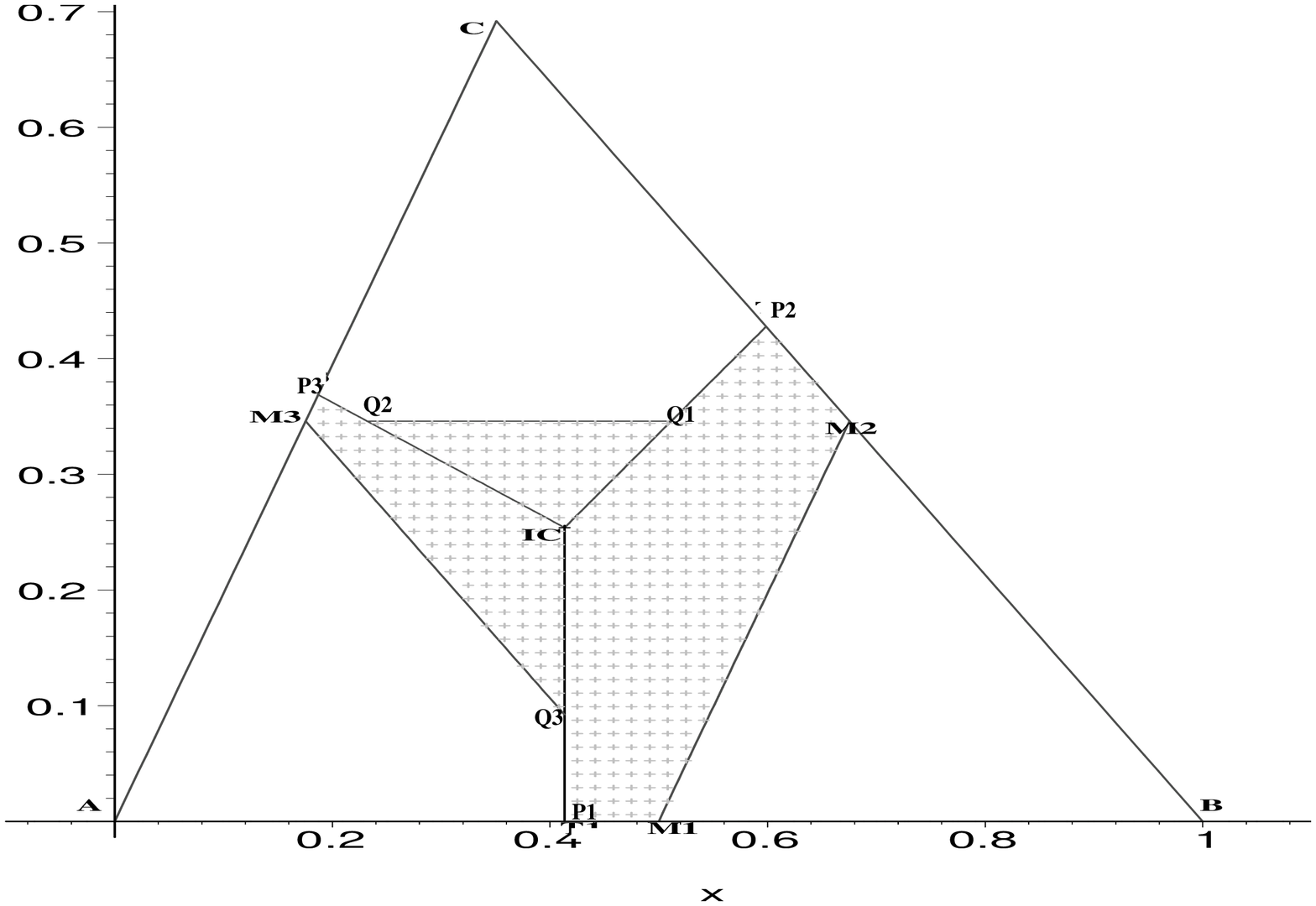, height=140pt , width=200pt}
\end{center}
\caption{
\label{fig:superset-NDE}
Superset region $\RS\left(\NPE^2,M_{CC}\right)$ in an acute
triangle (left), superset region, $\RS^{\perp}(\NPE^2,M_I)$ (right)
}
\end{figure}


Let $\RS^{\perp}\left( \NPE^r,M \right)$ be the superset region for
$\NPE^r$ based on $M$-vertex regions with orthogonal projections.
Then the superset region with the incenter $\RS^{\perp}\left(\NPE^2,M_I\right)$
is as in Figure \ref{fig:superset-NDE} (right).
Let $M_i$ be the midpoint of edge $e_i$ for $i=1,2,3$.
Then $T(M_1,M_2,M_3) \subseteq \RS^{\perp}(\NPE^2,M_I)$ for all $\TY$
with equality holding when $\TY$ is an equilateral triangle.

For $\NPE^2(\cdot,M_C)$ constructed using the median lines $\RS\left(\NPE^2,M_C\right)=T(M_1,M_2,M_3)$,
and for $\NPE^2(\cdot,M_C)$ constructed by the orthogonal projections,
$\RS^{\perp} \left(\NPE^2,M_C\right) \supseteq T(M_1,M_2,M_3)$ with equality holding when $\TY$ is an equilateral triangle.

\begin{figure}
\begin{center}
\psfrag{A}{\scriptsize{$\y_1$}}
\psfrag{B}{\scriptsize{$\y_2$}}
\psfrag{C}{\scriptsize{$\y_3$}}
\psfrag{CC}{\scriptsize{$M_{CC}$}}
\psfrag{IC}{\scriptsize{$M_{I}$}}
\psfrag{x}{\scriptsize{$x$}}
\psfrag{D}{}
\psfrag{E}{}
\psfrag{F}{}
\psfrag{S}{}
\psfrag{Ab}{}
\psfrag{Ac}{}
\psfrag{Ba}{}
\psfrag{Bc}{}
\psfrag{Ca}{}
\psfrag{Cb}{}
\epsfig{figure=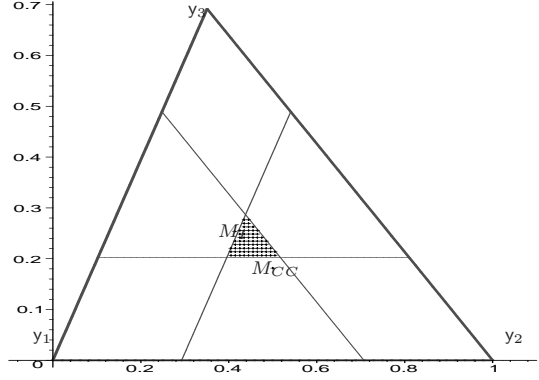, height=140pt , width=200pt}
\end{center}
\caption{The triangle $\mathscr T^{r=\sqrt{2}}$.
}
\label{fig:Tr-RS-NDA-CC}
\end{figure}

In $\TY$, drawing the lines $\xi_i(r,x)$ such that
$d(\y_i,e_i)=r\,d(\xi_i(r,x),\y_i)$ for $j\in \{1,2,3\}$  yields a triangle,
$\Tr$, for $r<3/2$ .
See Figure \ref{fig:Tr-RS-NDA-CC} for $\Tr$ with $r=\sqrt{2}$.

The functional form of $\Tr$ in $T_b$ is
{\small
\begin{align}
\label{eqn:T^r-def}
& \Tr = \left \{(x,y) \in T_b: y \ge \frac{c_2\,(r-1)}{r};\;
y \le \frac{c_2\,(1-r\,x)}{r\,(1-c_1)};\; y \le \frac{c_2\,(r\,(x-1)+1)}{r\,c_1} \right\}\\
&=T\Biggl( \left(\frac{(r-1)\,(1+c_1)}{r},\frac{c_2\,(r-1)}{r}\right),
\left(\frac{2-r+c_1\,(r-1)}{r},\frac{c_2\,(r-1)}{r}\right),
\left(\frac{c_1\,(2-r)+r-1}{r},\frac{c_2\,(r-2)}{r}\right) \Biggr) \nonumber.
\end{align}
}
There is a crucial difference between $\Tr$ and $T(M_1,M_2,M_3)$:
$T(M_1,M_2,M_3) \subseteq \RS\left( \NPE^r,M \right)$ for all $M$ and $r \ge 2$,
but $(\Tr)^o$ and $\RS\left( \NPE^r,M \right)$ are disjoint for all $M$ and $r$.

So if $M \in (\Tr)^o$, then $\RS\left( \NPE^r,M \right)=\emptyset$;
if $M \in \partial(\Tr)$, then $\RS\left( \NPE^r,M \right)=\{M\}$;
and if $M \not\in \Tr$, then $\RS\left( \NPE^r,M \right)$ has positive area.
The triangle $\Tr$ defined above plays a crucial role in the analysis of the distribution
of the domination number of the PCD based on proportional-edge proximity maps.
The superset region $\RS\left( \NPE^r,M \right)$ will be important for
both the domination number and the relative density of the
corresponding PCDs.

The functional forms of the superset region, $\RS\left( \NPE^r,M \right)$,
and $T(M_1,M_2,M_3)$ in $T_b$ are provided in \cite{ECarXivPCDGeo:2009}.

$\NPE^r(\cdot,M)$ is geometry invariant if
$M$-vertex regions are constructed with $M \in \TY^o$ by
using the extensions of the line segments joining $\y$ to $M$
for all $\y \in \Y_3$.
But when the vertex regions are constructed by orthogonal projections,
$\NPE^r(\cdot,M)$ is not geometry invariant (\cite{ECarXivPCDGeo:2009}),
hence such vertex regions are not considered here.

\subsubsection{Extension of $\NPE^r$ to Higher Dimensions}
\label{sec:NYr-higher-D}
The extension to $\R^d$ for $d > 2$ is straightforward.
The extension with $M=M_C$ is given her,
but the extension for general $M$ is similar.
Let $\Y_{d+1} = \{\y_1,\y_2,\ldots,\y_{d+1}\}$ be $d+1$ points
that do not lie on the same $d-1$-dimensional hyperplane.
Denote the simplex formed by these $d+1$ points as $\mathfrak S (\Y_{d+1})$.
A simplex is the simplest polytope in $\R^d$
having $d+1$ vertices, $d\,(d+1)/2$ edges and $d+1$ faces of dimension $(d-1)$.
For $r \in [1,\infty]$, define the proximity map as follows.
Given a point $x$ in $\mathfrak S (\Y_{d+1})$,
let $v := \argmin_{\y \in \Y_{d+1}} \mbox{V}(Q_{\y}(x))$
where $Q_{\y}(x)$ is the polytope with vertices
being the $d\,(d+1)/2$ midpoints of the edges,
the vertex $v$ and $x$ and $\mbox{V}(\cdot)$ is the
$d$-dimensional volume functional.
That is, the vertex region for vertex $v$ is the polytope with vertices
given by $v$ and the midpoints of the edges.
Let $v(x)$ be the vertex in whose region $x$ falls.
If $x$ falls on the boundary of two vertex regions,
$v(x)$ is assigned arbitrarily.
Let $\varphi(x)$ be the face opposite to vertex $v(x)$,
and $\Upsilon(v(x),x)$ be the hyperplane parallel to $\varphi(x)$ which contains $x$.
Let $d(v(x),\Upsilon(v(x),x))$ be the (perpendicular)
Euclidean distance from $v(x)$ to $\Upsilon(v(x),x)$.
For $r \in [1,\infty)$, let $\Upsilon_r(v(x),x)$ be
the hyperplane parallel to $\varphi(x)$
such that
\begin{gather*}
d(v(x),\Upsilon_r(v(x),x))=r\,d(v(x),\Upsilon(v(x),x))
\text{ and }
d(\Upsilon(v(x),x),\Upsilon_r(v(x),x))< d(v(x),\Upsilon_r(v(x),x)).
\end{gather*}
Let $\mathfrak S_r(x)$ be the polytope similar to and with the same
orientation as $\mathfrak S$ having $v(x)$
as a vertex and $\Upsilon_r(v(x),x)$ as the opposite face.
Then the proximity region
$\NPE^r(x,M_C):=\mathfrak S_r(x) \cap \mathfrak S(\Y_{d+1})$.
Notice that $r \ge 1$ implies $x \in \NPE^r(x,M_C)$.

\subsubsection{$\G_1$-Regions for Proportional-Edge Proximity Maps}
\label{sec:Gamma1-NYr}
\begin{figure} [ht]
\centering
\scalebox{.35}{\input{Gammaofnu2.pstex_t}}
\caption{Construction of the $\G_1$-region, $\G_1\left(x,\NPE^{r=2},M_C\right)$ (shaded region). }
\label{fig:ProxMapDef2}
\end{figure}
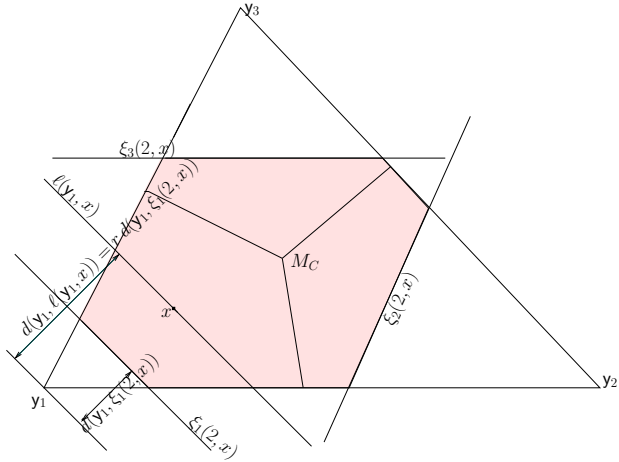

For $\NPE^r(\cdot,M)$, the $\G_1$-region is constructed as follows;
see also Figure \ref{fig:ProxMapDef2}.
Let $\xi_i(r,x)$ be the line parallel to $e_i$
such that $\xi_i(r,x)\cap \TY \not=\emptyset$ and
$r\,d(\y_i,\xi_i(r,x))=d(\y_i,\ell(\y_i,x))$ for $i \in \{1,2,3\}$.
Then
$$\G_1\left(x,\NPE^r,M\right)=\bigcup_{i=1}^3 \bigl[ \G_1\left(x,\NPE^r,M\right)\cap R_M(\y_i)\bigr]$$
where
$$\G_1\left(x,\NPE^r,M\right)\cap R_M(\y_i)=\{z \in R_M(\y_i):
d(\y_i,\ell(\y_i,z)) \ge d(\y_i,\xi_i(r,x)\} \text{ for } i \in \{1,2,3\}. $$
Notice that $r \ge 1$ implies $x \in \G_1\left(x,\NPE^r,M\right)$.
Furthermore, $\lim_{r \rightarrow \infty} \G_1\left(x,\NPE^r,M\right) = \TY$
for all $x \in \TY \setminus \Y_3$
and so we define $\G_1\left(x,\NPE^{r=\infty},M\right) = \TY$ for all such $x$.
For $x \in \Y_3$, $\G_1\left(x,\NPE^r,M\right)= \{x\}$ for all $r \in [1,\infty]$.

The functional form of $\G_1\left(x=(x_0,y_0),\NPE^r,M\right)$ in the basic triangle $T_b$ is given by
$$\G_1\left(x=(x_0,y_0),\NPE^r,M\right)=\bigcup_{i=1}^3 \left[ \G_1\left(x=(x_0,y_0),\NPE^r,M\right) \cap R_M(y_i) \right]$$
where
\begin{align*}
\G_1\left(x=(x_0,y_0),\NPE^r,M\right) \cap R_M(\y_1)&=\left \{(x,y) \in R_M(\y_1):
y \ge \frac{y_0}{r}-\frac{c_2\,(r\,x-x_0)}{(1-c_1)\,r} \right \},\\
\G_1\left(x=(x_0,y_0),\NPE^r,M\right) \cap R_M(\y_2)&=\left \{(x,y) \in R_M(\y_1):
y \ge \frac{y_0}{r}-\frac{c_2\,(r\,(x-1)+1-x_0}{c_1\,r} \right \},\\
\G_1\left(x=(x_0,y_0),\NPE^r,M\right) \cap R_M(\y_3)&=\left \{(x,y) \in R_M(\y_1):
y \le \frac{y_0-c_2\,(1-r)}{r} \right\}.
\end{align*}
Notice that $\G_1\left(x,\NPE^r,M_C\right)$ is a convex hexagon for all $r \ge 2$
and $x \in \TY \setminus \Y_3$, (since for such an $x$,
$\G_1\left(x,\NPE^r,M_C\right)$ is bounded by $\xi_i(r,x)$
and $e_i$ for all $i \in \{1,2,3\}$, see also Figure \ref{fig:ProxMapDef2})
else it is either a convex hexagon or a non-convex polygon depending on
the location of $x$ and the value of $r$.

Furthermore, in $\R^d$ with $d>2$ let $\zeta_i(x)$ be the hyperplane such that
$\zeta_i(x) \cap \mathfrak S\left(\Y_m\right) \not=\emptyset$ and
$r\,d(\y_i,\zeta_i(x))=d(\y_i,\eta(\y_i,x))$ for $i=1,2,\ldots,(d+1)$.
Then $\G_1\left(x,\NPE^r\right)\cap R_{CM}(\y_i)=\{z \in R_{CM}(\y_i):
d(\y_i,\eta(\y_i,z)) \ge d(\y_i,\zeta_i(x)\}$, for $i=1,2,3$.  Hence
$\G_1\left(x,\NPE^r\right)=\bigcup_{i=1}^{d+1} (\G_1\left(x,\NPE^r\right)\cap R_{CM}(\y_i))$.
Notice that $r \ge 1$ implies $x \in \NPE^r(x)$ and $x \in
\G_1\left(x,\NPE^r\right)$.

So far, we have described the $\G_1$-region for a point in $x \in \TY$.
For a set $\X_n$ of size $n$ in $\TY$, the region $\G_1\left(\X_n,\NPE^r,M\right)$
can be exactly described by the edge extrema.

\begin{definition}
The (closest) \emph{edge extrema} of a set $B$ in $\TY$
are the points closest to the edges of $\TY$, denoted $x_{e_i}$ for $i \in \{1,2,3\}$;
that is, $x_{e_i} \in \arginf_{x \in B}d(x,e_i)$.
\end{definition}
Note that if $B=\X_n$ is a random sample of size $n$ from $F$ then the edge extrema,
denoted $X_{e_i}(n)$, are random variables.

\begin{figure}
\begin{center}
\psfrag{A}{\scriptsize{$\y_1$}}
\psfrag{B}{\scriptsize{$\y_2$}}
\psfrag{C}{\scriptsize{$\y_3$}}
\psfrag{(x1,y1)}{\scriptsize{$x$}}
\psfrag{Q3}{\scriptsize{$B_3$}}
\psfrag{Q1}{\scriptsize{$B_1$}}
\psfrag{Q2}{\scriptsize{$B_2$}}
\psfrag{x}{}
\psfrag{CC}{\scriptsize{$M_{CC}$}}
\psfrag{G(A)}{}
\psfrag{G(B)}{}
\psfrag{G(C)}{\scriptsize{$\G_1\left(x,\NPE^2,M_{CC}\right) \cap R_{CC}(\y_3)$}}
\epsfig{figure=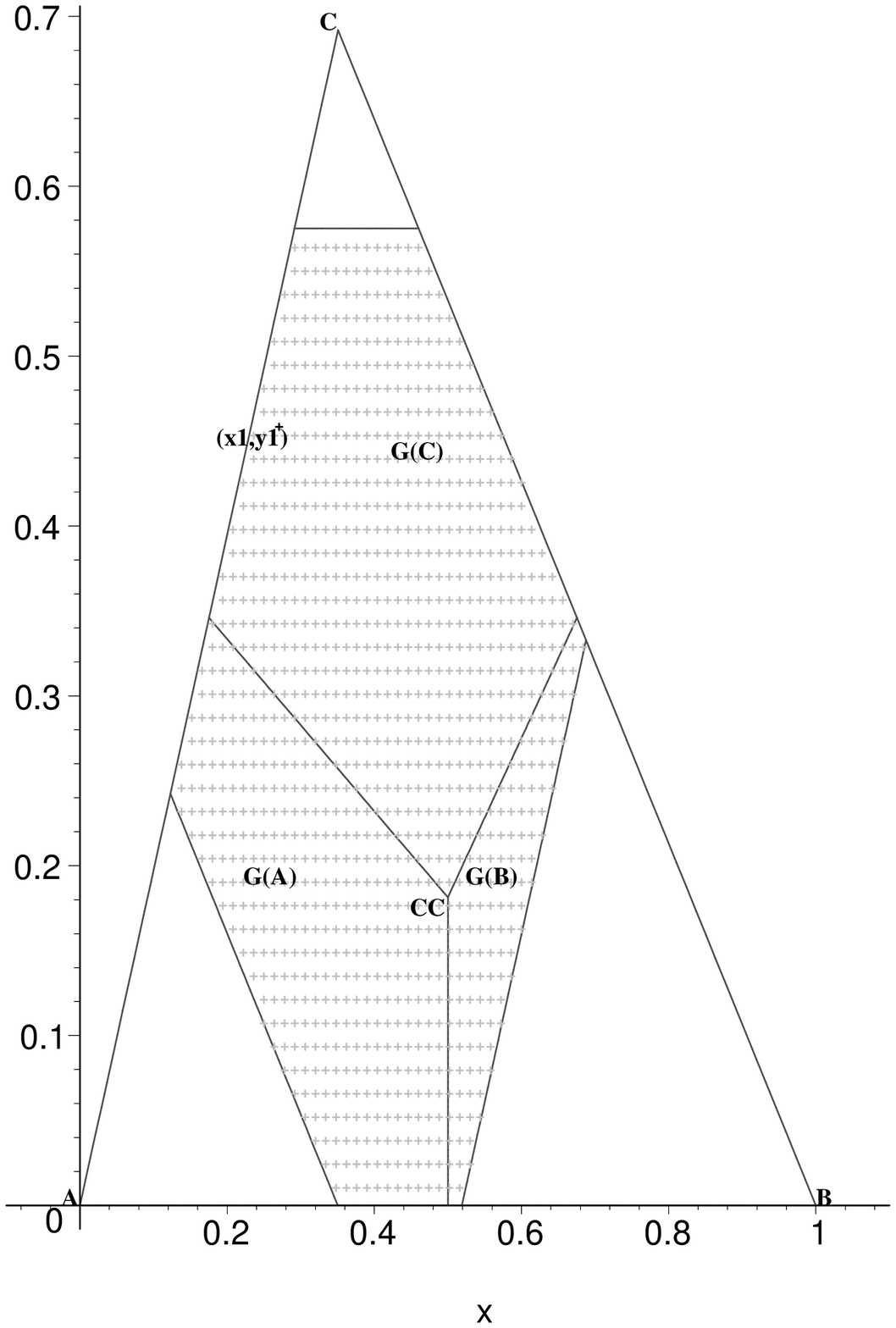, height=140pt , width=200pt}
\psfrag{IC}{\scriptsize{$M_I$}}
\psfrag{G(A)}{}
\psfrag{G(B)}{}
\psfrag{G(C)}{\scriptsize{$\G_1\left(x,\NPE^2,M_I\right) \cap R^{\perp}_{IC}(\y_3)$}}
\epsfig{figure=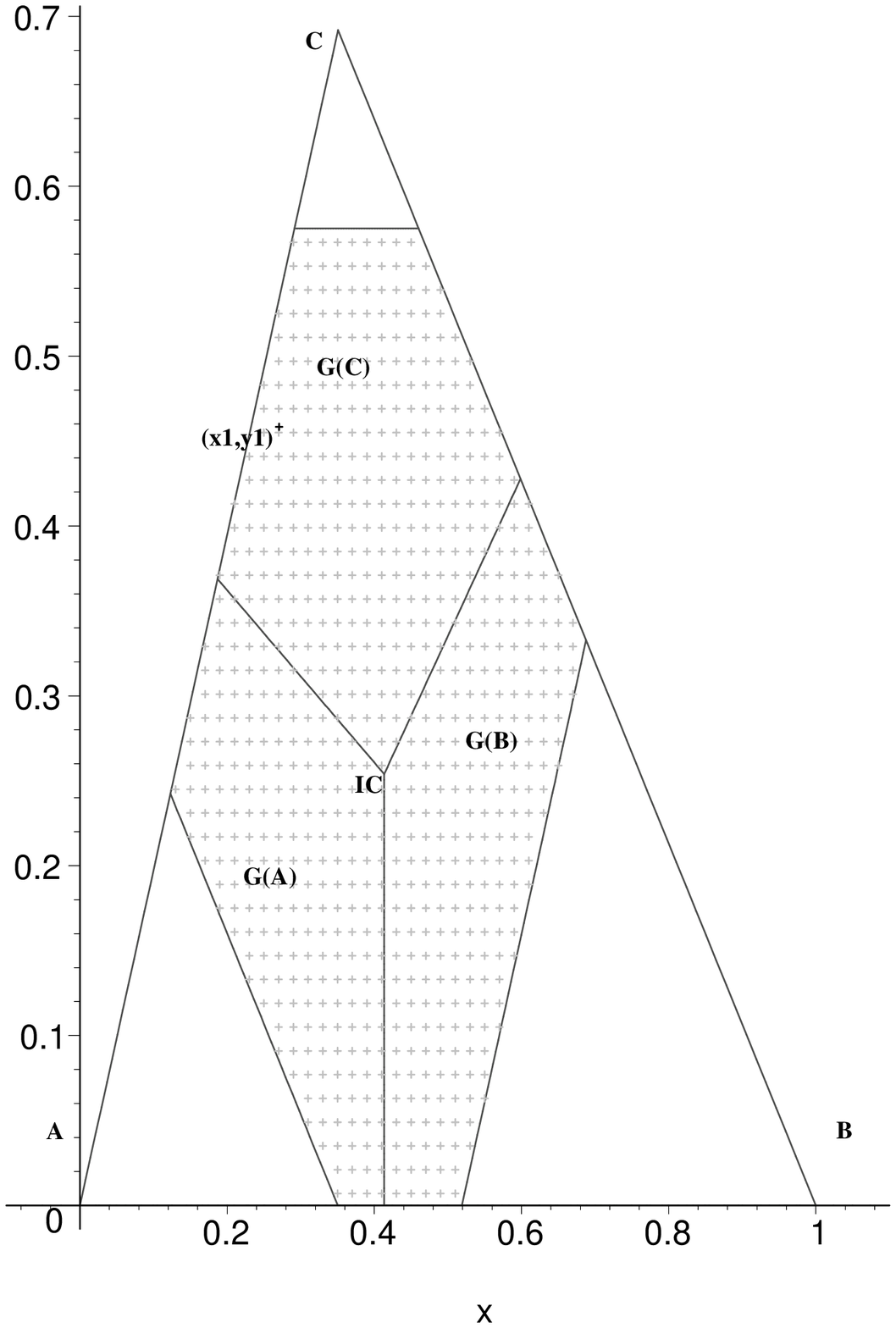, height=140pt , width=200pt}
\end{center}
\caption{
\label{fig:G1-NDE-1}
$\G_1\left(x,\NPE^2,M_{CC}\right)$ with $x \in R_{CC}(\y_3)$ (left), $x \in
R^{\perp}_{IC}(\y_3)$ (right).}
\end{figure}

\begin{proposition}
\label{prop:SMA-edge-extrema}
Let $B$ be any set of $n$ distinct points in $\TY$
and $x_{e_i} \in \arginf_{x \in B}d(x,e_i)$.
For proportional-edge proximity maps with $M$-vertex regions,
$\G_1\left(B,\NPE^r,M\right)=\bigcap_{i=1}^3\,\G_1\left(x_{e_i},\NPE^r,M\right)$.
\end{proposition}
\noindent \textbf{Proof:}
Given $B = \{x_1,x_2,\ldots,x_n\}$ in $\TY$.
Note that
$$\G_1\left(B,\NPE^r,M\right) \cap R_M(\y_i)=
\left[\bigcap_{i=1}^n\,\G_1\left(x_i,\NPE^r,M\right)\right]\cap R_M(\y_i),$$
and if $d(\y_i,\ell(\y_i,x)) \le d(\y_i,\ell(\y_i,x'))$
then $\NPE^r(x,M) \subseteq \NPE^r(x',M)$ for all $x,x' \in R_M(y_i)$.
Further, by definition  $x_{e_i} \in \argmax_{x\in B}d(\y_i,\xi_i(r,x))$,
so
$$\G_1\left(B,\NPE^r,M\right) \cap R_M(\y_i)=
\G_1\left(x_{e_i},\NPE^r,M\right)\cap R_M(\y_i) \text{ for }i \in \{1,2,3\}. $$
Furthermore, $\G_1\left(B,\NPE^r,M\right) =
\bigcup_{i=1}^3 \bigl[ \G_1\left(x_{e_i},\NPE^r,M\right) \cap R_M(\y_i) \bigr]$,
and
$$\G_1\left(x_{e_i},\NPE^r,M\right) \cap R_M(\y_i)=
\left[\bigcap_{j=1}^3  \G_1\left(x_{e_j},\NPE^r,M\right)\right] \cap R_M(\y_i)
\text{ for }i \in \{1,2,3\}.$$
Combining these two results, we obtain
$\G_1\left(B,\NPE^r,M\right)=\bigcap_{j=1}^3\,\G_1\left(x_{e_j},\NPE^r,M\right)$. $\blacksquare$

From the above proposition, we see that the $\G_1$-region for $B$ as in the proposition
can also be written as the union of three regions of the form
$$\G_1\left(B,\NPE^r,M\right)\cap R_M(\y_i)=\{z \in R_M(\y_i):\;
d(\y_i, \ell(\y_i,z)) \ge d(\y_i,\xi_i(r,x_{e_i}))\} \text{ for } i \in \{1,2,3\}.$$

\begin{corollary}
Let $\X_n$ be a random sample from a continuous distribution $F$ on $\TY$.
For proportional-edge proximity maps with $M$-vertex regions,
$\eta_n\left(\NPE^r\right) \le 3$ with equality
holding with positive probability for $n \ge 3$.
\end{corollary}
\noindent \textbf{Proof:}
From Proposition \ref{prop:SMA-edge-extrema},
$\eta_n\left(\NPE^r\right) \le 3$.
Furthermore, $X_e(n)$ is unique for each edge $e$ a.s. since $F$ is continuous,
and there are three distinct edge extrema with positive probability.
Hence $P(\eta_n\left(\NPE^r\right)= 3)>0$ for $n \ge 3$. $\blacksquare$

Then $\G_1\left(\X_n,\NPE^r,M\right)=\bigcap_{i=1}^3\G_1\left( X_{e_i},\NPE^r\right)$,
where $e_i$ is the edge opposite vertex $\y_i$, for $i=1,2,3$.
So $\G_1\left(\X_n,\NPE^r,M\right)\cap R_M(\y_i)=
\{z \in R_M(\y_i):\; d(\y_i, \ell(\y_i,z)) \ge d(\y_i,\xi_i(r,x_{e_i}))\}$,
for $i=1,2,3$.


Note that $P(\eta_n\left(\NPE^r\right) = 3) \rightarrow 1$
as $n \rightarrow \infty$ for $\X_n$ a random sample from $\UT$,
since edge extrema are distinct with probability 1
as $n \rightarrow \infty$ as shown in the following theorem.
\begin{theorem}
\label{thm:distinct-edge-ext}
Let $\X_n$ be a random sample from $\UT$ and
let $E_{c,3}(n)$ be the event that (closest) edge extrema are distinct.
Then $P(E_{c,3}(n)) \rightarrow 1$ as $n \rightarrow \infty$.
\end{theorem}
\noindent \textbf{Proof:}
Using the uniformity preserving transformation $\phi_e$
without loss of generality, one can assume $\X_n$ is a random sample from $\U(T_e)$.
Observe also that the edge extrema in $T_b$ are mapped into the edge extrema in $T_e$.
Note that the probability of edge extrema all being equal to each other is
$P(X_{e_1}(n)=X_{e_2}(n)=X_{e_3}(n))=\I(n=1)$.
Let $E_{c,2}(n)$ be the event that there are only two distinct (closest) edge extrema.
Then for $n>1$,
$$P(E_{c,2}(n))=P(X_{e_1}(n)=X_{e_2}(n))+P(X_{e_1}(n)=X_{e_3}(n))+P(X_{e_2}(n)=X_{e_3}(n))$$
since the intersection of events $X_{e_1}(n)=X_{e_2}(n))$,
$X_{e_1}(n)=X_{e_3}(n)$, and $X_{e_2}(n)=X_{e_3}(n)$ is equivalent to
the event $X_{e_1}(n)=X_{e_2}(n)=X_{e_3}(n)$.
Notice also that $P(E_{c,2}(n=2))=1$.
So, for $n>2$, there are two or three distinct edge extrema with probability 1;
i.e., $P(E_{c,3}(n))+P(E_{c,2}(n))=1$ for $n > 2$.

We will show that $P(E_{c,2}(n)) \rightarrow 0$ as $n \rightarrow \infty$,
which will imply the desired result.

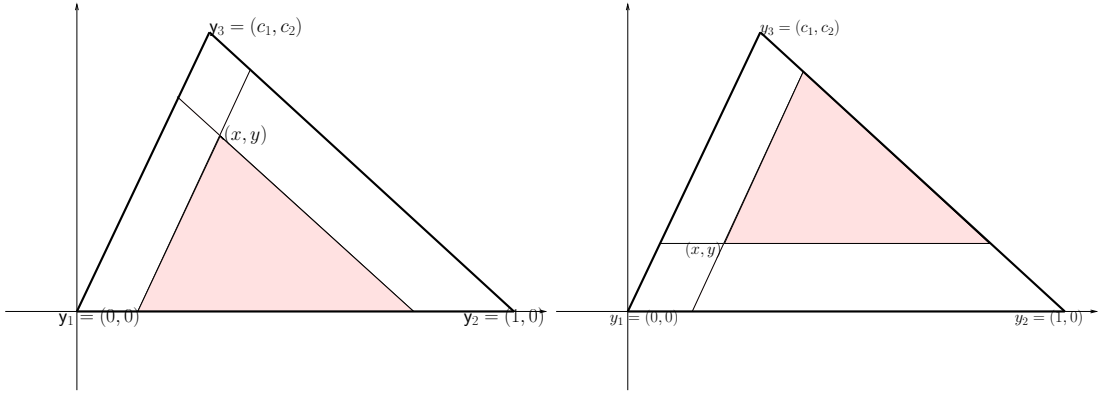
\begin{figure} [ht]
\centering
\scalebox{.3}{\input{two_extrema.pstex_t}}
\scalebox{.3}{\input{two_ext_F.pstex_t}}
\caption{
\label{fig:edge-extrema}
The figure for $X_{e_1}(n)=X_{e_2}(n)=(x,y)$
(left) and $X_{e_2}(n)=X_{e_3}(n)=(x,y)$ (right).}
\end{figure}

\begin{figure}
\centering
\psfrag{A}{\scriptsize{$\y_1$}}
\psfrag{B}{\scriptsize{$\y_2$}}
\psfrag{C}{\scriptsize{$\y_3$}}
\psfrag{(x,y)}{\scriptsize{$(x,y)$}}
\psfrag{x}{}
\epsfig{figure=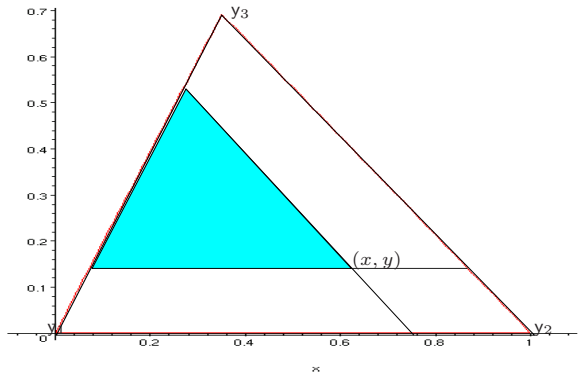, height=140pt , width=240pt}
\caption{
\label{fig:two-extrema}
The figure for $X_{e_1}(n)=X_{e_3}(n)=(x,y)$.}
\end{figure}

First consider $P(X_{e_1}(n)=X_{e_2}(n))$.
The event $X_{e_1}(n)=X_{e_2}(n)=X_e=(X,Y)$ is equivalent to the event that
$\X_n \subset \{U \in T_e: \; d(\y_1,\ell(\y_1,U)) \le d(\y_1,\ell(\y_1,X_e)),\;
d(\y_2,\ell(\y_2,U)) \le d(\y_2,\ell(\y_2,X_e))\}$.
For example, if given $X_{e_1}(n)=X_{e_2}(n)=(x,y)$ the
remaining $n-1$ points will lie in the shaded region in Figure \ref{fig:edge-extrema} (left).
For other pairs of edge extrema, see Figures \ref{fig:edge-extrema} (right) and \ref{fig:two-extrema}.

The pdf of such $X_e=(X,Y)$ is
$f(x,y)=n\,\left(4/\sqrt{3}\right)\,\left(y^2/\sqrt{3}\right)^{n-1}$.
Let $\ve > 0$, by Markov's inequality,
$P\left(\sqrt{3}/2-Y > \ve\right) \le \E \left[ \sqrt{3}/2-Y \right]/\ve$.
But,
\begin{eqnarray*}
\E\left[ \sqrt{3}/2-Y \right]&=&\int_0^{1/2}\int_0^{\sqrt{3}\,x}
\left(\sqrt{3}/2-y\right)\,n\,\left(y^2/\sqrt{3}\right)^{n-1}\,\left(4/\sqrt{3}\right) dy\,dx\\
& &+\int_{1/2}^1\int_0^{\sqrt{3}\,(1-x)} \left(\sqrt{3}/2-y\right)\,n\,
\left(y^2/\sqrt{3}\right)^{n-1}\,\left(4/\sqrt{3}\right) dy\,dx\\
&=&4\,\left(\sqrt{3}/4\right)^n\frac{1}{n\,(4\,n^2-1)},
\end{eqnarray*}
which converges to 0 as $n \rightarrow \infty$.
So if $X_{e_1}(n)=X_{e_2}(n)$ were the case,
then geometric locus of this point goes to $\y_3$.
That is, for each $\ve >0$,
$P\left(\sqrt{3}/2-Y > \ve\right) \rightarrow 0$ as $n \rightarrow \infty$.
Hence $P(X_{e_1}(n)=X_{e_2}(n)\not= \y_3)\rightarrow 0$ as $n \rightarrow \infty$.
Furthermore, $P(X_{e_1}(n)=X_{e_2}(n)= \y_3) \le P(X_{e_2}(n)\in e_2)=0$ for all $n \ge 1$.
So $P(X_{e_1}(n)=X_{e_2}(n))\rightarrow 0$ as $n \rightarrow \infty$.

Likewise, by symmetry, it follows that
$\lim_{n\rightarrow \infty}P(X_{e_1}(n)=X_{e_3}(n))=\lim_{n\rightarrow \infty}P(X_{e_2}(n)=X_{e_3}(n))=0$.
Hence $P(E_{c,2}(n))\rightarrow 0$ as $n \rightarrow \infty$.
Thus $P(E_{c,3}(n)) \rightarrow 1$ as $n \rightarrow \infty$. $\blacksquare$

The above theorem implies that the asymptotic distribution
of $\eta_n\left(\NPE^r\right)$ is degenerate with
$P(\eta_n\left(\NPE^r\right)=3) \rightarrow 1$ as $n\rightarrow \infty$.
But for finite $n$, $\eta_n\left(\NPE^r\right)$ for $X_i \stackrel{iid}{\sim}\UT$
has the following non-degenerate distribution.
\[
\eta_n\left(\NPE^r\right) ~=~ \left\{
        \begin{array} {ll}
           2 ~~~ & \mbox{wp $\pi_2(n)$} \\[1ex]
           3 & \mbox{wp $\pi_3(n)=1-\pi_2(n)$,}
        \end{array}
          \right.
\]
where $\pi_2(n) \in (0,1)$ is the probability
of edge extrema for any two distinct edges being concurrent.

\begin{remark}
If $\X_n$ is a random sample from $F$ such that
$\mS(F) \cap \{x \in \TY:\;d(x,e_i) \le \ve_1 \}$ has
positive measure and $\mS(F) \cap B(\y_i,\ve_2)=\emptyset$
for some $\ve_1,\ve_2 >0$, then $P(E_{c,3}(n)) \rightarrow 1$
as $n \rightarrow \infty$ follows trivially.
However, the case that $F$ has positive density around the vertices $\Y_3$
requires more work to prove as shown below. $\square$
\end{remark}

\begin{theorem}
\label{thm:gen-F-edge-ext}
Let $\X_n$ be a random sample from $F$
such that $B(\y_i,\ve) \subseteq \mS(F)$ for some
$\ve >0$ and for each $i=1,2,3$, then $P(E_{c,3}(n)) \rightarrow
1$ as $n \rightarrow \infty$.
\end{theorem}
\noindent \textbf{Proof:}
Using the transformation $\phi_e:\,(x,y)\rightarrow (u,v)$,
above, we can without loss of generality assume
$\X_n$ is a random sample from $F$ with support $\mS(F)\subseteq T_e$.
After $\phi_e$ is applied, suppose $F$ becomes $F_e$,
then $\phi_e(\X_n)$ becomes a random sample from $F_e$
such that $B(\phi_e(\y_i),\ve_e) \subseteq \mS(F_e)$ for some
$\ve_e >0$ and for each $i=1,2,3$.
First consider $P(X_{e_1}=X_{e_3})$.
Given $X_{e_1}=X_{e_3}=(x,y)$ the remaining $n-1$ points will lie in the
shaded region in Figure \ref{fig:edge-extrema} (right).
$X_{e_1}=X_{e_3}=(X,Y)$ is equivalent to the event that
$$\X_n \subset S_R(X,Y):=\{(U,W) \in T_e: \; \ell(\y_1,(U,W)) \le \ell(\y_1,(X,Y)),\;\ell(\y_3,(U,W)) \le \ell(\y_3,(X,Y))\}.$$
The pdf of such $(X,Y)$ is
$f(x,y)=n\,G(x,y)^{n-1}\,f(x,y)$ where $G(u,v):=P_F(X \in S_R(u,v))$.
Note that $P(X_{e_2}=X_{e_3}=\y_1)=0$ for all $n$, and
$X_{e_2}=X_{e_3}\not=\y_1$ is equivalent to $d((X,Y),\y_1)>0$.
Let $\ve > 0$, by Markov's inequality,
$P(d((X,Y),\y_1) > \ve) \le \E\left[ d((X,Y),\y_1) \right]/\ve=
\E\left[ \sqrt{X^2+Y^2} \right]/\ve$.
Switching to the polar coordinates as $X=R\,\cos \theta$ and $Y=R\,\sin \theta$,
we get $\sqrt{X^2+Y^2}=R$.
But, $\displaystyle \E[R]=\int_0^{\ve}\int_0^{\pi/3}n\,r\,G(r,\theta)^{n-1}\,f(r,\theta)r\,dr\,d\theta$.
Integrand is critical at $r=0$, since for $r>0$
it converges to zero as $n \rightarrow \infty$.
So we use the Taylor series expansion around $r=0$ as
\begin{align*}
f(r,\theta)&=f(0,\theta)+\frac{\partial{f(0,\theta)}}{\partial r}r+O\left(r^2\right),\\
G(r,\theta)&=G(0,\theta)+\frac{\partial{G(0,\theta)}}{\partial
r}r+O\left(r^2\right)=1+\frac{\partial{G(0,\theta)}}{\partial
r}r+O\left(r^2\right).
\end{align*}
Note that $\displaystyle \frac{\partial{G(0,\theta)}}{\partial r}<0$,
since area of $S_R(u,v)$ decreases as $r$ increases for fixed $\theta$.
So let $r=w/n$, then
\begin{eqnarray*}
\E[R]&\sim& \int_0^{n\,\ve}\int_0^{\pi/3} n\,
\frac{w}{n}\,\left(1+\frac{\partial{G(0,\theta)}}{\partial r}\frac{w}{n}+
O\left(n^{-2}\right)\right)^{n-1}\,
\left(f(0,\theta)+\frac{\partial{f(0,\theta)}}{\partial r}\frac{w}{n}+
O\left(n^{-2}\right)\right)\frac{w}{n^2}\,dw\,d\theta\\
&=&\frac{1}{n^2}\int_0^{\infty}\int_0^{\pi/3} w^2\,\exp\left(
\frac{\partial{G(0,\theta)}}{\partial
r}w\right)\,f(0,\theta)w\,dw\,d\theta=O\left(n^{-2}\right).
\end{eqnarray*}
Hence $P(X_{e_2}=X_{e_3}\not= \y_1)\rightarrow 0$ as $n \rightarrow
\infty$. Then $P(X_{e_2}=X_{e_3})\rightarrow 0$ as $n \rightarrow
\infty$.

Likewise, it follows that $\lim_{n\rightarrow
\infty}P(X_{e_1}=X_{e_2})=\lim_{n\rightarrow
\infty}P(X_{e_1}=X_{e_3})=0$. Hence $P(E_{c,2}(n))\rightarrow 0$ as
$n \rightarrow \infty$. Thus $P(E_{c,3}(n)) \rightarrow 1$ as $n
\rightarrow \infty$. $\blacksquare$

Notice that Theorem \ref{thm:distinct-edge-ext} follows
as a corollary from Theorem \ref{thm:gen-F-edge-ext}.
For $r \ge 3/2$ and $M \in \R^2 \setminus \Y_3$, $\G_1\left(\X_n,\NPE^r,M\right) \not=  \emptyset$ a.s.,
 since $\RS\left( \NPE^r,M \right)\not= \emptyset$ and $\RS\left( \NPE^r,M \right) \subseteq \G_1\left(\X_n,\NPE^r,M\right)$.

Now, for $n>1$, let $X_{e_i}(n)=x_{e_i}=(u_i,w_i)$ be given
for $i \in \{1,2,3\},$ be the edge extrema in
a given realization of $\X_n$.
Then the functional form of $\G_1$-region in $T_b$ is given by
$$\G_1\left(\X_n,\NPE^r,M\right)=\bigcup_{i=1}^3 \bigl[\G_1\left(\X_n,\NPE^r,M\right) \cap R_M(\y_i)\bigr]$$
where
\begin{align*}
\G_1\left(\X_n,\NPE^r,M\right) \cap R_M(\y_1)&=\left \{(x,y) \in R_M(\y_1):
y \ge \frac{w_1}{r}-\frac{c_2\,(r\,x-u_1)}{(1-c_1)\,r} \right \},\\
\G_1\left(\X_n,\NPE^r,M\right) \cap R_M(\y_2)&=\left \{(x,y) \in R_M(\y_1):
y \ge \frac{w_2}{r}-\frac{c_2\,(r\,(x-1)+1-u_2}{c_1\,r} \right\},\\
\G_1\left(\X_n,\NPE^r,M\right) \cap R_M(\y_3)&=\left \{(x,y) \in R_M(\y_1):
y \le \frac{w_3-c_2\,(1-r)}{r} \right\}.
\end{align*}

\begin{figure}
\begin{center}
\psfrag{A}{\scriptsize{$\y_1$}}
\psfrag{B}{\scriptsize{$\y_2$}}
\psfrag{C}{\scriptsize{$\y_3$}}
\psfrag{(x1,y1)}{\scriptsize{$x_{e_1}$}}
\psfrag{(x2,y2)}{\scriptsize{$x_{e_2}$}}
\psfrag{(x3,y3)}{\scriptsize{$x_{e_3}$}}
\psfrag{Q3}{\scriptsize{$Q_3$}}
\psfrag{Q1}{\scriptsize{$Q_1$}}
\psfrag{Q2}{\scriptsize{$Q_2$}}
\psfrag{x}{}
\psfrag{CC}{\scriptsize{$M_{CC}$}}
\psfrag{G(A)}{}
\psfrag{G(B)}{}
\psfrag{G(C)}{\scriptsize{$\G_1\left(\X_n,\NPE^2,M\right) \cap R_{CC}(\y_3)$}}
\epsfig{figure=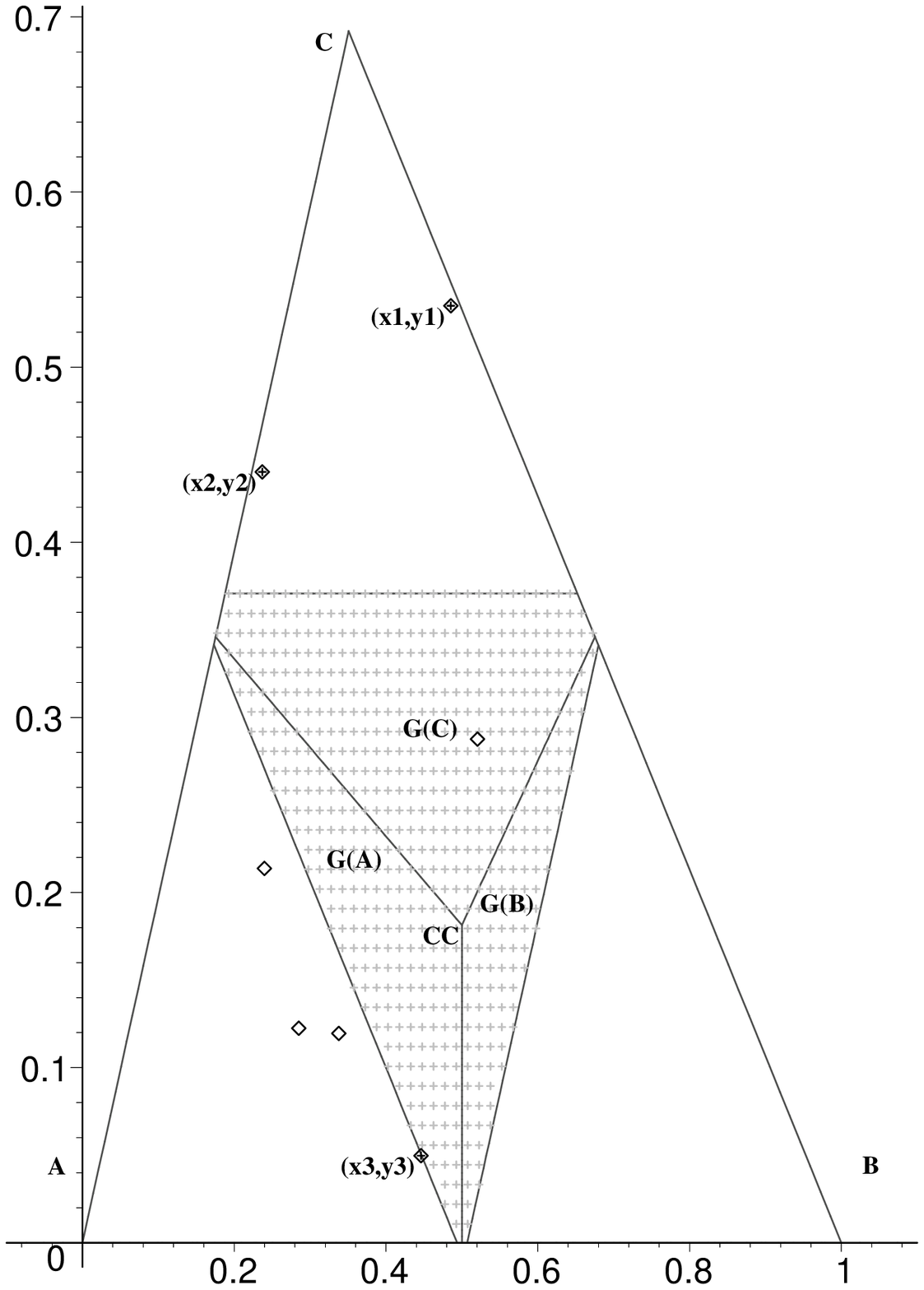, height=140pt , width=200pt}
\psfrag{IC}{\scriptsize{$M_I$}}
\psfrag{G(A)}{}
\psfrag{G(B)}{}
\psfrag{G(C)}{\scriptsize{$\G_1\left(x,\NPE^2\right) \cap R^{\perp}_{IC}(\y_3)$}}
\epsfig{figure=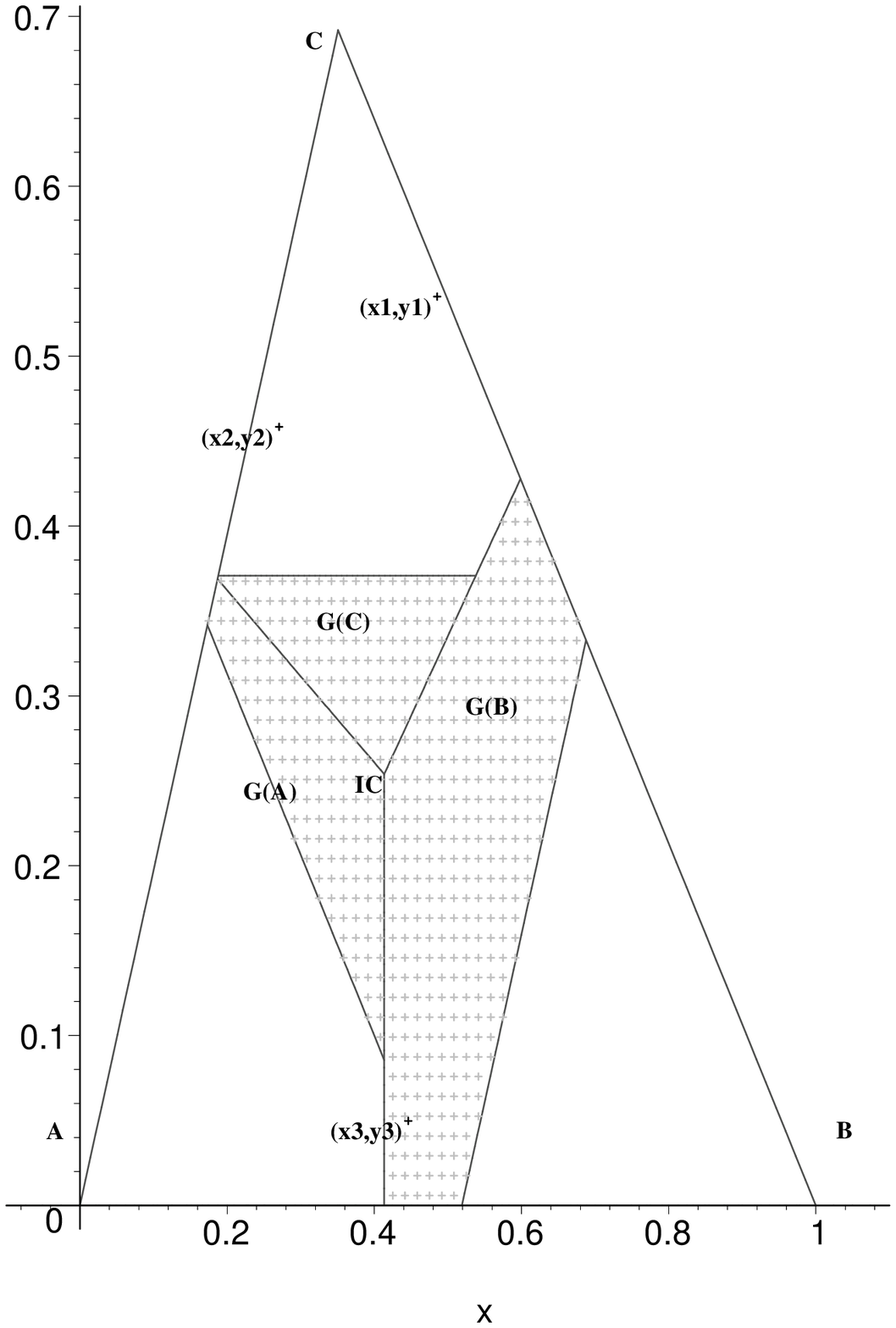, height=140pt , width=200pt}
\end{center}
\caption{$\G_1\left(\X_n,\NPE^2,M\right)$ for $M=M_{CC}$ (left) and $M=M_I$
(right) with three distinct edge extrema.} \label{fig:G1-NDE-n}
\end{figure}

See Figure \ref{fig:G1-NDE-n}  for $\G_1\left(\X_n,\NPE^2,M\right)$ with $n \ge
3$ where $M$-vertex regions for $M=M_{CC}$ and $M=M_I$ with
orthogonal projections are used.  Note that only the edge extrema
are shown in Figure \ref{fig:G1-NDE-n} (right).

Note that, for $\X_n$ a random sample from $\UT$,
 $P(\eta_n\left(\NPE^r\right) = 3) \rightarrow 1$ as $n \rightarrow \infty$,
since the edge extrema are distinct with probability 1 as
$n \rightarrow \infty$.
However, for $r<3/2$, the region $\G_1\left(x,\NPE^r,M\right)\cap R_M(\y_i)$
might be empty for some $i \in \{1,2,3\}$.
Furthermore, if $M \in (\Tr)^o$ (see Equation \eqref{eqn:T^r-def}
for $\Tr$)
with $r<3/2$, then $\G_1\left(\X_n,\NPE^r,M\right)$ will be empty with
probability 1 as $n \rightarrow \infty$.
In such a case, there is no $\G_1$-region to construct.
But the definition of the $\eta$-value still works in the sense that
$\G_1\left(X_n,\NPE^r,M\right)=\emptyset=\bigcap_{x \in S_M} \G_1\left(x,\NPE^r,M\right)$
(see Definition \ref{def:active-set} for $S_M$) and
$\G_1\left(x,\NPE^r,M\right) \not=\emptyset$
for all $x \in \X_n$ since $x \in \G_1\left(x,\NPE^r,M\right)$.
To determine whether the $\G_1$-region is empty or not,
it suffices to check the intersection of the
$\G_1$-regions of the edge extrema.
If $M \notin (\Tr)^o$, the $\G_1$-region
is guaranteed to be nonempty.

Note that $\eta_n \left(\NPE^{r_1}\right) \stackrel{d}{=}\eta_n \left(\NPE^{r_2}\right)$
for all $(r_1,\,r_2) \in [1,\infty)^2$,
where $\stackrel{d}{=}$ stands for ``equality in distribution".

\begin{figure}
\begin{center}
\psfrag{A}{\scriptsize{$\y_1$}}
\psfrag{B}{\scriptsize{$\y_2$}}
\psfrag{C}{\scriptsize{$\y_3$}}
\psfrag{(x1,y1)}{\scriptsize{$x_{e_1}$}}
\psfrag{(x2,y2)}{\scriptsize{$x_{e_2}$}}
\psfrag{(x3,y3)}{\scriptsize{$x_{e_3}$}}
\psfrag{Q3}{\scriptsize{$Q_3$}}
\psfrag{Q1}{\scriptsize{$Q_1$}}
\psfrag{Q2}{\scriptsize{$Q_2$}}
\psfrag{x}{}
\psfrag{CC}{\scriptsize{$M_{CC}$}}
\psfrag{G(A)}{}
\psfrag{G(B)}{}
\psfrag{G(C)}{\scriptsize{$\G_1\left(\X_n,\NPE^2,M\right) \cap R_{CC}(\y_3)$}}
\epsfig{figure=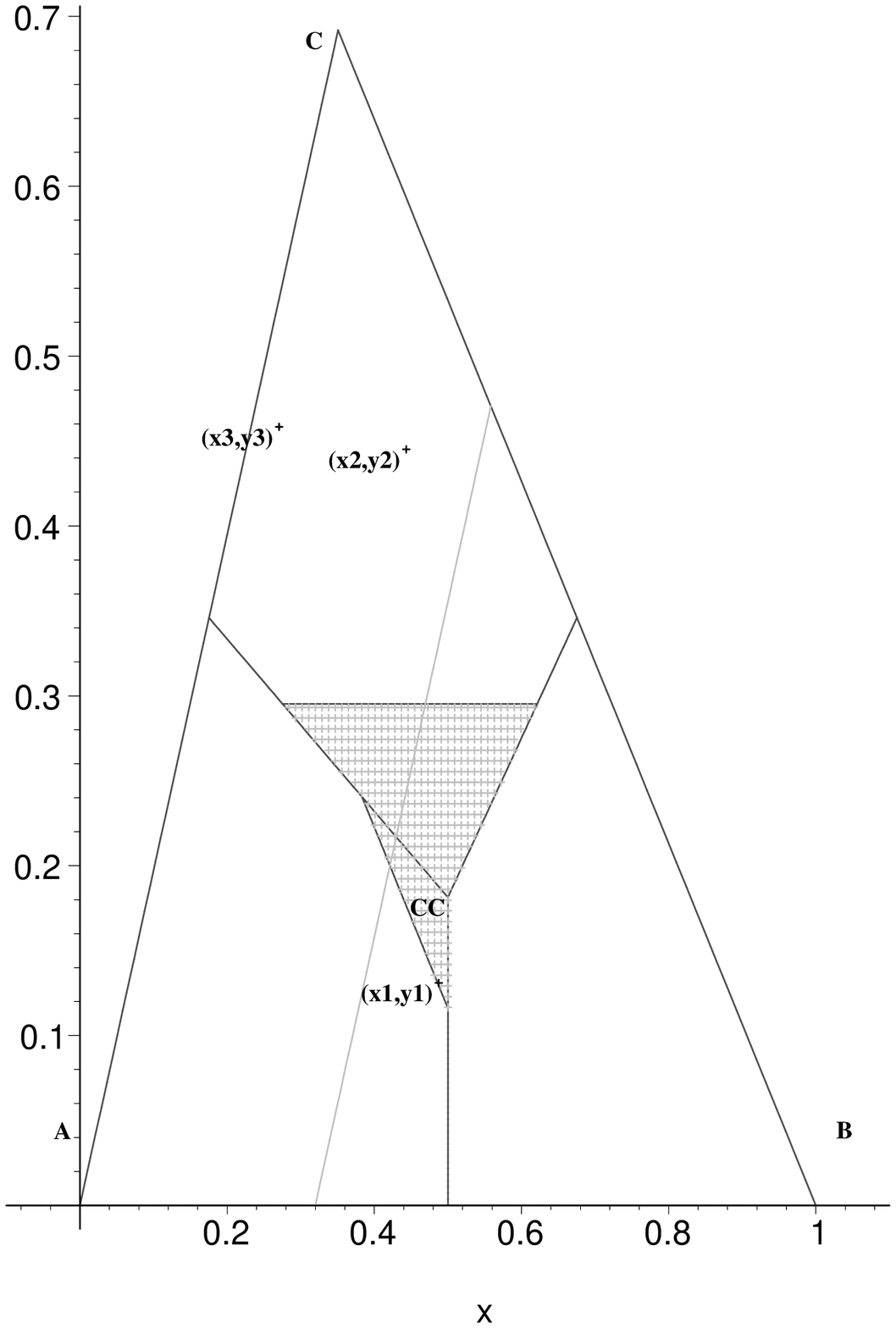, height=140pt , width=200pt}
\psfrag{IC}{\scriptsize{$M_I$}}
\psfrag{G(A)}{}
\psfrag{G(B)}{}
\psfrag{G(C)}{\scriptsize{$\G_1\left(\X_n,\NPE^2,M\right) \cap R^{\perp}_{IC}(\y_3)$}}
\epsfig{figure=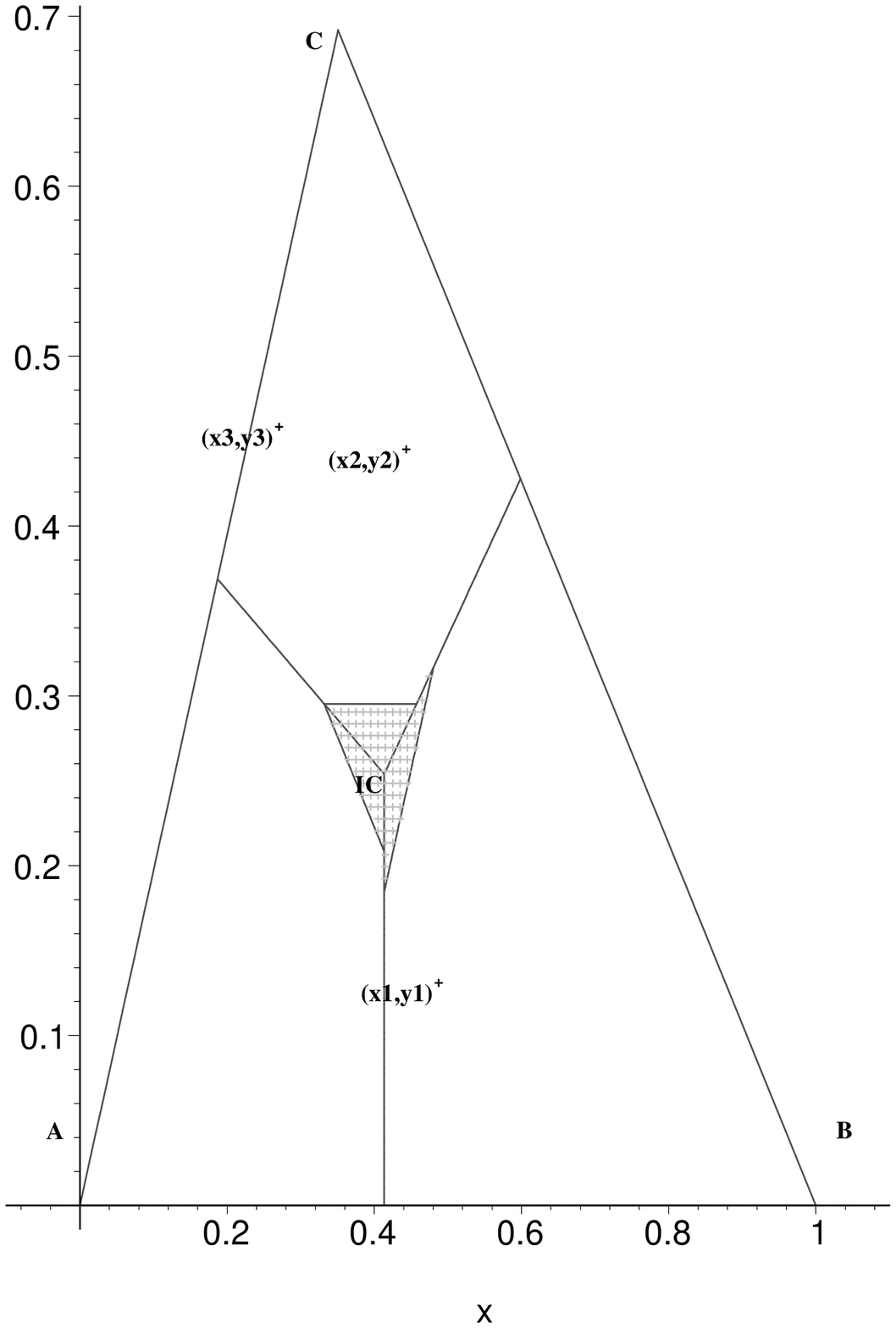, height=140pt , width=200pt}
\end{center}
\caption{$\G_1\left(\X_n,\NPE^{\sqrt{2}},M\right)$ for $M=M_{CC}$ (left) and
$M=M_I$ (right) with three distinct edge extrema.}
\label{fig:G1-NDA-n}
\end{figure}

See Figure \ref{fig:G1-NDA-n} for $\G_1\left(\X_n,\NPE^{\sqrt{2}},M\right)$ where
$M$-vertex regions for $M=M_{CC}$ and $M=M_I$ with orthogonal
projections are used.

\begin{figure}
\begin{center}
\psfrag{A}{\scriptsize{$\y_1$}}
\psfrag{B}{\scriptsize{$\y_2$}}
\psfrag{C}{\scriptsize{$\y_3$}}
\psfrag{(x1,y1)}{\scriptsize{$x$}}
\psfrag{x}{}
\psfrag{CC}{\scriptsize{$M_{CC}$}}
\psfrag{G(A)}{}
\psfrag{G(B)}{}
\psfrag{G(C)}{\scriptsize{$\G_1\left(x,\NPE^2\right) \cap R_{CC}(\y_3)$}}
\epsfig{figure=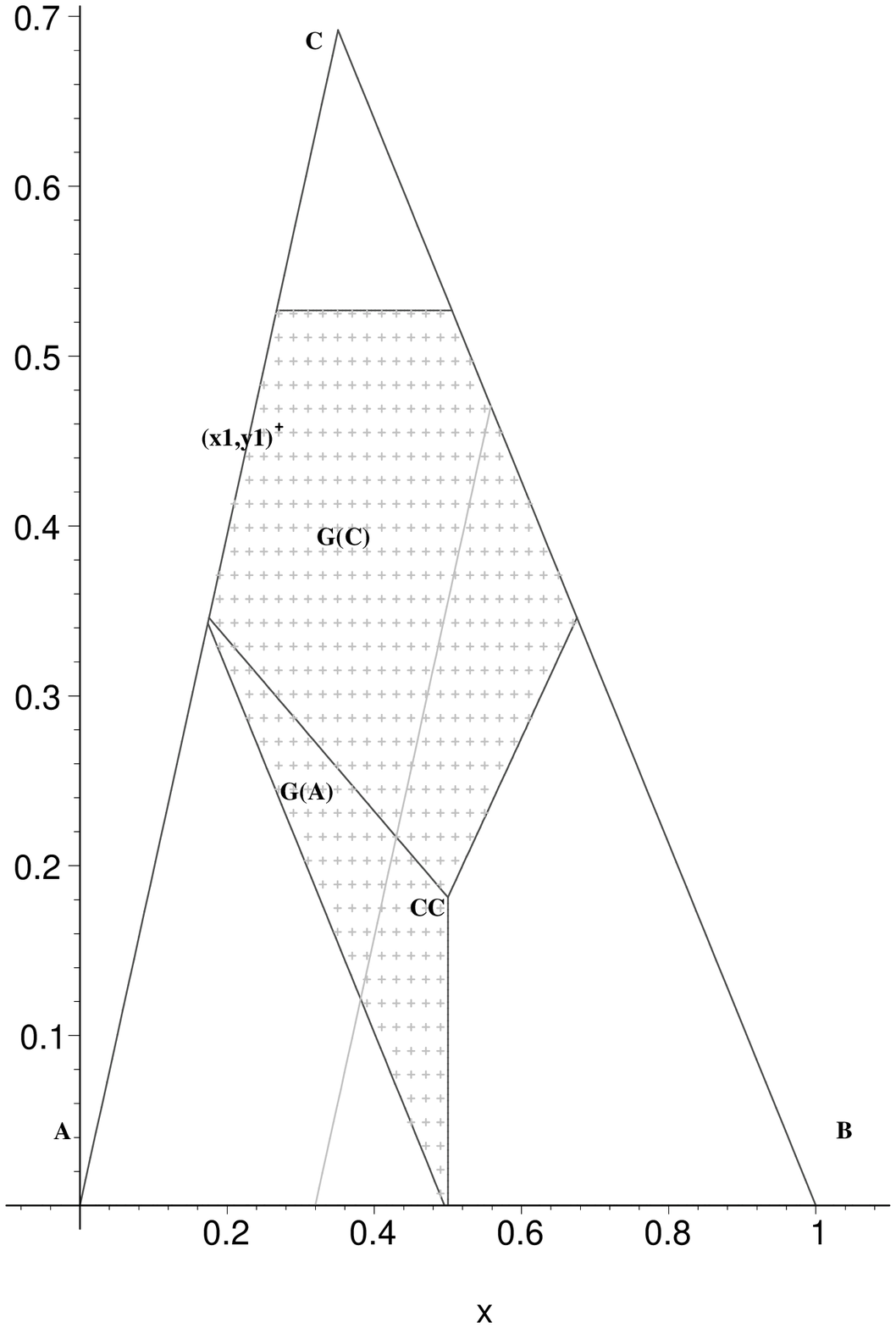, height=140pt , width=200pt}
\psfrag{IC}{\scriptsize{$M_I$}}
\psfrag{G(A)}{}
\psfrag{G(B)}{}
\psfrag{G(C)}{\scriptsize{$\G_1\left(x,\NPE^2\right) \cap R^{\perp}_{IC}(\y_3)$}}
\epsfig{figure=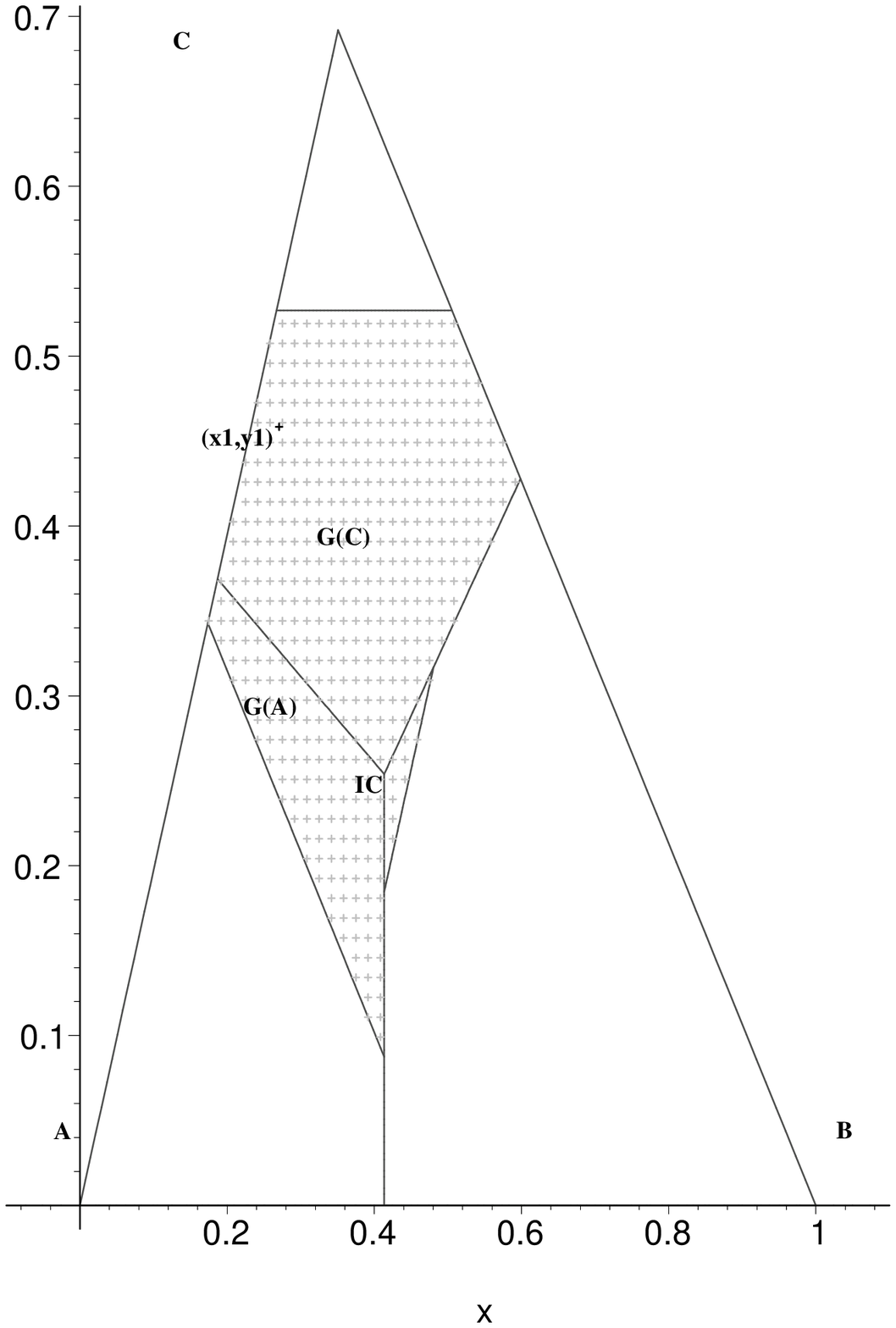, height=140pt , width=200pt}
\end{center}
\caption{$\G_1\left(x,\NPE^{\sqrt{2}},M\right)$ with $x \in R_{CC}(\y_3)$ (left),
$x \in R^{\perp}_{IC}(\y_3)$.} \label{fig:G1-NDA-1}
\end{figure}

\begin{remark}
\begin{itemize}
\item[]
\item
For $r_1 < r_2$, $\G_1\left(x,\NPE^{r_1},M\right) \subseteq
\G_1\left(x,\NPE^{r_2},M\right)$ for all $x \in \TY$ with equality holding only when $x \in \Y_3$.
\item
For $r_1 < r_2 $, $\G_1\left(\X_n,\NPE^{r_1},M\right) \subseteq
\G_1\left(\X_n,\NPE^{r_2},M\right)$ with equality holding only when $\X_n \subseteq \Y_3$ or
$\G_1\left(\X_n,\NPE^{r},M\right)=\emptyset$ for $r=r_1,r_2$.
\item
Suppose $X,Y$ are iid from a continuous distribution $F$ whose support is $\mS(F)\subseteq \TY$.
Then for $r_1 < r_2$, $A\left(\G_1\left(X,\NPE^{r_1},M\right)\right) \le^{ST}
A\left(\G_1\left(Y,\NPE^{r_2},M\right)\right)$.
\item
Suppose $\X_n$ and $\X'_n$ are two random samples from a continuous distribution
$F$ whose support is $\mS(F)\subseteq \TY$.
Then for $r_1 < r_2$, $A(\G_1\left(\X_n,\NPE^{r_1},M\right)) \le^{ST} A(\G_1\left(\X'_n,\NPE^{r_2},M\right))$. $\square$
\end{itemize}
\end{remark}

\begin{remark}
In $\R^d$ with $d>2$, recall $\mathfrak S\left(\Y_m\right)$,
the simplex based on $d+1$ points that do not lie on the same hyperplane.
Furthermore, let $\varrho_i(r,x)$ be the hyperplane such that
$\varrho_i(x) \cap \mathfrak S\left(\Y_m\right) \not=\emptyset$ and
$r\,d(\y_i,\varrho_i(r,x))=d(\y_i,\Upsilon(\y_i,x))$ for $i \in \{1,2,\ldots,d+1$.
Then
$$\G_1\left(x,\NPE^r,M_C\right)\cap R_{CM}(\y_i)=\{z \in R_{CM}(\y_i): d(\y_i,\Upsilon(\y_i,z)) \ge d(\y_i,\varrho_i(r,x)\} \text{
for } i \in \{1,2,3\}. $$
Hence  $\G_1\left(x,\NPE^r,M_C\right)=\bigcup_{i=1}^{d+1} (\G_1\left(x,\NPE^r,M_C\right)\cap R_{CM}(\y_i))$.
Furthermore, it is easy to see that
$\G_1\left(\X_n,\NPE^r,M_C\right)=\bigcap_{i=1}^{d+1}\G_1\left(X_{\varphi_i}(n),\NPE^r,M_C\right)$,
where $X_{\varphi_i}(n)$ is one of the closest points in $\X_n$ to face $\varphi_i$.
$\square$
\end{remark}

\subsection{Central Similarity Proximity Maps}
\label{sec:NCS-region}
The other type of triangular proximity map we introduce is the
central similarity proximity map.
Furthermore, the relative density of the corresponding
PCD will have mathematical tractability.
Alas, the distribution of the domination number of the associated PCD
is still an open problem.

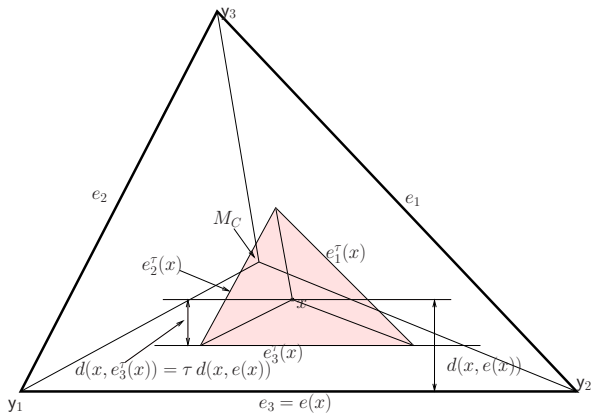
\begin{figure} [ht]
\centering
\scalebox{.35}{\input{N_CSexample2.pstex_t}}
\caption{Construction of central similarity proximity region,
$\NCS^{\tau=1/2}(x,M_C)$ (shaded region).}
\label{fig:ProxMapDefCS}
\end{figure}

For $\tau \in [0,1]$,
define $\NCSt(\cdot,M):=N(\cdot,M;\tau,\Y_3)$ to be the
\emph{central similarity proximity map} with $M$-edge regions as follows;
see also Figure \ref{fig:ProxMapDefCS} with $M=M_C$.
For $x \in \TY\setminus \Y_3$, let $e(x)$ be the
edge in whose region $x$ falls; i.e., $x \in R_M(e(x))$.
If $x$ falls on the boundary of two edge regions,
we assign $e(x)$ arbitrarily.
For ${\tau} \in(0,1]$, the parametrized
central similarity proximity region
$\NCSt(x,M)$ is defined to be the triangle
$T_{\tau}(x)$ with the following properties:
\begin{itemize}
\item[(i)] $T_{\tau}(x)$ has edges $e^{\tau}_i(x)$ parallel to $e_i$
for each $i \in \{1,2,3\}$, and for $x \in R_M(e(x))$,
$d(x,e^{\tau}(x))=\tau\, d(x,e(x))$ and $d(e^{\tau}(x),e(x)) \le d(x,e(x))$
where $d(x,e(x))$ is the Euclidean (perpendicular) distance from $x$ to $e(x)$;
\item[(ii)] $T_{\tau}(x)$ has the same orientation as and is similar to $\TY$;
\item[(iii)] $x$ is the same type of center of $T_{\tau}(x)$ as $M$ is of $\TY$.
\end{itemize}
Note that (i) implies the parametrization of the proximity region,
(ii) explains ``similarity", and
(iii) explains ``central" in the name,
{\em central similarity proximity map}.
For $\tau=0$, we define $\NCS^{\tau=0}(x,M)=\{x\}$ for all $x \in \TY$.
For $x \in \partial(\TY)$,
we define $\NCSt(x,M)=\{x\}$ for all $\tau \in [0,1]$.

Notice that by definition $x \in \NCSt(x,M)$ for all $x \in \TY$.
Furthermore, $\tau \le 1$ implies that $\NCSt(x,M)\subseteq \TY$
for all $x \in \TY$ and $M \in \TY^o$.
For all $ x\in \TY^o \cap R_M(e(x))$,
the edges $e^{\tau}(x)$ and $e(x)$ are coincident iff $\tau=1$.

Notice that $X_i \stackrel{iid}{\sim} F$,
with the additional assumption
that the non-degenerate two-dimensional
probability density function $f$ exists
with support $\mS(F)\subseteq \TY$,
implies that the special case in the construction
of $\NCSt(\cdot)$ ---
$X$ falls on the boundary of two edge regions ---
occurs with probability zero.
Note that for such an $F$, $\NCSt(X,M)$ is a triangle a.s.

Central similarity proximity maps are defined
with $M$-edge regions for $M \in \TY^o$.
In general, for central similarity proximity regions with $M$-edge regions,
the similarity ratio of $\NCSt(x,M)$ to $\TY$ is
$d(x,e^{\tau}(x))/d(M,e(x))$.
See Figure \ref{fig:ProxMapDefCS} for $\NCS^{\tau=1/2}(x,M_C)$ with $e=e_3$.
The functional form of $\NCSt(x,M_C)$ is provided in \cite{ECarXivPCDGeo:2009}.

\subsubsection{Extension of $N^{\tau}_{CS}$ to Higher Dimensions}
\label{sec:NCS-higher-D}
The extension of $N^{\tau}_{CS}$ to $\R^d$ for $d > 2$ is straightforward.
the extension for $M=M_C$ is described, the extension for general $M$ is similar.
Let $\Y_{d+1} = \{\y_1,\y_2,\ldots,\y_{d+1}\}$ be $d+1$ points
that do not lie on the same $(d-1)$-dimensional hyperplane.
Denote the simplex formed by these $d+1$ points as $\mathfrak S(\Y_{d+1})$.
For $\tau \in (0,1]$,
define the central similarity proximity map as follows.
Let $\varphi_i$ be the face opposite vertex $\y_i$
for $i \in \{1,2,\ldots,(d+1)\}$, and ``face regions''
$R_{CM}(\varphi_1),R_{CM}(\varphi_2),\ldots,R_{CM}(\varphi_{d+1})$
partition $\mathfrak S(\Y_{d+1})$ into $d+1$ regions, namely the $d+1$
polytopes with vertices being the center of mass together
with $d$ vertices chosen from $d+1$ vertices.
For $x \in \mathfrak S(\Y_{d+1}) \setminus \Y_{d+1}$, let $\varphi(x)$ be the
face in whose region $x$ falls; $x \in R(\varphi(x))$.
If $x$ falls on the boundary of two face regions,
$\varphi(x)$ is assigned arbitrarily.
For ${\tau} \in (0,1]$,
the central similarity proximity region
$\NCSt(x,M_C)=\mathfrak S_{\tau}(x)$ is defined to be the simplex
$\mathfrak S_{\tau}(x)$ with the following properties:
\begin{itemize}
\item[(i)] $\mathfrak S_{\tau}(x)$ has faces $\varphi^{\tau}_i(x)$ parallel to
$\varphi_i(x)$ for $i \in \{1,2,\ldots,(d+1)\}$, and for $x \in R_{CM}(\varphi(x))$,
$\tau\, d(x,\varphi(x))=d(\varphi^{\tau}(x),x)$
 where $d(x,\varphi(x))$ is the Euclidean (perpendicular)
distance from $x$ to $\varphi(x)$;
\item[(ii)] $\mathfrak S_{\tau}(x)$ has the same orientation
as and similar to $\mathfrak S(\Y_{d+1})$;
\item[(iii)] $x$ is the center of mass of
$\mathfrak S_{\tau}(x)$, as $M_C$ is of $\mathfrak S (\Y_{d+1})$.
Note that $\tau>1$ implies that $x \in \NCSt(x)$.
\end{itemize}

\subsubsection{$\G_1$-Regions for Central Similarity Proximity Maps}
\label{sec:Gamma1-NCS}
\begin{figure} [ht]
\centering
\scalebox{.3}{\input{N_CSGam1.pstex_t}}
\scalebox{.3}{\input{N_CSGam2.pstex_t}}
\scalebox{.3}{\input{N_CSGam3.pstex_t}}
\scalebox{.3}{\input{N_CSGam4.pstex_t}}
\caption{
\label{fig:Gam1Def}
Examples of the four types of the $\G_1$-region,
$\G_1\left(x,\NCS^{\tau=1/2},M_C\right)$ with four distinct
$x \in R_{CM}(e_3)$ (shaded regions).}
\end{figure}
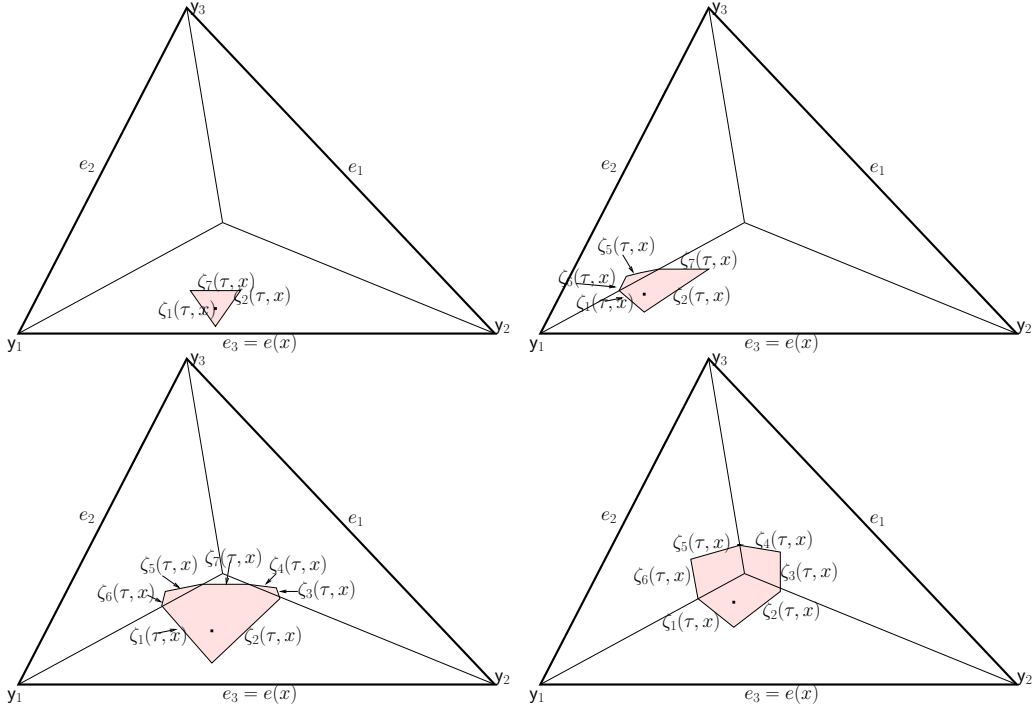

For $\NCSt$, the $\G_1$-region is constructed as follows.
Let $e^{\tau}_i(x)$ be the edge of $T_{\tau}(x)$ parallel to edge $e_i$ for $i \in \{1,2,3\}$.
Now, suppose $u \in R_M(e_3)$ and let $\zeta_i(\tau,x)$ for $i \in \{1,2,\ldots,7\}$ be the lines such that
\begin{align}
\label{eqn:boundarv-G1-NCS}
v \in \zeta_1(\tau,u)\cap R_M(e_3)& \Longrightarrow u \in e^{\tau}_1(v),&
v \in \zeta_5(\tau,u)\cap R_M(e_2)& \Longrightarrow u \in e^{\tau}_3(v), \nonumber\\
v \in \zeta_2(\tau,u)\cap R_M(e_3)& \Longrightarrow u \in e^{\tau}_2(v),&
v \in \zeta_6(\tau,u)\cap R_M(e_2)& \Longrightarrow u \in e^{\tau}_2(v),\\
v \in \zeta_3(\tau,u)\cap R_M(e_1)& \Longrightarrow u \in e^{\tau}_2(v),&
v \in \zeta_7(\tau,u)\cap R_M(e_3)&\Longrightarrow u \in e^{\tau}_3(v). \nonumber\\
v \in \zeta_4(\tau,u)\cap R_M(e_1)& \Longrightarrow u \in e^{\tau}_3(v), & \nonumber
\end{align}
Then $\G_1\left(x,\NCSt,M\right)$ is the region bounded by these lines.
See also Figure \ref{fig:Gam1Def}.
$\G_1\left(x,\NCSt,M\right)$ for $x \in R_M(e_i)$ for $i \in \{1,2\}$ can be described similarly.

Notice that $\tau>0$ implies that $x \in \G_1\left(x,\NCSt,M\right)$.
Furthermore, $\G_1\left(x,\NCSt,M\right)=\{x\}$ iff
(i) $\tau=0$ or (ii) $x \in \partial(\TY)$.

\begin{figure}
\begin{center}
\psfrag{A}{\small{$\y_1$}}
\psfrag{B}{\small{$\y_2$}}
\psfrag{C}{\small{$\y_3$}}
\psfrag{x1}{}
\psfrag{x}{}
\psfrag{D}{}
\psfrag{E}{}
\psfrag{(x1,y1)}{\small{$x$}}
\psfrag{CM}{\small{$M_C$}}
\psfrag{G(AB)}{}
\psfrag{G(BC)}{}
\psfrag{G(AC)}{}
\epsfig{figure=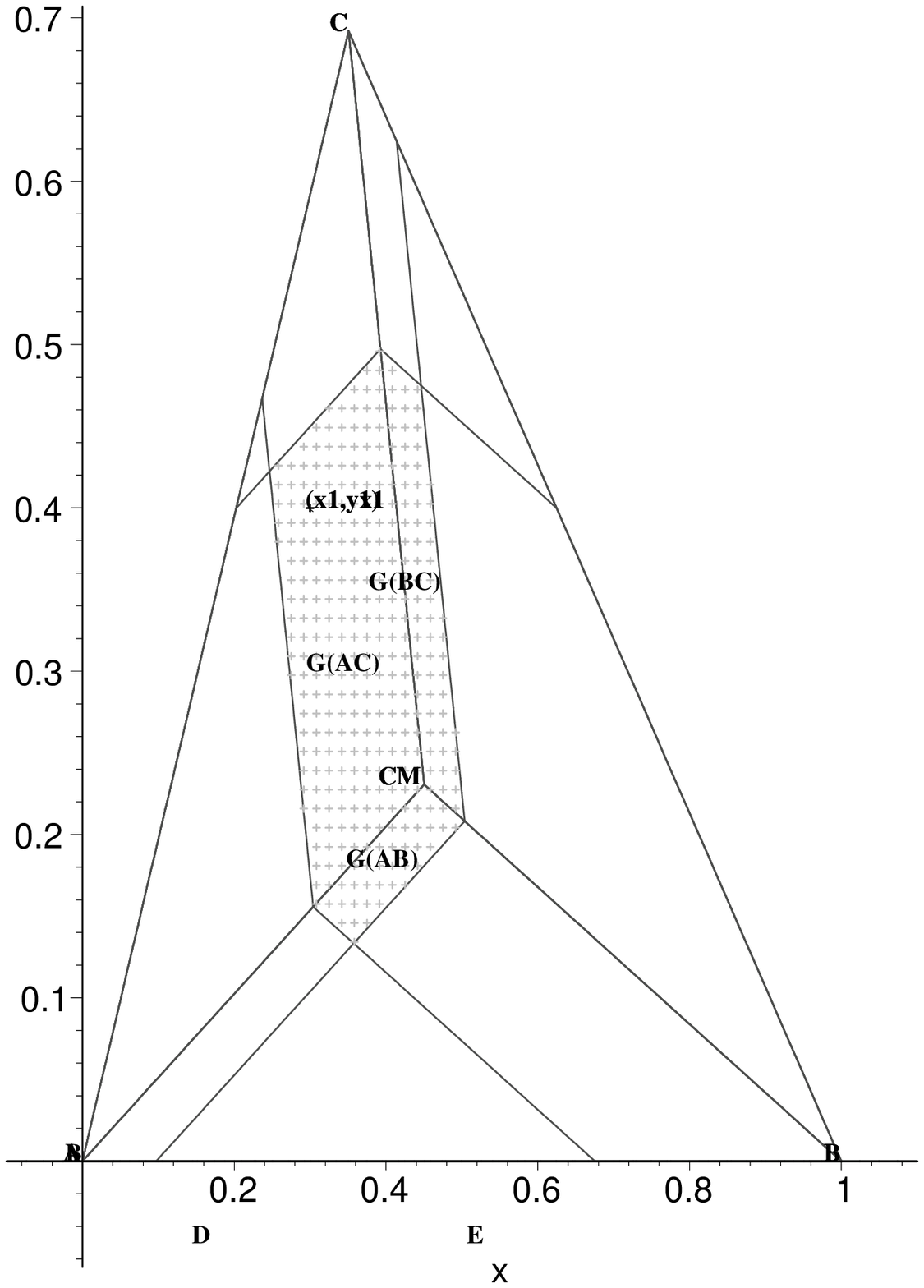, height=140pt , width=200pt}
\psfrag{xa}{\small{$x_{e_1}$}}
\psfrag{xb}{\small{$x_{e_2}$}}
\psfrag{xc}{\small{$x_{e_3}$}}
\epsfig{figure=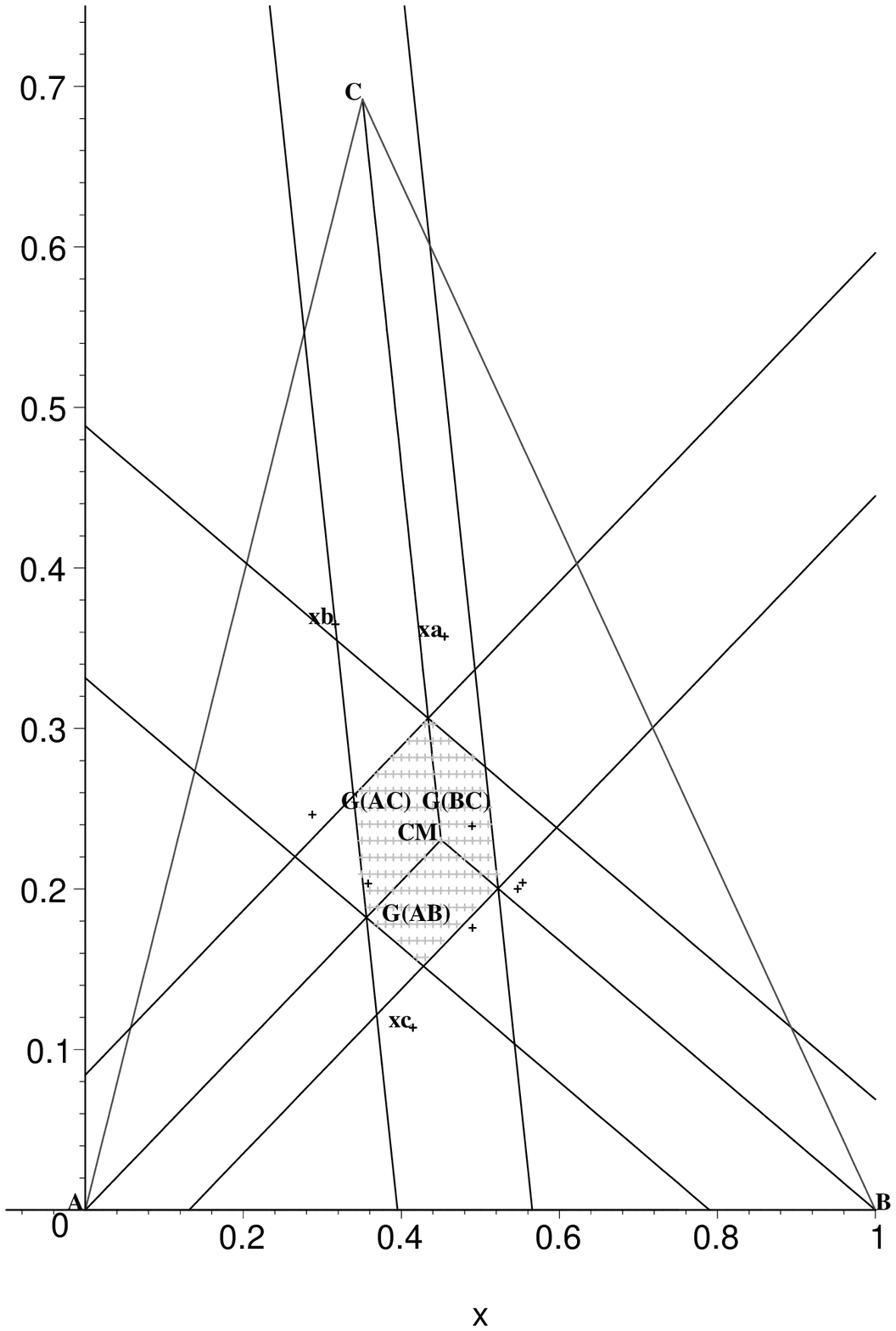, height=140pt , width=200pt}
\end{center}
\caption{
\label{fig:G1-NCS-CM-1n}
The $\G_1$-region $\G_1\left(x,\NCSt,M_C\right)$
with $x \in R_M(e_2)$ (left)
and $\G_1\left(X_n,\NCSt,M_C\right)$ with $n > 1$ (right).}
\end{figure}

The $\G_1$-region $\G_1\left(x,\NCSt,M\right)$ is a convex
$k$-gon with $3 \le k \le 7$ vertices.
In particular, for $\tau=1$, $\G_1\left(x,\NCS^{\tau=1},M\right)$ is a convex hexagon.
See Figure \ref{fig:G1-NCS-CM-1n} (left).

Then the functional form of $\zeta_i(x)$ in Equation \eqref{eqn:boundarv-G1-NCS} that define the boundary of
$\G_1\left(x=(x_0,y_0),\NCSt\right)$ with $M$-central proximity regions
in the basic triangle $T_b$ for $x \in R_M(e_3)$ is
\begin{align*}
\zeta_1(x)&=\frac{m_2\,(c_2\,(x-x_0)+y_0\,c_1)}{m_2\,c_1\,(1-\tau)+c_2\,\tau\,m_1}, \;\;
\zeta_2(x)=\frac{m_2\,(c_2\,(x_0-x)+y_0\,(1-c_1))}{m_2\,(1-c_1)\,(1-\tau)+c_2\,\tau\,(1-m_1)},\\
\zeta_3(x)&=\frac{m_2\,c_2\,(x-1)+y_0\,(m_2\,(1-c_1)-\,c_2\,(1-m_2))}{c_1\,(1-c_1)\,(1-\tau)},\\
\zeta_4(x)&=\Bigl[c_1\,(m_2\,(1-c_1)+c_2\,(1-m_1))\,y_0+c_2\,(m_2\,(1-c_1)+c_2\,(1-m_1))\,x_0+c_2\,(\tau\,(m_2\,c_1-m_1\,c_2)\\
&+c_2\,(1-m_1)+m_2\,(1-c_1))\,x+c_2\,\tau\,(m_1\,c_2-m_2\,c_1)\Bigr]\Big/\Bigl[c_2\,(\tau-c_1\,(1-\tau))\,m_1\\
&+c_1\,(1-c_1)\,(1-\tau)\,m_2+c_1\,c_2\Bigr],\\
\zeta_5(x)&=\frac{\tau\,m_2\,c_2\,x-y_0\,(c_1\,m_2-c_2\,m_1)}{c_2\,m_1-c_1\,m_2\,(1-\tau)},\;\;
\zeta_7(x)= \frac{y_0}{1-\tau}\\
\zeta_6(x)&=\Bigl[(1-c_1)\,c_1\,(m_2\,(1-c_1)-c_2\,(1-m_1))\,y_0+c_1\,c_2\,(m_2\,(1-c_1)-c_2\,(1-m_1))\,x_0+(c_2^2\,c_1\,(1-\tau)\,m_1\\
&+c_2\,\tau\,(1-c_1)+(1-\tau)\,(c_1\,c_2-c_1^2\,c_2)\,m_2-c_2^2\,\tau\,(1-c_1)-c_2^2\,c_1\,(1-\tau))\,x+
c_2^2\,\tau\,(m_1-c_1)\Bigr]\Big/\Bigl[(1-c_1)\,\\
&(c_2\,(c_2\,(1-\tau)-\tau)\,m_1+c_1\,(1-c_1)\,(1-\tau)\,m_2+c_1\,c_2\,(1-2\,\tau))\Bigr]\\
\end{align*}

\begin{proposition}
Let $\X_n$ be a random sample from $\UT$.
For central similarity proximity maps with $M$-edge regions
(by definition $M \in \TY^o$) and $\tau>0$, $\eta_n\left(\NCSt\right) \le 3$ with equality holding
with positive probability for $n \ge 3$.
\end{proposition}
\noindent \textbf{Proof:}
Let $M \in \TY^o$ and $\tau >0$.
Then given $\X_n$, for $i \in \{1,2,3\}$,
we have
\begin{multline*}
\G_1\left(\X_n,\NCSt,M)\right)\cap R_M(e_i)=\bigcap_{i=1}^n \left[ \G_1\left(X_i,\NCSt,M\right)\cap R_M(e_i) \right]\\
=\G_1\left(X_{e_j}(n),\NCSt,M\right)\cap \G_1 \left(X_{e_k}(n),\NCSt ,M\right) \cap R_M(e_i)
\end{multline*}
where $j,k \not=i$, since for $x \in R_M(e_i)$,
if $\{X_{e_k}(n),X_{e_l}(n)\} \subset \NCSt(x,M)$,
then $\X_n \subset \NCSt(x,M)$.
Hence for each edge we need the edge extrema with respect to the other edges,
then the minimum active set is $S_M=\{X_{e_1}(n),X_{e_2}(n),X_{e_3}(n)\}$,
hence $\eta_n\left(\NCSt\right) \le 3$.
Furthermore, for the random sample $\X_n$, $X_e(n)$ is unique for each
edge $e$ with probability 1 and there are
three distinct edge extrema with positive probability (see Theorem \ref{thm:distinct-edge-ext}).
Hence $P(\eta_n\left(\NCSt\right) = 3)>0$ for $n \ge 3$.
$\blacksquare$

\begin{figure}
\centering
\psfrag{A}{\small{$\y_1$}}
\psfrag{B}{\small{$\y_2$}}
\psfrag{C}{\small{$\y_3$}}
\psfrag{x1}{}
\psfrag{x}{}
\psfrag{D}{}
\psfrag{E}{}
\psfrag{F}{}
\psfrag{G}{}
\psfrag{CM}{\small{$M_C$}}
\psfrag{G(AB)}{\small{$G(\overline{AB})$}}
\epsfig{figure=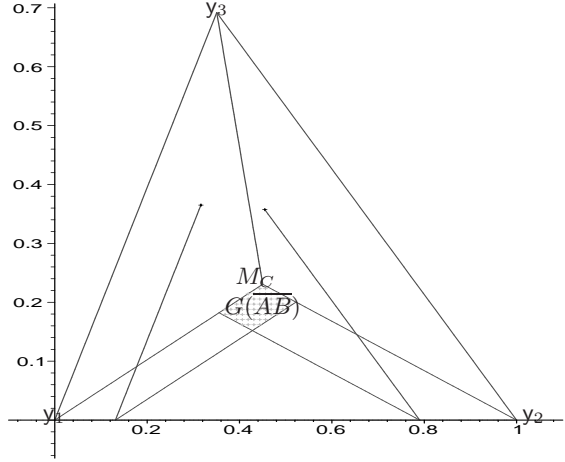, height=175pt , width=210pt}
\caption{
\label{fig:CS-G1-region}
$G(\overline{AB})$ region for $x_1$ and $\G_1$-region for $n \ge 1$}
\end{figure}

Note that  for $\tau>0$ and $\X_n$ a random sample from $\UT$,
$P(\eta(\X_n,\NCSt)=3) \rightarrow 1$ as
$n \rightarrow \infty$,
since the edge extrema are distinct with probability 1 as
$n \rightarrow \infty$.
For $\tau=1$, the $\G_1$-region can be determined by the edge extrema,
in particular $\G_1\left(\X_n,\NCS^{\tau=1},M\right) \cap R_M(\y_i)$ can be determined by
the edge extrema $X_{e_j}$ and $X_{e_k}$ for $i,j,k$ are all distinct.
See Figure \ref{fig:CS-G1-region} for an example.
Here $\G_1\left(\X_n,\NCS^{\tau=1},M\right) \not=  \emptyset$ for all $n$,
because by construction $M \in \G_1\left(\X_n,\NCS^{\tau=1},M\right)$ since $\NCS^{\tau=1}(M)=\TY$.
Furthermore, $\eta_n\left(\NCSt\right) \stackrel{d}{=}\eta_n\left(\NPE^r\right)$ for
all $(r,\tau) \in [1,\infty) \times (0,1]$.

\begin{proposition}
The $\G_1$-region  $\G_1\left(B, \NCS^{\tau=1},M\right)$ is a convex hexagon
for any set $B$ of size $n$ in $\TY^o$.
\end{proposition}
\noindent \textbf{Proof:}
This follows from the fact that each $\zeta_i(x)$
for $i \in \{ 1,2,\ldots,6\}$
is parallel to a line joining $\y_i$ to $M$ for $\tau=1$
($\zeta_7(x)$ is not used in construction of $\G_1\left(B,\NCS^{\tau=1},M\right)$).
See Figure \ref{fig:G1-NCS-CM-1n} (right) for an example.
$\blacksquare$

With $M=M_C$, the functional forms of the lines that determine the
boundary of  $\G_1\left(x=(x_0,y_0),\NCSt\right)$ in
Equation \eqref{eqn:boundarv-G1-NCS} for $T_b$ with $x \in R_{CM}(e_3)$ are given by
\begin{align*}
\zeta_1(x)&=\frac{y_0\,c_1+c_2\,(x-x_0)}{c_1+\tau},\;
\zeta_2(x)=\frac{c_2\,(x_0-x)+y_0\,(1-c_1)}{\tau+1-c_1},\\
\zeta_3(x)&=\frac{c_2\,\tau\,(1-x)+y_0}{1+\tau\,(1-c_1)},\;\;\;\;
\zeta_4(x)=\frac{\tau\,c_2\,(1-x)+c_2\,(x_0-x)-y_0\,c_1}{\tau\,(1-c_1)-c_1},\\
\zeta_5(x)&=\frac{c_2\,x+y_0}{1+\tau\,c_1},\;\;\;\;\;\;\;\;\;\;\;\;\;\;\;\;
\zeta_7(x)= \frac{y_0}{1-\tau},\\
\zeta_6(x)&=\frac{c_1\,(1-c_1)\,y_0+c_1\,c_2\,x_0-c_2\,(2\,\tau\,(1-c_1)+c_1\,(1-\tau))\,x+\tau\,c_2\,(1-2\,c_1)}{c_1\,(1-c_1)\,(1-\tau)}.
\end{align*}

For $n > 1$ with $M=M_C$, let $X_{e_i}(n)=x_{e_i}=(u_i,w_i)$ be given
 for $i \in \{1,2,3\}$.
Then the functional forms of the lines that determine the boundary
of $\G_1\left(\X_n,\NCSt,M_C\right)$ are given by

\begin{tabular}{ll}
$\displaystyle \zeta_1(x)=\frac{w_2\,c_1+c_2\,(x-u_2)}{c_1+\tau}$,   & $\displaystyle \zeta_3(x)=\frac{c_2\,\tau\,(1-x)+y_0}{1+\tau\,(1-c_1)}$,\\
$\displaystyle \zeta_2(x)=\frac{c_2\,(u_1-x)+w_1\,(1-c_1)}{\tau+1-c_1}$, & $\displaystyle \zeta_5(x)=\frac{c_2\,x+w_3}{1+\tau\,c_1}$,\\
$\displaystyle \zeta_4(x)=\frac{\tau\,c_2\,(1-x)+c_2\,(u_3-x)-w_3\,c_1}{\tau\,(1-c_1)-c_1}$, &  $\displaystyle \zeta_7(x)= \frac{w_3}{1-\tau}$,  \\
\multicolumn{2}{c}{$\displaystyle \zeta_6(x)=\frac{c_1\,(1-c_1)\,w_1+c_1\,c_2\,u_1-c_2\,(2\,\tau\,(1-c_1)+c_1\,(1-\tau))\,x+\tau\,c_2\,(1-2\,c_1)}{c_1\,(1-c_1)\,(1-\tau)}$.} \\\\
\end{tabular}


See Figure \ref{fig:G1-NCS-CM-1n} for  $\G_1\left(x,\NCS^{\tau=1},M_C\right)$
and $\G_1\left(\X_n,\NCS^{\tau=1},M_C\right)$.

\begin{remark}
\begin{itemize}
\item[]
\item
For $\tau_1 < \tau_2$,
$\G_1\left(x,\NCS^{\tau_1},M\right) \subseteq \G_1\left(x,\NCS^{\tau_2},M\right)$
for all $x \in \TY$ with equality holding only when $x \in \partial(\TY)$.
\item
Let $\X_n$ be a set of $n$ points in $\TY$.
Then for $\tau_1 < \tau_2$,
$\G_1\left(\X_n,\NCS^{\tau_1},M\right) \subseteq \G_1\left(\X_n,\NCS^{\tau_2},M\right)$
with equality holding only when $\X_n \subset \partial(\TY)$.
\item
Suppose $X,Y$ are iid from a continuous distribution $F$ whose support is $\mS(F)\subseteq \TY$.
Then for $\tau_1 < \tau_2$, $A\left(\G_1\left(X,\NCS^{\tau_1},M\right)\right) \le^{ST}
A\left(\G_1\left(Y,\NCS^{\tau_2},M\right)\right)$.
\item
Suppose $\X_n$ and $\X'_n$ are two random samples from a continuous distribution
$F$ whose support is $\mS(F)\subseteq \TY$.
Then for $\tau_1 < \tau_2$, $A(\G_1\left(\X_n,\NCS^{\tau_1},M\right)) \le^{ST}
A(\G_1\left(\X'_n,\NCS^{\tau_2},M\right))$. $\square$
\end{itemize}
\end{remark}


\section{Investigation of the Proximity Regions and the Associated PCDs Using the Auxiliary Tools}
\label{sec:characterize-PCDs}

\subsection{Characterization of Proximity Maps Using $\eta$-Values}
By definition, it is trivial to show that that the minimum number of
points to describe the $\G_1$-region,
$\eta_n(N) \le n$ for any proximity map.
We have improved the upper bound for
$\NPE^r$ and $\NCSt$: $\eta(\X_n,\NPE^r) \le 3$
and $\eta(\X_n,\NCSt) \le 3$.
However, finding such an improvement does not hold for $\eta_n(\NAS)$;
that is, finding a $k<n$ such that $\eta_n(\NAS)\le k$
for all $\X_n$ is still an open problem.

Below we state a condition for $N(\cdot,M)$ defined with $M$-vertex regions to have
$\eta_n(N) \le 3$ for $\X_n$ with support in $\TY$.
\begin{theorem}
\label{thm:eta<=3-vertex}
Suppose $N(\cdot,M)$ is a proximity region defined with $M$-vertex regions
and $B$ is a set of $n$ distinct points in $\TY$.
Then $\eta(B,N) \le 3$ if
\begin{itemize}
\item[(i)]
for each $\y_i \in \Y_3$ there exists a point $x(\y_i) \in B$ (i.e., related to $\y_i$)
 such that $\G_1(B,N) \cap R_M(\y_i)=\G_1\left(x(\y_i),N\right) \cap R_M(\y_i)$,
\item[] or
\item[(ii)]
there exist points $x(\y_j),\,x(\y_k) \in B$ such that
$\G_1(B,N) \cap R_M(\y_i)=\G_1\left(x(\y_j),N\right) \cap \G_1\left(x(\y_k),N\right) \cap R_M(\y_i)$
for $j,k\not=i$ with $i \in \{1,2,3\}$ and $(j,k) \in \{(1,2),(1,3),(2,3)\}$.
\end{itemize}
\end{theorem}
\noindent \textbf{Proof:}
Let $B=\{x_1,x_2,\ldots,x_n\} \subset \TY$.
\begin{itemize}
\item[(i)]
Suppose there exists a point $x(\y_i) \in B$ such that
$\G_1(B,N) \cap R_M(\y_i)=\G_1\left(x(\y_i),N\right) \cap R_M(\y_i)$ for each $i \in \{1,2,3\}$.
Then
\begin{multline*}
\G_1(B,N) \cap R_M(\y_i)=\G_1\left(x(\y_i),N\right) \cap R_M(\y_i)=\\
\bigcap_{i=1}^n[\G_1\left(x_i,N\right) \cap R_M(\y_i)]=
\bigcap_{j=1}^3[\G_1\left(x(\y_j),N\right) \cap R_M(\y_i)]=
\left[ \bigcap_{j=1}^3\G_1\left(x(\y_j),N\right) \right] \cap R_M(\y_i)
\end{multline*}
and
$$\G_1(B,N)= \bigcup_{i=1}^3[\G_1(B,N) \cap R_M(\y_i)]=
\bigcup_{i=1}^3 \left(\left[ \bigcap_{j=1}^3\G_1\left(x(\y_j),N\right) \right] \cap R_M(\y_i)\right).$$
Then, we get
$$\G_1(B,N)= \bigcap_{i=1}^3\G_1\left(x(\y_i),N\right).$$
Hence, the minimum active set $S_M \subseteq \{x(\y_1),x(\y_2),x(\y_3)\}$,
which implies $\eta(B,N) \le 3$.
The $\eta$-value $\eta(B,N)<3$ will hold if $x(\y_i)$ are not all distinct.
\item[(ii)]
Suppose there exist points $x(\y_j)$ and $x(\y_k)$ such that
$\G_1\left(B,N\right) \cap R_M(\y_i)=
\G_1\left(x(\y_j),N\right) \cap \G_1\left(x(\y_k),N\right) \cap R_M(\y_i)$ for $j,k\not=i$.
Then
\begin{multline*}
\G_1\left(B,N\right) \cap R_M(\y_i)=
\G_1\left(x(\y_j),N\right) \cap \G_1\left(x(\y_k),N\right) \cap R_M(\y_i)=\\
\bigcap_{i=1}^n[\G_1\left(x_i,N\right) \cap R_M(\y_i)]=
\bigcap_{q=1}^3[\G_1\left(x(\y_q),N\right) \cap R_M(\y_i)]=
\left[\bigcap_{q=1}^3\G_1\left(x(\y_q),N\right)\right] \cap R_M(\y_i)
\end{multline*}
and
$$\G_1(B,N)= \bigcup_{i=1}^3[\G_1(B,N) \cap R_M(\y_i)]=
\bigcup_{i=1}^3 \left( \left[\bigcap_{q=1}^3\G_1\left(x(\y_q),N\right)\right] \cap R_M(\y_i) \right).$$
Then, we get
$\G_1(B,N)= \bigcap_{i=1}^3\G_1\left(x(\y_i),N\right).$
Hence, the minimum active set $S_M \subseteq \{x(\y_1),x(\y_2),x(\y_3)\}$
which implies $\eta(B,N) \le 3$. $\blacksquare$
\end{itemize}

Notice that $\NPE^r$ satisfies condition (i) Theorem \ref{thm:eta<=3-vertex}.

Below we state some conditions for $N(\cdot,M)$ defined with $M$-edge regions to have
$\eta$-value less than equal to 3.
\begin{theorem}
\label{thm:eta<=3-edge}
Suppose $N(\cdot,M)$ is  a proximity region defined with $M$-edge regions
 and $B$ is set of $n$ distinct points in $\TY$.
Then $\eta(B,N) \le 3$ if
\begin{itemize}
\item[(i)]
for each $e_i \in \{e_1,e_2,e_3\}$,
there exists a point $x(e_i) \in B$ such that
$\G_1(B,N) \cap R_M(e_i)=\G_1\left(x(e_i),N\right) \cap R_M(e_i)$,
\item[] or
\item[(ii)]
there exist points $x(e_j),\,x(e_k) \in B$ such that
$\G_1(B,N) \cap R_M(e_i)=\G_1\left(x(e_j),N\right) \cap \G_1\left(x(e_k),N\right) \cap R_M(e_i)$
for $j,k\not=i$ with $i \in \{1,2,3\}$ and $(j,k) \in \{(1,2),(1,3),(2,3)\}$.
\end{itemize}
\end{theorem}
\noindent \textbf{Proof:}
Let $B=\{x_1,x_2,\ldots,x_n\} \subset \TY$.
\begin{itemize}
\item[(i)]
Suppose there exists a point $x(e_i) \in B$ such that
$\G_1(B,N) \cap R_M(e_i)=\G_1\left(x(e_i),N\right) \cap R_M(e_i)$
for each $i \in \{1,2,3\}$.
Then
\begin{multline*}
\G_1(B,N) \cap R_M(e_i)=\G_1\left(x(e_i),N\right) \cap R_M(e_i)=\\
\bigcap_{i=1}^n[\G_1\left(x_i,N\right) \cap R_M(e_i)]=
\bigcap_{j=1}^3[\G_1\left(x(e_j),N\right) \cap R_M(e_i)]=
\left[ \bigcap_{j=1}^3\G_1\left(x(e_j),N\right) \right] \cap R_M(e_i)
\end{multline*}
and
$$\G_1(B,N)= \bigcup_{i=1}^3[\G_1(B,N) \cap R_M(e_i)]=
\bigcup_{i=1}^3 \left(\left[ \bigcap_{j=1}^3\G_1\left(x(e_j),N\right) \right] \cap R_M(e_i)\right).$$
Then, we get
$$\G_1\left(\X_n,N\right)= \bigcap_{j=1}^3 \G_1\left(x(e_j),N\right).$$
Hence, the minimum active set $S_M \subseteq \{x(e_1),x(e_2),x(e_3)\}$
which implies $\eta_n(N) \le 3$.
\item[(ii)]
Suppose there exist points $x(e_j)$ and $x(e_k)$ such that
$\G_1\left(\X_n,N\right) \cap R_M(e_i)=
\G_1\left(x(e_j),N\right) \cap \G_1\left(x(e_k),N\right) \cap R_M(e_i)$ for $j,k\not=i$.
Then
\begin{multline*}
\G_1\left(B,N\right) \cap R_M(e_i)=
\G_1\left(x(e_j),N\right) \cap \G_1\left(x(e_k),N\right) \cap R_M(e_i)=\\
\bigcap_{i=1}^n[\G_1\left(x_i,N\right) \cap R_M(e_i)]=
\bigcap_{q=1}^3[\G_1\left(x(e_q),N\right) \cap R_M(e_i)]=
\left[\bigcap_{q=1}^3\G_1\left(x(e_q),N\right)\right] \cap R_M(e_i)
\end{multline*}
and
$$\G_1(B,N)= \bigcup_{i=1}^3[\G_1(B,N) \cap R_M(e_i)]=
\bigcup_{i=1}^3 \left( \left[\bigcap_{q=1}^3\G_1\left(x(e_q),N\right)\right] \cap R_M(e_i) \right).$$
Then, we get
$\G_1(B,N)= \bigcap_{i=1}^3\G_1\left(x(e_i),N\right).$
Hence, the minimum active set $S_M \subseteq \{x(e_1),x(e_2),x(e_3)\}$
which implies $\eta(B,N) \le 3$. $\blacksquare$
\end{itemize}

Notice that $\NCSt$ satisfies condition (ii) in Theorem \ref{thm:eta<=3-edge}.

\subsection{The Behavior of $\G_1\left(\X_n,N\right)$ for the Proximity Maps in $\TY$}
\label{sec:behaviour-G1-Regions}
In Section \ref{sec:gamma1-regions}, we have investigated the behavior
of $\G_1\left(\X_n,N\right)$ for general proximity maps in $\Omega$.
The assertions made about $\G_1$-regions will be stronger
for the proximity regions we have defined,
i.e., for $\NAS$, $\NPE^r$, and $\NCSt$, compared
to the general assertions in Section \ref{sec:gamma1-regions}.
One property enjoyed by these proximity maps is that the region
$N(x)$ gets larger as $x$ moves along a line from $\partial(\TY)$ to $\RS(N)$
in a region with positive $\R^2$-Lebesgue measure.
So the modifications of the assertions in Section \ref{sec:gamma1-regions}
also hold for $\left\{ \NS, \NAS, \NPE^r, \NCSt \right\}$.
In particular, we have a stronger result than the one in
Proposition \ref{prop:RSsubsetG1} in the sense that,
$\RS(N)$ is a proper subset of $\G_1\left(\X_n,N\right)$ as shown below.

\begin{proposition}
For each type of proximity map $N \in \left\{ \NS, \NAS, \NPE^r, \NCSt \right\}$
and any random sample $\X_n=\bigl\{ X_1,X_2,\ldots,X_n \bigr\}$ from
a continuous distribution $F$ on $\TY$,
if $\RS(N) \not=\emptyset$, then $\RS(N) \subsetneq \G_1\left(\X_n,N\right)$ a.s. for each $n<\infty$.
\end{proposition}
\noindent \textbf{Proof:}
We have shown that $\RS(N) \subseteq \G_1\left(\X_n,N\right)$ (see Proposition \ref{prop:RSsubsetG1}).
Moreover,
$\RS(N) \not=\emptyset$ for $\NAS$, $\NPE^r$ with $r>3/2$, and $\NCS^{\tau=1}$.
For these proximity regions, $\RS(N)=\G_1\left(\X_n,N\right)$ with probability 0 for each
finite $n$ since
\begin{itemize}
\item[(i)]
for $N \in \left\{ \NS, \NAS \right\}$,
$\G_1\left(\X_n,N\right)=\RS\left(N\right)$ iff $X_{\y}(n)=\y$ for each $\y \in \Y_3$
which happens with probability 0,
\item[(ii)]
for $N \in \left\{\NCS^{\tau=1},\, \NPE^{r>3/2} \right\}$,
$\G_1\left(\X_n,N,M\right)=\RS\left(N,M\right)$
iff $X_{e_i}(n) \in e_i$ for each $i \in \{1,2,3\}$
which happens with probability 0.
\end{itemize}
Furthermore, for $\NPE^{r}$ with $r<3/2$, $\RS\left( \NPE^r,M \right)\not=\emptyset$
iff $M \notin (\Tr)^o$ (see Equation \eqref{eqn:T^r-def} for $\Tr$),
say $M=(m_x,m_y)$ is such that $d(\zeta(m_x,\y_2),\y_2) < d(\y_2,e_2)/r$.
Then $\G_1\left(\X_n,\NPE^r ,M\right)=\RS\left( \NPE^r,M \right)$
iff $X_{e_2}(n)\in e_2$ which happens with probability 0.
Similarly the same result also holds for edges $e_1$ and $e_3$.
$\blacksquare$

Note that, if $\RS(N) =\emptyset$ and $\X_n$ is a random sample from a continuous distribution
on $\TY$, then $\G_1\left(\X_n,N\right) = \emptyset$ a.s. as $n \rightarrow \infty$.
In particular, this holds for $\NPE^r(\cdot,M)$  with $r<3/2$ and $M \in (\Tr)^o$.
Lemma \ref{lem:gamma1-nonincreasing} holds as stated.
In Lemma \ref{lem:gamma1-nonincreasing},
we have shown that $\G_1\left(\X(n),N\right)$ is non-increasing.
Furthermore, for the proximity regions $\left\{\NS,\NAS,\NPE^r,\NCSt\right\}$,
we can state that $\G_1\left(\X(n+1),N\right) \subsetneq \G_1\left(\X(n),N\right)$
with positive probability,
since the new point in $\X(n+1)$ has positive probability
to fall closer to the subset of $\TY$
that defines $\RS(N)$ (e.g., $\partial(\TY)$).
The general results in Theorems \ref{thm:G1-as-conv1} and \ref{thm:AG1-converge-in-p}
and Proposition \ref{prop:stoch-order-G1} hold for the proximity maps
$\left\{\NS,\NAS,\NPE^r,\NCSt\right\}$ also.
In particular we have the following corollaries to these results.

\noindent
\emph{\textbf{Corollary to Theorem \ref{thm:G1-as-conv1}:}
Given a sequence of random variables $X_1,X_2,\ldots \stackrel{iid}{\sim}\UT$,
let $\X(n):=\X(n-1) \cup \{X_n\}$ for $n=0,1,2,\ldots$ with $\X(0):= \emptyset$.
Then for each $N \in \left\{\NS,\NAS,\NPE^r,\NCSt\right\}$,
$\G_1\left(\X(n),N\right) \downarrow \RS(N)$ as $n \rightarrow \infty$ a.s.,
in the sense that $\G_1\left(\X(n+1),N\right) \subseteq \G_1\left(\X(n),N\right)$ and
$A(\G_1\left(\X(n),N\right)\setminus \RS(N)) \downarrow 0$ a.s.
}

\noindent
\emph{\textbf{Corollary to Proposition \ref{prop:stoch-order-G1}:}
For positive integers $m > n$,
let $\X_m$ and $\X_n$ be two random samples from $\UT$.
Then
$A\left(\G_1\left(\X_m,N,M\right)\right) \le^{ST}A\left(\G_1\left(\X_n,N,M\right)\right)$
for proximity maps $N \in \left\{\NPE^r,\NCSt\right\}$,
for $r \in [1,\infty)$ and $\tau \in (0,1]$.
}

For $N \in \left\{ \NS, \NAS \right\}$,
the stochastic ordering of  $A\left(\G_1\left(\X_n,N,M\right)\right)$
is still an open problem, although we conjecture that
$A\left(\G_1\left(\X_m,N,M\right)\right) \le^{ST} A\left(\G_1\left(\X_n,N,M\right)\right)$
for $\X_m$ and $\X_n$ two random samples from $\UT$ with $m>n$.


\noindent
\emph{\textbf{Corollary to Theorem \ref{thm:AG1-converge-in-p}:}
Let $\{\X_n\}_{n=1}^{\infty}$ be a sequence of data sets which are iid  $\UT$.
Then $\G_1\left(\X_n,N\right) \xrightarrow{a.s.} \RS(N)$ as $n \rightarrow \infty$
for $N \in \left\{\NS,\NAS,\NPE^r,\NCSt\right\}$.
}

\subsection{Expected Measure of $\G_1$-Regions}
Let $\lambda(\cdot)$ be the $\R^d$-Lebesgue measure on $\R^d$ with $d \ge 1$.
In $\R$, $\lambda(\cdot)$ is the length $|\cdot|$,
in $\R^2$, $\lambda(\cdot)$ is the area $A(\cdot)$,
and in $\R^d$, $\lambda(\cdot)$ is the $d$-dimensional volume $V(\cdot)$.
In $\R$ with $\Y_2=\{0,1\}$,
let $\X_n$ be a random sample from $\U(0,1)$, and
$\NS(x)=B(x,r(x))$ where $r(x)=\min(x,1-x)$.
Then,
$$\G_1\left(\X_n,\NS\right)=\left[ X_{n:n}/2,(1+X_{1:n})/2 \right] \Longrightarrow
\lambda(\G_1\left(\X_n,\NS\right))=|\G_1\left(\X_n,\NS\right)|=\left(1+X_{1:n}-X_{n:n}\right)/2.$$
 Hence the expected length of the $\G_1$-region is
\begin{eqnarray*}
\E[\lambda(\G_1\left(\X_n,\NS\right)]&=&\E\left[\frac{1+X_{1:n}-X_{n:n}}{2}\right]=\frac{1+\E[X_{1:n}]-\E[X_{n:n}]}{2}\\
&=&\frac{1+\frac{1}{n+1}-\frac{n}{n+1}}{2}=\frac{1}{n+1} \rightarrow 0 \text{ as }n \rightarrow \infty.
\end{eqnarray*}
In $\R^2$, with three non-collinear points $\Y_3=\{\y_1,\y_2,\y_3\}$,
let $\X_n$ be a random sample from $\UT$.
Then $A(\G_1\left(\X_n,\NPE^r,M\right))>0$ a.s. for all $n<\infty$,
$r>3/2$, $M \in \R^2 \setminus \Y_3$.
Furthermore, $A(\G_1\left(\X_n,\NCS^{\tau=1},M\right))>0$ a.s. for all $n<\infty$ and
$M \in (\TY)^o$;
and for $N \in \left\{\NS,\NAS\right\}$
$A(\G_1\left(\X_n,N,M\right))>0$ a.s. for all $n<\infty$ and $M \in \R^2 \setminus \Y_3$.

The $\G_1$-region, $\G_1\left(\X_n,N\right)$, is closely related to the distribution of the
domination number of the PCD associated with $N$.
Hence we study the asymptotic behavior of the expected area
$\E[A(\G_1\left(\X_n,N,M\right)]$, as $n \rightarrow \infty$,
for $N \in \left\{\NAS,\NPE^r,\NCSt\right\}$.

In $\R$, $\NS$ and $\NAS$ are equivalent functions,
with the extension that $\NAS$ is defined as $\NS$
for $\X$ points outside $(\y_{1:m},\y_{m:m})$.
In higher dimensions, determining the areas of $\G_1\left(\X_n,N,M\right)$
for for $N \in \left\{\NS,\NAS\right\}$ and for general $\X_n$ and
hence finding its expected areas are both open problems.

\subsubsection{The Limit of Expected Area of $\G_1\left(\X_n,N,M\right)$ for $\NPE^r$ and $\NCSt$}
Recall that for $N \in \{\NPE^r,\,\NCSt \}$,
$\G_1\left(\X_n,N,M\right)$ is determined by the (closest) edge extrema
$X_{e}(n) \in \argmin_{X \in \X_n}d(X,e)$, for $e \in \{e_1,e_2,e_3\}$.
So, to find the expected area of $\G_1\left(\X_n,N,M\right)$,
we need to find the expected locus of $X_{e}(n)$;
i.e., the expected distance of $d(X_e(n)$ from $e$.
For example, for $\X_n$ a random sample from a continuous distribution $F$,
$\argmin_{X \in \X_n}d(X,e)$ is unique a.s.,
and if $d(X_e(n),e)=u$, then $X_e(n)$ falls
on a line parallel to $e$ whose distance from $e$ is $u$ a.s.

\begin{lemma}
\label{lem:EXe}
Let $D_i(n):=d(X_{e_i}(n),e_i)$ for $i \in \{1,2,3\}$
and $\X_n$ be a random sample from $\UT$.
Then $\E[D_i(n)] \rightarrow 0$
(i.e., the expected locus of $X_{e_i}(n)$ is on $e_i$)
for each $i \in \{1,2,3\}$, as $n\rightarrow \infty$.
\end{lemma}
\noindent \textbf{Proof:}
Given $Z_i =(X_i,Y_i) \stackrel{iid}{\sim}\UT$.
Then for $e=e_3$, $D_3(n)=Y_{1:n}$ (the minimum $y$-coordinate of $Z_i \in \X_n$).
First observe that $\displaystyle P(Y_i \le y)={\frac {y\,(2\,c_2-y)}{c_2^2}}$,
hence
$$F_Y(y)={\frac {y\,(2\,c_2-y)}{c_2^2}}\I(0 \le y<c_2)+\I(y \ge c_2).$$
So the pdf of $Y_i$ is $\displaystyle f_Y(y)=2\,{\frac {c_2-y}{c_2^2}}\I(0 \le y \le c_2)$.
Then the pdf of $Y_{1:n}$ is
$$f_{1:n}(y)=2\,n(c_2-y)\left(1-{\frac {y\,(2\,c_2-y)}{c_2^2}}\right)^{n-1}{c_2}^{-2}\I(0 \le y \le c_2).$$
Therefore,
$$\E\,[Y_{1:n}]=\int_0^{c_2} \! 2\,y\,n(c_2-y)
\left(1-{\frac {y\,(2\,c_2-y)}{c_2^2}}\right)^{n-1}{c_2}^{-2}{dy}=
\frac{c_2}{2\,n+1} \rightarrow 0, \text{ as }n \rightarrow \infty.$$
 Hence $\E[Y_{1:n}]=\E[D_3(n)] \rightarrow 0$.
Similarly, $\E[D_i(n)] \rightarrow 0$ for $i \in \{1,2\},$ as $n \rightarrow \infty$.
$\blacksquare$


\begin{theorem}
Let $\X_n$ be a random sample from $\UT$ and $M \in \TY^o$.
For $N \in \{\NPE^r,\, \NCSt\}$,
$\E[A(\G_1\left(\X_n,N,M\right))] \rightarrow A(\RS(N,M))$ as $n \rightarrow \infty$.
\end{theorem}
\noindent \textbf{Proof:}
Recall that for $N \in \{\NPE^r,\,\NCSt\}$,
$\G_1\left(\X_n,N,M\right) = \bigcap_{i=1}^3\G_1\left(X_{e_i}(n),N,M\right)$.
Moreover, $\G_1\left(\X_n,N,M\right)=\RS(N,M)$ iff $X_{e_i}(n)\in e_i$
for $i \in \{1,2,3\}$.
In Lemma \ref{lem:EXe},
we have shown that expected locus of $X_{e}(n)$
converges to edge $e$ as $n \rightarrow \infty$.
Hence the expected locus of $\partial (\G_1\left(\X_n,N,M\right)) \cap R_M(e_i)$
converges to the $\partial(\RS(N,M)) \cap R_M(e_i)$ for each $i \in \{1,2,3\}$.
Hence
$$\E[A(\G_1\left(\X_n, N,M\right)] \rightarrow A(\RS(N,M)) \text{ as }n\rightarrow \infty.
\;\blacksquare $$


\begin{remark}
In particular,
\begin{itemize}
\item[i-]
$\E\left[A\left(\G_1\left(\X_n,\NPE^2,M_C\right)\right)\right] \rightarrow 1/4$
as $n \rightarrow \infty$,
since $ \RS\left(\NPE^2,M_C\right)=T(M_1,M_2,M_3)$.
\item[ii-]
$\E\left[ A\left(\G_1 \left(\X_n,\NPE^{r},M\right)\right)\right] \rightarrow 0$
as $n \rightarrow \infty$ if $M \in \Tr$,
since $ \RS\left(\NPE^{r},M\right) = \emptyset$ for $M \in \Tr$.
\item[iii-]
Furthermore, $\E\left[ A \left(\G_1 \left(\X_n,\NPE^{3/2},M_C\right)\right) \right] \rightarrow 0$
since $\RS\left(\NPE^{3/2},M_C\right)=\{M_C\}$.
\item[iv-]
For any $M \in \TY^o$,
$\E\left[ A\left(\G_1\left(\X_n,\NCSt,M)\right)\right)\right] \rightarrow 0$ as
$n \rightarrow \infty$, since $\RS\left(N^{\tau}_{CS},M\right)=\{M\}$.
\item[v-]
We also have $\E\left[ A\left(\G_1\left(\X_n,\NPE^r,M_C\right)\right) \right] \rightarrow 0$
for $r \in [1,3/2]$ as $n \rightarrow \infty$.
\item[vi-]
Furthermore, by careful geometric calculations, we get
$\E\left[ A\left(\G_1\left(\X_n,\NPE^r,M_C\right)\right) \right] \rightarrow \sqrt{3}\left[1-3/(2\,r)\right]^2$
for $ r \in (3/2,2]$.
\item[vii-]
$\E\left[ A\left(\G_1\left(\X_n,\NPE^r,M_C\right)\right) \right] \rightarrow \sqrt{3}/\left[4\,(1-3/r^2)\right]$ for
$r \in (2,\infty]$, as $n \rightarrow \infty$.
\end{itemize}
\end{remark}

We also derive the rate of convergence of $\E[A(\G_1\left(\X_n,\NPE^r,M_C\right))]$ for $r=3/2$.
\begin{theorem}
\label{thm:rate-of-EANYr->0}
Let $\X_n$ be a random sample from $\UT$.
For $r=3/2$, the expected area of the the $\G_1$-region,
$\E\left[ A\left(\G_1\left(\X_n,\NPE^r,M_C\right)\right) \right]$, converges to zero,
at rate $O\left(n^{-2}\right)$.
\end{theorem}

\begin{figure} [ht]
\centering
\scalebox{.5}{\input{Tri_Gamma.pstex_t}}
\caption{
\label{fig:TriGamma1}
The shaded regions are the triangular
$\G_1^{3/2}(X_{e_i}) \cap R_{CM}(\y_i)$ regions for $i=1,2,3$. }
\end{figure}
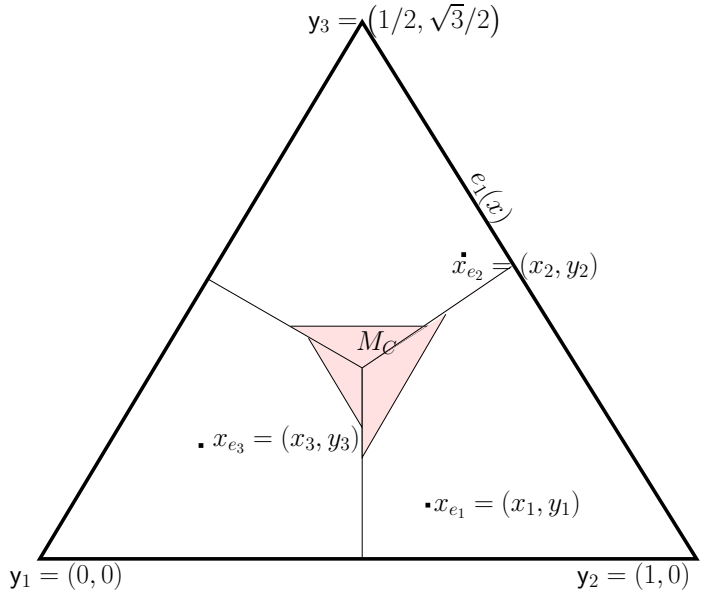

\noindent \textbf{Proof:} For $r=3/2$ and $M=M_C$, and sufficiently large
$n$, $\G_1\left(X_{e_i},\NPE^r\right) \cap R_{CM}(\y_i)$ is a
triangle for $i=1,2,3$ w.p. 1. See Figure \ref{fig:TriGamma1}. With
the realization of the edge extrema denoted as $x_{e_i}=(x_i,y_i)$
close enough to $e_i$, for $i=1,2,3$,

$\G_1\left(x_{e_1},\NPE^{3/2}\right) \cap R_{CM}(\y_1)$ is the
triangle with vertices
$$\left(\frac{\sqrt{3}}{3}\,y_2+x_2-\frac{1}{2},-\frac{\sqrt{3}}{18}\,
\left(-9+2\,\sqrt{3}\,y_2+6\,x_2\right)\right),
\left(\frac{1}{2},\frac{\sqrt{3}}{9}\left(-\frac{9}{2}+2\,\sqrt{3}\,y_2+6\,x_2\right)\right),
\left(\frac{1}{2},\frac{\sqrt{3}}{6}\right),$$
$\G_1\left(x_{e_2},\NPE^{3/2}\right) \cap R_{CM}(\y_2)$ is the triangle with vertices
$$\left(\frac{1}{2},\frac{\sqrt{3}}{18}\,\left(3+4\,\sqrt{3}\,y_3-12\,x_3\right)\right),
\left(\frac{1}{2}-\frac{\sqrt{3}}{3}\,y_3+x_3,-
\frac{\sqrt{3}}{18}\left(-3+2\,\sqrt{3}\,y_3-6\,x_3\right)\right),
\left(\frac{1}{2},\frac{\sqrt{3}}{6}\right),$$
 and $\G_1\left(x_{e_3},\NPE^{3/2}\right) \cap R_{CM}(\y_3)$ is the triangle with vertices
$$\left(-\frac{\sqrt{3}}{6}\,\left(-\sqrt{3}+4\,y_1\right),
\frac{\sqrt{3}}{6}+\frac{2}{3}\,y_1\right),\left(\frac{1}{2},
\frac{\sqrt{3}}{6}\right),\left(\frac{\sqrt{3}}{6}\,\left(\sqrt{3}+4\,y_1\right),
\frac{\sqrt{3}}{6}+\frac{2}{3}\,y_1\right).$$

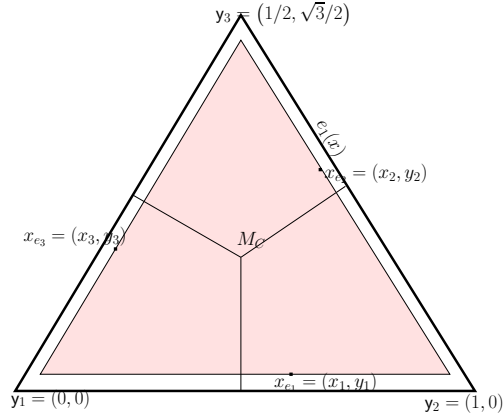
\begin{figure} [ht]
\centering
\scalebox{.35}{\input{jntpdf1.pstex_t}}
\caption{
\label{fig:jntpdf1}
The figure for the joint pdf of $X_{e_i}$. The shaded region is $T(\zeta)$.}
\end{figure}

Then for sufficiently large $n$,
\begin{align*}
&A\left(\G_1\left(\X_n,\NPE^{3/2},M\right)\right) =
\frac{\sqrt{3}}{27}\left(3\,x_2-3+\sqrt{3}\,y_2\right)^2 +
\frac{\sqrt{3}}{27}\left(-3\,x_3+\sqrt{3}\,y_3\right)^2 + \frac{4\,\sqrt{3}}{9}y_1^2\\
       &= \frac{\sqrt{3}}{9}\left(3\,x_2^2-6\,x_2+2\,\sqrt{3}\,y_2\,x_2-
       2\,\sqrt{3}\,y_2+y_2^2+3+y_3^2-2\,\sqrt{3}\,y_3\,x_3+3\,x_3^2+4\,y_1^2\right).
\end{align*}
To find the expected area, we need the joint density of the $X_{e_i}$.
The edge extrema are all distinct with probability 1 as
$n \rightarrow \infty$ (see Theorem \ref{thm:distinct-edge-ext}).
Let $T(\zeta)$ be the triangle formed by the lines at $x_{e_i}$ parallel
to $e_i$ for $i=1,2,3$ where $\zeta=(x_1,y_1,x_2,y_2,x_3,y_3)$.
See Figure \ref{fig:jntpdf1}.

Then the asymptotically accurate joint pdf of $X_{e_i}$ (\cite{ceyhan:dom-num-NPE-MASA})
is
\begin{eqnarray*}
f_3(\zeta) &=& n(n-1)(n-2)\left(\frac{A(T(\zeta))}{A(\TY)}\right)^{n-3}\frac{1}{A(\TY)^3}\\
       &=& n(n-1)(n-2) \left(\sqrt{3}/36\left(-2\,\sqrt{3}\,y_1+
       \sqrt{3}\,y_3-3\,x_3+\sqrt{3}\,y_2+3\,x_2\right)^2\right)^{n-3}\Big/\left(\sqrt{3}/4\right)^n.
\end{eqnarray*}
with the support $D_S=\bigl\{\zeta \in \Re^6:\; (x_i,y_i) \text{'s are distinct} \bigr\}.$

Then for sufficiently large $n$,
\begin{eqnarray*}
\E\left[ A\left(\G_1\left(\X_n,\NPE^{3/2},M\right)\right) \right] &\sim&
\int A\left(\G_1 \left(\X_n,\NPE^{3/2}\right)\right)f_3(\zeta)d\zeta\\
 & = &\int A\left(\G_1 \left(\X_n,\NPE^{3/2}\right)\right)n(n-1)(n-2)
 \left(\frac{A(T(\zeta))}{A(\TY)}\right)^{n-3}\frac{1}{A(\TY)^3}d\zeta.
\end{eqnarray*}
Let $G(\zeta)=A(T(\zeta))/A(\TY)$.
Notice that the integrand is
critical when $x_{e_i} \in e_i$ for $i=1,2,3$, since $G(\zeta)$ when
$x_{e_i} \in e_i$ for each $i=1,2,3$.
So we make the change of
variables  $y_1=z_1$, $y_2=\sqrt{3}\,(1-x_2)-z_2$, and
$y_3=\sqrt{3}\,x_3-z_3$, then $G(\zeta)$ and $A(\G_1\left(\X_n,N\right))$ becomes
$$G(\zeta)=G(z_1,z_2,z_3)=\left(2\,z_1+z_3-\sqrt{3}+z_2\right)^2\Big/3\text{ and }
A(\G_1\left(\X_n,N\right))=\sqrt{3}\,\left(z_2^2+z_3^2+4\,z_1^2\right)/9,$$
respectively.
Hence the integrand does not depend on $x_1,x_2,x_3$
and integrating with respect to $x_1,x_2,x_3$ yields a constant $K$.
Now, the integrand is critical at $(z_1,z_2,z_3)=(0,0,0)$, since $G(0,0,0)=1$.
So let $E^{\ve}_{u}$ be the event that $0 \le
z_i \le \ve$ for $i=1,2,3$ for $\ve>0$ small enough.
Then making the change of variables $z_i =w_i/n$ for $i=1,2,3$, we
get $A\left(\G_1 \left(\X_n,\NPE^{3/2}\right)\right)=O\left(n^{-2}\right)$ and $G(z_1,z_2,z_3)$ becomes
$G(w_1,w_2,w_3)=1-\frac{1}{n}\left(2/\sqrt{3}\,(2\,w_1+w_2+w_3)\right)+O\left(n^{-2}\right)$, hence
\begin{align*}
& \E\left[ A\left(\G_1\left(\X_n,\NPE^{3/2},M\right)\right) \right] \sim\\
& K \int_0^{n\ve}\int_0^{n\ve}\int_0^{n\ve}
A\left(\G_1 \left(\X_n,\NPE^{3/2}\right)\right)n(n-1)(n-2)\frac{1}{n^3}\,
G(w_1,w_2,w_3)^{n-3} dw_1dw_2dw_3,\\
&  \text{  letting $n \rightarrow \infty$}\\
& \approx K\int_0^{\infty}\int_0^{\infty}\int_0^{\infty}
O\left(n^{-2}\right)\exp\left(-2/\sqrt{3}\,(2\,w_1+w_2+w_3)\right)dw_1dw_2dw_3 =
O\left(n^{-2}\right),
\end{align*}
since $\int_0^{\infty}\int_0^{\infty}\int_0^{\infty}
\exp\left(-2/\sqrt{3}\,(2\,w_1+w_2+w_3)\right)dw_1dw_2dw_3=
3\,\sqrt{3}/16$ which is a finite constant.
Hence $\E\left[ A\left(\G_1\left(\X_n,\NPE^{3/2},M\right)\right) \right] \rightarrow 0$
as $n \rightarrow \infty$ at the rate $O\left(n^{-2}\right)$.
$\blacksquare$


\section{$\G_k$-Regions for Proximity Maps in $\TY$}
\label{sec:Gammak-regions}
We can also define the regions associated with $\g_n(N)=k$ for $k
\le n$.

In $\R$ with $\Y_2=\{0,1\}$,  $\g_n(N) \le 2$, hence we can only
define $\G_2$-regions.  Recall that $\g_n(N)=2$ iff $\X_n \cap
\left[ \frac{x_{n:n}}{2},\frac{1+x_{1:n}}{2} \right] = \emptyset$
iff $\X_n \subset [0,1] \setminus \G_1\left(\X_n,N\right)$.  So
\begin{eqnarray*}
\G_2\left(\X_n,N\right)&=&\left\{ (x,y) \in [0,1]^2:\;
\X_n \subset N(x) \cup N(y);\;x,y \notin \G_1\left(\X_n,N\right) \right\}\\
&=&\left\{ (x,y) \in [0,1]^2 \setminus \G_1\left(\X_n,N\right)^2:\; \X_n
\subset N(x) \cup N(y) \right\}.
\end{eqnarray*}
Notice that $\G_2\left(\X_n,N\right) \subseteq [0,1]^2$.
Let
$$\widetilde{x}_1:=\argmin_{x \in \X_n \cap (0,1/2)} (1/2-x)\text{ and }
\widetilde{x}_2:=\argmin_{x \in \X_n \cap (1/2,1)} (x-1/2),$$
then $\g_n(N)=2$ iff $\widetilde{x}_1, \widetilde{x}_2 \notin \G_1\left(\X_n,N\right)$.
In such a case $\X_n \subset N(\widetilde{x}_1) \cup N(\widetilde{x}_2)$ by construction.

In general,
\begin{definition}
The \emph{$\G_2$-region} for proximity map $\NY(\cdot)$ and set
$B \subset \Omega$ is $\G_2(B,N)=\{(x,y) \in [\Omega \setminus \G_1(B)]^2:\;B \subseteq \NY(x)\cup \NY(y)\}$.
In general, \emph{$\G_k$-region} for proximity map $\NY(\cdot)$ and
set $B \subset \Omega$ for $k=1,2,\ldots,n$ is
\begin{multline*}
\G_k(B,N)=\Biggl\{ (x_1,x_2,\ldots,x_k) \in \Omega^k: B \subseteq
\bigcup_{i=1}^k \NY\left(x_i\right) \;\text{ and all possible $m$-permutations $(u_1,u_2,\ldots,u_m)$}\\
\text{ of $(x_1,x_2,\ldots,x_k)$
satisfy }(u_1,u_2,\ldots,u_m) \not\in \G_m(B,N) \text{ for each } m=1,2,\ldots,k-1 \Biggr\}.
\square
\end{multline*}
\end{definition}
Note that $\G_k$-regions are defined for $k \le n$
and they might be empty.
Moreover, $\G_k$-regions are in $\Omega^k$, not in $\Omega$.

Let $\displaystyle \left \{ {n \atop m} \right \}$ denote the Stirling partition
number for a set of size $n$ into $m$ blocks and let $\left \{ {A
\atop m} \right \}$ denote all Stirling partitions of a set $A$ into
$m$ blocks; that is,
$$\left\{{A \atop m}\right\}:=\left\{ \{\mathcal B_1,\mathcal B_2, \ldots, \mathcal B_m\}:\;
\mathcal B_i \not= \emptyset, \mathcal B_i \subset A;\; A=\bigcup_i \mathcal B_i \right\}.$$

In particular, $\left \{ { \X_n \atop 2} \right \}$ is the unordered
pair of blocks $\mathcal B_1$ and $\mathcal B_2$ such that
$\mathcal B_i \not= \emptyset$ and $\mathcal B_i \subset \X_n$ for $i=1,2$,
and $\mathcal B_1 \cup \mathcal B_2 =\X_n$.
Note that $\mathcal B_2 = \X_n \setminus \mathcal B_1$.
Then
\begin{proposition}
$\G_2\left(\X_n, N\right)=
\bigcup_{\{ \mathcal B_1,\mathcal B_2\} \in \left\{ {\X_n \atop 2} \right\}}
[\G_1\left(\mathcal B_1,N\right) \times \G_1\left(\mathcal B_2,N\right)] \setminus \G_1\left(\X_n,N\right)^2$
for any nonempty Stirling blocks $\mathcal B_1$ and $\mathcal B_2$ in $\left \{ { \X_n \atop 2} \right \}$.
\end{proposition}
\noindent \textbf{Proof:}
Given $\X_n$, suppose $(u,v) \in \G_2\left(\X_n,N\right)$, then
$\X_n \subset N(u) \cup N(v)$ and $u,v \notin \G_1\left(\X_n,N\right)$.
Let $\mathcal B_1=\X_n \cap N(u)$, and $\mathcal B_2=[\X_n \cap N(v)] \setminus N(u)$.
Then $\mathcal B_1$ and $\mathcal B_2$ are two
Stirling blocks in $\left\{ {\X_n \atop 2} \right\}$.
Hence $\mathcal B_1 \subset N(u) \Rightarrow u \in \G_1\left(\mathcal B_1,N\right)
\setminus \G_1\left(\X_n,N\right)$ and $\mathcal B_2 \subset N(v) \Rightarrow v
\in \G_1\left(\mathcal B_2,N\right) \setminus \G_1\left(\X_n,N\right)$, hence $\G_2\left(\X_n, N\right)
\subseteq \bigcup_{\{ \mathcal B_1,\mathcal B_2\} \in \left\{ {\X_n
\atop 2} \right\}} [\G_1\left(\mathcal B_1,N\right) \times \G_1\left(\mathcal B_2,N\right)]
\setminus (\G_1\left(\X_n,N\right))^2$.
The other direction is trivial, hence the desired result follows.
$\blacksquare$

In $\R$ with $N=\NS$, we can exploit the natural ordering available.

\begin{proposition} In $\R$ with $\Y_2=\{0,1\}$,
$\displaystyle \G_2\left(\X_n,\NS\right)=\bigcup_{k=1}^{n-1} \left(\frac{x_{k:n}}{2},\frac{x_{n:n}}{2}\right) \times
\left(\frac{1+x_{1:n}}{2},\frac{1+x_{(k+1):n}}{2}\right)$.
\end{proposition}
\noindent \textbf{Proof:}
Recall that $\G_2\left(\X_n,\NS\right)=
\bigcup_{\{ \mathcal B_1,\mathcal B_2\} \in \left\{ {\X_n \atop 2} \right\}}
[\G_1\left(\mathcal B_1,\NS\right) \times \G_1\left(\mathcal B_2,\NS\right)] \setminus \G_1\left(\X_n,\NS\right)^2$.
Let $\{\mathcal B_1,\mathcal B_2\}$ be a Stirling
partition in $\left\{ { \X_n \atop 2}\right\}$.
First observe that
$\mathcal B_1=\{x_{1:n},\ldots,x_{k:n} \}$ and
$\mathcal B_2=\{x_{(k+1):n}, \ldots,x_{n:n} \}$ forms a Stirling partition of
$\X_n $, and
$$\G_1 \left(\{x_{1:n},\ldots,x_{k:n} \},\NS\right) \setminus
\G_1\left(\X_n,\NS\right)=\left[ \frac{x_{k:n}}{2},\frac{1+x_{1:n}}{2} \right] \setminus
\left[ \frac{x_{n:n}}{2},\frac{1+x_{1:n}}{2} \right]=
\Biggl[\frac{x_{k:n}}{2},\frac{x_{n:n}}{2}\Biggr)$$
and
$$\G_1\left(\{x_{(k+1)}, \ldots,x_{n:n} \},\NS\right) \setminus
\G_1\left(\X_n,\NS\right)=\left[ \frac{x_{n:n}}{2},\frac{1+x_{(k+1)}}{2} \right] \setminus
\left[ \frac{x_{n:n}}{2},\frac{1+x_{1:n}}{2} \right]=
\Biggl(\frac{1+x_{1:n}}{2},\frac{1+x_{(k+1)}}{2}\Biggr].$$
Furthermore,
\begin{equation}
\label{eqn:G1B1B2}
\G_1\left(\{x_{1:n}, \ldots,x_{k:n} \},\NS\right)
\times \G_1\left(\{x_{(k+1)}, \ldots,x_{n:n} \},\NS\right) \supseteq
\bigcup_{ \{\mathcal B_1,\mathcal B_2\} \in \Lambda}\G_1\left(\mathcal B_1,\NS\right)
\times \G_1\left(\mathcal B_2,\NS\right)
\end{equation}
where
$$\Lambda =\left\{ \{\mathcal B_1,\mathcal B_2\} \in \left\{ {\X_n \atop 2} \right\}:\;
x_{1:n} \in \mathcal B_1,\; x_{n:n} \in \mathcal B_2,\;
\max\left(\mathcal B_1\right)=x_{k:n},\; \min\left(\mathcal B_2\right)=x_{(l)},\; l <k \right\},$$
since the left hand side in Equation \eqref{eqn:G1B1B2} is
$\displaystyle \left[ \frac{x_{k:n}}{2},\frac{x_{n:n}}{2} \right] \times
\left[ \frac{1+x_{1:n}}{2},\frac{1+x_{(k-1)}}{2} \right] $.
Hence the desired result follows.
$\blacksquare$


Similarly, for $\NPE^r$,
$$\G_2\left(\X_n,\NPE^r\right)=
\bigcup_{\{ \mathcal B_1,\mathcal B_2\} \in \left\{ {\X_n \atop 2} \right\}}
\left[ \G_1\left(\mathcal B_1,\NPE^r,M\right) \times
\G_1\left(\mathcal B_2,\NPE^r,M\right) \right] \setminus \G_1\left(\X_n,\NPE^r,M\right)^2.$$
Note that $\G_1\left(\mathcal B_i,\NPE^r,M\right)$ is determined by the edge extrema in
$\mathcal B_i$, $i=1,2$.
Furthermore, if $(u,v) \in \G_2\left(\X_n,\NPE^r,M\right)$, then $(u,v)\notin R_M(\y)^2$,
since either $\NPE^r(u) \subseteq \NPE^r(v)$ or $\NPE^r(v) \subseteq \NPE^r(u)$
should hold if $(u,v)\in R_M(\y)^2$.

For $\NCSt$,
$$\G_2\left(\X_n,\NCSt,M\right)=\bigcup_{\{ \mathcal B_1,\mathcal B_2\} \in
\left\{ {\X_n \atop 2} \right\}} \left[ \G_1\left(\mathcal B_1,\NCSt,M\right)
\times \G_1\left(\mathcal B_2,\NCSt,M\right) \right] \setminus
\G_1\left(\X_n,\NCSt,M\right)^2.$$
Note that $\G_1\left(\mathcal B_i,\NCSt,M\right)$ is determined by the
edge extrema of $\mathcal B_i$, $i=1,2$.
But if $(u,v) \in \G_2\left(\X_n,\NCSt,M\right)$,
then $(u,v) \in R_M(e)^2$ can hold.
In these proximity regions,
not all edge extrema should fall in the same partition set
for $\G_k\left(\X_n,\NCSt,M\right)$ to be nonempty.

For any proximity map $N$,
\begin{align*}
&P(\g_n(N)=2)=P(\X_n^2 \cap \G_2\left(\X_n, N\right) \not=\emptyset,\,\g_n(N) \not= 1)\\
&=P\left(\X_n^2 \bigcap \left[ \bigcup_{\{ \mathcal B_1,\mathcal B_2\} \in
\left\{ {\X_n \atop 2} \right\}} \left[ \G_1\left(\mathcal B_1,N\right) \times
\G_1\left(\mathcal B_2,N\right) \right] \setminus (\G_1\left(\X_n,N\right))^2 \right]
\not= \emptyset ,\,\X_n \cap \G_1\left(\X_n, N\right) =\emptyset\right).
\end{align*}
A more compact way to write this is as
$$P(\g_n(N)> 2)=P\left(\X_n^2 \cap \G_{\le 2}\left(\X_n, N\right) =\emptyset\right)$$
where $\G_{\le 2}\left(\X_n, N\right):=\bigcup_{\{ \mathcal B_1,\mathcal B_2\} \in
\left\{ {\X_n \atop 2} \right\}} \G_1\left(\mathcal B_1,N\right) \times
\G_1\left(\mathcal B_2,N\right)$.

Furthermore, for $k \ge 3$, the $\G_k$-regions are defined similarly
as
$$\G_k\left(\X_n, N\right)=\bigcup_{\{ \mathcal B_1,\ldots, \mathcal B_k\} \in
\left\{ { \X_n \atop k} \right \}} [\G_1\left(\mathcal B_1,N\right) \times
\ldots \times \G_1\left(\mathcal B_k,N\right)] \setminus \G_1\left(\X_n,N\right)^k.$$
Hence,
\begin{eqnarray*}
P(\g_n(N)=k)&=&P\left(\X_n^k \cap \G_k\left(\X_n, N\right) \not=\emptyset,\g_n(N)>k-1\right)\\
&=&P\left(\X_n^k \cap \left[ \bigcup_{\{ \mathcal B_1,\ldots, \mathcal B_k\} \in
\left\{ { \X_n \atop k} \right \}} \G_1\left(\mathcal B_1,N\right)
\ldots \times \G_1\left(\mathcal B_k,N\right) \setminus \G_1\left(\X_n,N\right)^k \right]
\not= \emptyset\right).
\end{eqnarray*}
A more compact way to write this is as
$$P(\g_n(N)> k)=P\left(\X_n^k \cap \G_{\le k}\left(\X_n, N\right) \not=\emptyset\right)$$
where $\G_{\le k}\left(\X_n, N\right):=\bigcup_{\{ \mathcal B_1, \ldots,\mathcal
B_k\} \in \left\{ { \X_n \atop k} \right \}} \G_1\left(\mathcal B_1,N\right)
\times \ldots \times \G_1\left(\mathcal B_k,N\right)$.

\section{$\kappa$-Values for the Proximity Maps in $\TY$}
\label{sec:kappa-value}
Recall that the domination number, $\g_n(N)$ is the cardinality of a
minimum dominating set of the PCD based on $N$.
So by definition, $\g_n(N) \le n$.
We will seek an a.s. least upper bound for $\g_n(N)$
which suggests the following concept.

\begin{definition}
Let $\X_n$ be a random sample from $F$ on $\TY$
and let $\g_n(N)$ be the domination number for the
PCD based on a proximity map $N$.
The general a.s. least upper bound for $\g_n(N)$ that works for all
$n \ge 1$ is called the \emph{$\kappa$-value}; i.e.,
$\kappa_n(N) := \min \{k(n): \g\left(\X_n, N\right) \le k(n) \text{ a.s. for all } n \ge 1\}$.
$\square$
\end{definition}
It is more desirable to have a $\kappa$-value that is independent of $n$.
Further, if $\kappa_n(N)=\kappa$ exists for $N$ and is independent of $n$,
then the domination number has the following discrete probability mass function:
\begin{equation*}
\label{eqn:}
\g\left(\X_n, N\right) =
\begin{cases}
1   &\text{w.p.} \quad p_1, \\
2   &\text{w.p.} \quad p_2, \\
\vdots   & \text{ $~\vdots$} \quad \quad \vdots, \\
\kappa   &\text{w.p.} \quad p_\kappa.
\end{cases}
\end{equation*}

In $\R$ with $\Y_2=\{0,1\}$, for $\X_n$ a random sample from $\U(0,1)$,
$\g_n(N) \le 2$ with equality holding with positive probability
for $N \in \left\{\NS,\NAS \right\}$.
Hence $\kappa_n(\NS)=2$.
But in $\R^d$ with $d>1$, finding $\kappa_n(N)$ for
$N \in \left\{\NS,\NAS \right\}$ is an open problem.
Next, we investigate the $\kappa$-values for
$\left\{\NAS,\NPE^r,\NCSt\right\}$ in $\R^2$.

\begin{theorem}
\label{thm:kappaNYr=3}
Let $\X_n$ be a random sample from $\UT$, and $M \in \R^2 \setminus \Y_3$.
For $\NPE^r(\cdot,M)$, $\kappa_n\left(\NPE^r\right)=3$.
\end{theorem}
\noindent \textbf{Proof:}
For $\NPE^r(\cdot,M)$, pick the point closest to edge $e_i$
in vertex region $R_M(\y_i)$;
that is, pick
$U_i \in \argmin_{X \in \X_n \cap R_M(\y_i)} d(X,e_i)
=\argmax_{X \in \X_n \cap R_M(\y_i)} d(\ell(\y,X),\y_i)$
in the vertex region for which $\X_n \cap R_M(\y_i) \not=\emptyset$ for $i \in \{1,2,3\}$
(note that as $n \rightarrow \infty$,
$\X_n \cap R_M(\y_i) \not=\emptyset$ for all $i \in \{1,2,3\}$ a.s.,
and also $U_i$ is unique a.s. for each $i$, since $X$ is from $\UT$).
Then $\X_n \cap R_M(\y_i) \subset \NPE^r(U_i,M)$.
Hence $\X_n \subset \bigcup_{i=1}^3 \NPE^r(U_i,M)$.
 So $\g_n\left(\NPE^r,M_C\right) \le 3$ with equality holding with positive probability.
Thus $\kappa_n\left(\NPE^r\right)=3$.
$\blacksquare$

There is no least upper bound for $\g_n(\NCSt,M)$
that works for all $n>0$ as shown below.

\begin{figure}
\begin{center}
\psfrag{A}{\scriptsize{$\y_1$}}
\psfrag{B}{\scriptsize{$\y_2$}}
\psfrag{C}{\scriptsize{$\y_3$}}
\psfrag{x}{}
\psfrag{M}{\scriptsize{$M$}}
\psfrag{(x1,y1)}{\scriptsize{$(x_1,y_1)$}}
\psfrag{(x2,y2)}{\scriptsize{$(x_2,y_2)$}}
\psfrag{(xn,yn)}{\scriptsize{$(x_n,y_n)$}}
\epsfig{figure=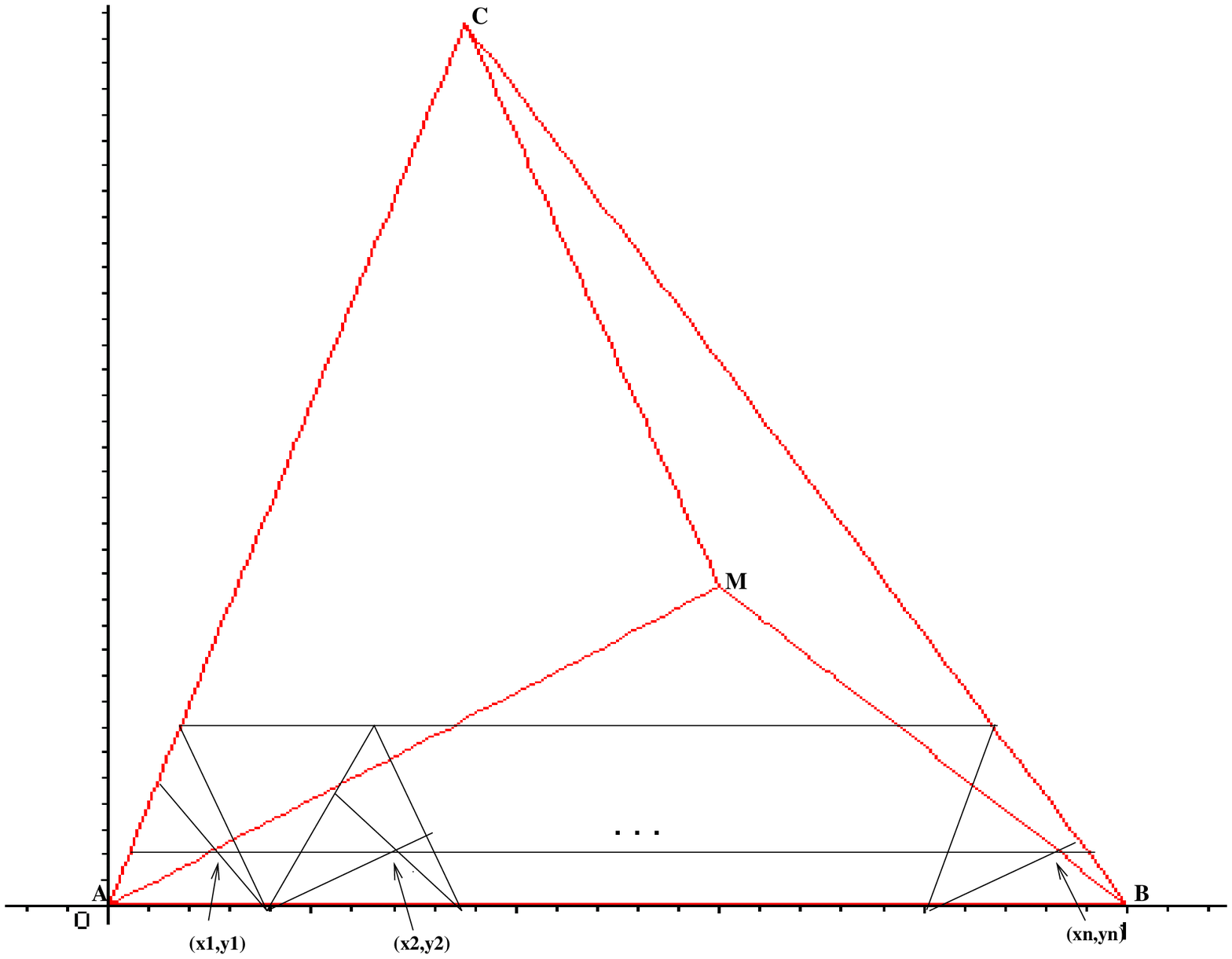, height=210pt,width=300pt}
\end{center}
\caption{
\label{fig:kappa_NCS}
The figure for $\kappa_n\left(\NCSt\right)=n$.}
\end{figure}

\begin{theorem}
\label{thm:kappa-NCS-n}
Let $\X_n$ be random sample from $\UT$.
Then $\g_n(\NCSt,M)=n$ holds with positive
probability for all $\tau \in [0,1]$.
\end{theorem}
\noindent \textbf{Proof:}
For $\tau=0$, the result follows trivially.
For $\tau=1$,
we will prove the theorem by showing that there is a
union of $n$ regions of positive area in $\TY$,
so that $u \in \NCS^{\tau=1}(v,M)$ iff $u=v$, for any $u,v \in \X_n$.
Let $M=(m_1,m_2) \in \TY^o$.
In $R_{M}(e_3)$ locate $n$ triangles evenly
on $e_3$ with base length $1/n$ and
similar to $\TY$ (with similarity ratio $1/n$).
See also Figure \ref{fig:kappa_NCS}.
Then locate $n$ points in each triangle at
$z_i=(x_i,y_i)$ such that $(x_i,y_i)$ is the same type of
center of $T_i$ as $M$ is of $\TY$.
Then using the similarity ratio of
$\NCS^{\tau=1}(z_i,M)$ to $\TY$, namely, $y_i/m_2=1/n$,
we get $y_i=m_2/n$ for all $i=1,2,\ldots,n$.
Moreover, $x_i-x_{i-1}=1/n$ for $i=2,3,\ldots,n$
with $x_1=m_1/n$ and $x_n=1-(1-m_1)/n$.
Then $(x_i,y_i) \in \NCS^{\tau=1}((x_j,y_j),M)$ iff $i=j$.
Furthermore, for sufficiently small $\ve>0$, the same holds for
the $\ve$ neighborhood of each $z_i=(x_i,y_i)$.
That is, $\NCS^{\tau=1}(x,M) \cap B(z_j,\ve)=\emptyset$
for all $x \in B(z_i,\ve)$ for any distinct pair $i,j \in \{1,2,\ldots,n\}$,
and probability of $\X_n$ being composed of $n$ points one from each
$B(z_i,\ve)$ is positive.
Then $\g_n(\NCSt,M)=n$ holds with positive probability.
  The result for $\tau \in (0,1)$ follows similarly.
$\blacksquare$

\subsection{Characterization of Proximity Maps Using $\kappa$-Values}
Notice that for the proximity maps we have considered,
$\kappa_n\left(\NPE^r\right)=3$ and
$\kappa_n\left(\NCSt\right)=n$.
One property of proximity maps that makes $\kappa_n(N)<n$
is that probability of having an $\X_n$ for
which $N(X) \cap \X_n =\{X\}$ for all $X \in \X_n$ is zero.

Below we state a condition for $\kappa_n(N(\cdot,M)) \le 3$
for $N(\cdot,M)$ defined with $M$-vertex regions.
\begin{theorem}
\label{thm:kappa<=3,vertex}
Suppose $N(\cdot,M)$ is defined with $M$-vertex regions
 with $M \in \R^2 \setminus \Y_3$ and
$N(x,M)$ gets larger as $d(\ell(\y,x),\y)$ increases for $x \in R_M(\y)$
in the sense that $N(x,M) \subseteq N(z,M)$
for all $x,z \in R_M(\y)$ when $d(\ell(\y,x),\y) \le d(\ell(\y,z),\y)$.
Then $\kappa_n(N) \le 3$.
\end{theorem}
\noindent \textbf{Proof:}
When $\X_n \cap R_M(\y_i)\not=\emptyset$, pick one of the points
$U_i(n) \in \argmax_{X \in \X_n \cap R_M(\y_i)} d(\ell(\y_i,X),\y_i)$,
then $\X_n \cap R_M(\y_i) \subset N(U_i(n))$ for each $i \in \{1,2,3\}$.
So $\g\left(\X_n,N,M\right) \le 3$, and hence $\kappa_n(N) \le 3$.
$\blacksquare$

Notice that $\NPE^r$ satisfies the conditions of Theorem \ref{thm:kappa<=3,vertex}.

\begin{theorem}
\label{thm:kappa<=3,edge}
Suppose $N(\cdot,M)$ is defined with $M$-edge regions
and $N(x,M)$ gets larger as $d(x,e)$ increases for $x \in R_M(e)$
in the sense that $N(x,M) \subseteq N(y,M)$
for all $x,y \in R_M(e)$ when $d(x,e) \le d(y,e)$.
Then $\kappa_n(N) \le 3$.
\end{theorem}
\noindent \textbf{Proof:}
When $\X_n \cap R_M(e_i) \not=\emptyset$,
pick one of the points
$U_{e_i}(n) \in \argmax_{X \in \X_n \cap R_M(e_i)} d(X,e_i)$.
Then $\X_n \cap R_M(e_i) \subset N(U_e(n))$ for each $i \in \{1,2,3\}$.
So $\g\left(\X_n,N,M\right) \le 3$, and hence $\kappa_n(N) = 3$.
$\blacksquare$

In Theorems \ref{thm:kappa<=3,vertex} and \ref{thm:kappa<=3,edge},
we have an upper bound for $\kappa_n(N)$.
To determine an exact value for $\kappa_n(N)$ we need further restrictions.
Let $A_1:=\{x \in \TY: \mu(N(x,M))=\mu(\TY)\}$ and
$A_2:=\{(x,y) \in \TY \times \TY: \mu(N(x,M)\cup N(y,M))=\mu(\TY)\}$.
If in addition to the hypothesis of Theorem
\label{thm:kappa<=3,vertex} (and \label{thm:kappa<=3,edge})
we have $A_1$ and $A_2$ have zero measure (e.g., zero area for continuous distributions
with support in $\TY$) then $\kappa_n(N)=3$ would hold.

\section{Discussion and Conclusions}
\label{sec:disc-conc}
In this article,
we provide a probabilistic characterization of proximity maps,
associated regions, and digraphs and related quantities.
In particular, we discuss the probabilistic behavior of
proximity regions, superset regions, and $\G_1$-regions;
construct digraphs (proximity catch digraphs(PCDs)) and
investigate related quantities such as domination number, $\eta$ and $kappa$ values.
We also provide auxiliary tools such as vertex and edge
regions for the construction of proximity regions.

Although $\G_1$-regions and superset regions were introduced before
(\cite{ceyhan:dom-num-NPE-SPL},\cite{ceyhan:arc-density-PE}, and \cite{ceyhan:dom-num-NPE-MASA})
a thorough investigation is only performed in this article.
$G_1$-regions are a sort of a ``dual" of proximity regions
and are associated with domination number being equal to 1.
We provide a probabilistic characterization of $\G_1$-regions for
general proximity maps $N$ and data points from distribution $F$.
We also extend this concept by introducing $\G_k$-regions,
which are associated with domination number being equal to $k$.

We introduce the quantities related to $\G_1$-regions and domination
number, namely, $\eta$-values and $\kappa$-values.
$\eta$-value is the minimum number of points in a set required to determine
the $\G_1$-region for that set.
We determine some general conditions that make $\eta_n(N)\le 3$
for data in the triangle $\TY$.
$\kappa$-value is the a.s. least upper bound for the domination number of the PCDs.
We also determine some general conditions that make $\kappa_n(N)\le 3$
for data in the triangle $\TY$.

We provide two PCD families, namely proportional-edge PCDs (\cite{ceyhan:arc-density-PE} and
central similarity PCDs (\cite{ceyhan:arc-density-CS}) as illustrative examples.
We discuss the construction of proximity regions and $\G_1$-regions
for these PCDs.
Furthermore, we calculate the expected area of $\G_1$-regions
for these PCDs in the limit.

Determining $\G_1$-regions, $\eta$ and $\kappa$ values for
spherical and arc-slice PCDs contain many open problems
and are subjects of ongoing research.

With the above characterizations,
given another PCD, then we can determine how it behaves
in terms of $\G_1$-regions and related quantities;
in particular we can determine a.s. least upper bounds
for the domination number of the new PCD.

\section*{Acknowledgments}
This research was supported by the research agency TUBITAK via the Kariyer Project \# 107T647.


\end{document}

%% file: Basic_Tri.pstex_t
\begin{picture}(0,0)%
\includegraphics{Basic_Tri.pstex}%
\end{picture}%
\setlength{\unitlength}{3947sp}%
\begingroup\makeatletter\ifx\SetFigFont\undefined%
\gdef\SetFigFont#1#2#3#4#5{%
  \reset@font\fontsize{#1}{#2pt}%
  \fontfamily{#3}\fontseries{#4}\fontshape{#5}%
  \selectfont}%
\fi\endgroup%
\begin{picture}(11349,8124)(289,-7348)
\put(4276,-661){\makebox(0,0)[lb]{\smash{{\SetFigFont{50}{16.8}{\rmdefault}{\mddefault}{\updefault}{\color[rgb]{0,0,0}$y_3=(c_1,c_2)$}%
}}}}
\put(4126,-3586){\makebox(0,0)[lb]{\smash{{\SetFigFont{50}{16.8}{\rmdefault}{\mddefault}{\updefault}{\color[rgb]{0,0,0}$(x,y)$}%
}}}}
\put(901,-1336){\makebox(0,0)[lb]{\smash{{\SetFigFont{50}{16.8}{\rmdefault}{\mddefault}{\updefault}{\color[rgb]{0,0,0}$y$}%
}}}}
\put(5026,-7261){\makebox(0,0)[lb]{\smash{{\SetFigFont{50}{16.8}{\rmdefault}{\mddefault}{\updefault}{\color[rgb]{0,0,0}$x$}%
}}}}
\put(9301,-6586){\makebox(0,0)[lb]{\smash{{\SetFigFont{50}{16.8}{\rmdefault}{\mddefault}{\updefault}{\color[rgb]{0,0,0}$y_2=(1,0)$}%
}}}}
\put(826,-6586){\makebox(0,0)[lb]{\smash{{\SetFigFont{50}{16.8}{\rmdefault}{\mddefault}{\updefault}{\color[rgb]{0,0,0}$y_1=(0,0)$}%
}}}}
\end{picture}%

%% file: Equal_Tri_phi.pstex_t
\begin{picture}(0,0)%
\includegraphics{Equal_Tri_phi.pstex}%
\end{picture}%
\setlength{\unitlength}{3947sp}%
\begingroup\makeatletter\ifx\SetFigFont\undefined%
\gdef\SetFigFont#1#2#3#4#5{%
  \reset@font\fontsize{#1}{#2pt}%
  \fontfamily{#3}\fontseries{#4}\fontshape{#5}%
  \selectfont}%
\fi\endgroup%
\begin{picture}(11349,8124)(289,-7348)
\put(826,-6586){\makebox(0,0)[lb]{\smash{{\SetFigFont{50}{16.8}{\rmdefault}{\mddefault}{\updefault}{\color[rgb]{0,0,0}$\phi_e(y_1)=(0,0)$}%
}}}}
\put(9301,-6586){\makebox(0,0)[lb]{\smash{{\SetFigFont{50}{16.8}{\rmdefault}{\mddefault}{\updefault}{\color[rgb]{0,0,0}$\phi_e(y_2)=(1,0)$}%
}}}}
\put(901,-1336){\makebox(0,0)[lb]{\smash{{\SetFigFont{50}{16.8}{\rmdefault}{\mddefault}{\updefault}{\color[rgb]{0,0,0}$v$}%
}}}}
\put(5026,-7261){\makebox(0,0)[lb]{\smash{{\SetFigFont{50}{16.8}{\rmdefault}{\mddefault}{\updefault}{\color[rgb]{0,0,0}$u$}%
}}}}
\put(5551,389){\makebox(0,0)[lb]{\smash{{\SetFigFont{50}{16.8}{\rmdefault}{\mddefault}{\updefault}{\color[rgb]{0,0,0}$\phi_e(y_3)$}%
}}}}
\put(4651,-3211){\makebox(0,0)[lb]{\smash{{\SetFigFont{50}{16.8}{\rmdefault}{\mddefault}{\updefault}{\color[rgb]{0,0,0}$\phi_e(x,y)$}%
}}}}
\end{picture}%

%% file: Nofnu2.pstex_t
\begin{picture}(0,0)%
\includegraphics{Nofnu2.pstex}%
\end{picture}%
\setlength{\unitlength}{3947sp}%
\begingroup\makeatletter\ifx\SetFigFont\undefined%
\gdef\SetFigFont#1#2#3#4#5{%
  \reset@font\fontsize{#1}{#2pt}%
  \fontfamily{#3}\fontseries{#4}\fontshape{#5}%
  \selectfont}%
\fi\endgroup%
\begin{picture}(11231,8043)(64,-7348)
\put(1051,-6511){\makebox(0,0)[lb]{\smash{\SetFigFont{22}{14.4}{\rmdefault}{\mddefault}{\updefault}{\color[rgb]{0,0,0}$\y_1=v(x)$}%
}}}
\put(3376,-4936){\makebox(0,0)[lb]{\smash{\SetFigFont{22}{14.4}{\rmdefault}{\mddefault}{\updefault}{\color[rgb]{0,0,0}$x$}%
}}}
\put(5326,-3736){\makebox(0,0)[lb]{\smash{\SetFigFont{22}{14.4}{\rmdefault}{\mddefault}{\updefault}{\color[rgb]{0,0,0}$M_C$}%
}}}
\put(751,-1861){\rotatebox{320.0}{\makebox(0,0)[lb]{\smash{\SetFigFont{22}{14.4}{\rmdefault}{\mddefault}{\updefault}{\color[rgb]{0,0,0}$\ell(v(x),x)$}%
}}}}
\put(2992,143){\rotatebox{315.0}{\makebox(0,0)[lb]{\smash{\SetFigFont{22}{14.4}{\rmdefault}{\mddefault}{\updefault}{\color[rgb]{0,0,0}$\ell_2(v(x),x)$}%
}}}}
\put(4876,539){\makebox(0,0)[lb]{\smash{\SetFigFont{22}{14.4}{\rmdefault}{\mddefault}{\updefault}{\color[rgb]{0,0,0}$\y_3$}%
}}}
\put(6975,-1389){\rotatebox{310.0}{\makebox(0,0)[lb]{\smash{\SetFigFont{22}{14.4}{\rmdefault}{\mddefault}{\updefault}{\color[rgb]{0,0,0}$e(x)$}%
}}}}
\put(11251,-6136){\makebox(0,0)[lb]{\smash{\SetFigFont{22}{14.4}{\rmdefault}{\mddefault}{\updefault}{\color[rgb]{0,0,0}$\y_2$}%
}}}
\put(2926,-7186){\rotatebox{45.0}{\makebox(0,0)[lb]{\smash{\SetFigFont{22}{14.4}{\rmdefault}{\mddefault}{\updefault}{\color[rgb]{0,0,0}$d(v(x),\ell_2(v(x),x))=2\,d(v(x),\ell(v(x),x))$}%
}}}}
\put(960,-5299){\rotatebox{45.0}{\makebox(0,0)[lb]{\smash{\SetFigFont{22}{14.4}{\rmdefault}{\mddefault}{\updefault}{\color[rgb]{0,0,0} $d(v(x),\ell(v(x),x))$}%
}}}}
\end{picture}

%% file: Gammaofnu2.pstex_t
\begin{picture}(0,0)%
\includegraphics{Gammaofnu2.pstex}%
\end{picture}%
\setlength{\unitlength}{3947sp}%
\begingroup\makeatletter\ifx\SetFigFont\undefined%
\gdef\SetFigFont#1#2#3#4#5{%
  \reset@font\fontsize{#1}{#2pt}%
  \fontfamily{#3}\fontseries{#4}\fontshape{#5}%
  \selectfont}%
\fi\endgroup%
\begin{picture}(10684,8043)(589,-7348)
\put(1051,-6511){\makebox(0,0)[lb]{\smash{\SetFigFont{22}{14.4}{\rmdefault}{\mddefault}{\updefault}{\color[rgb]{0,0,0}$\y_1$}%
}}}
\put(4876,539){\makebox(0,0)[lb]{\smash{\SetFigFont{22}{14.4}{\rmdefault}{\mddefault}{\updefault}{\color[rgb]{0,0,0}$\y_3$}%
}}}
\put(11300,-6136){\makebox(0,0)[lb]{\smash{\SetFigFont{22}{14.4}{\rmdefault}{\mddefault}{\updefault}{\color[rgb]{0,0,0}$\y_2$}%
}}}
\put(1405,-2508){\rotatebox{320.0}{\makebox(0,0)[lb]{\smash{\SetFigFont{22}{14.4}{\rmdefault}{\mddefault}{\updefault}{\color[rgb]{0,0,0}$\ell(\y_1,x)$}%
}}}}
\put(5701,-4036){\makebox(0,0)[lb]{\smash{\SetFigFont{22}{14.4}{\rmdefault}{\mddefault}{\updefault}{\color[rgb]{0,0,0}$M_C$}%
}}}
\put(3376,-4936){\makebox(0,0)[lb]{\smash{\SetFigFont{22}{14.4}{\rmdefault}{\mddefault}{\updefault}{\color[rgb]{0,0,0}$x$}%
}}}
\put(2626,-2011){\makebox(0,0)[lb]{\smash{\SetFigFont{22}{14.4}{\rmdefault}{\mddefault}{\updefault}{\color[rgb]{0,0,0}$\xi_3(2,x)$}%
}}}
\put(7533,-5107){\rotatebox{65.0}{\makebox(0,0)[lb]{\smash{\SetFigFont{22}{14.4}{\rmdefault}{\mddefault}{\updefault}{\color[rgb]{0,0,0}$\xi_2(2,x)$}%
}}}}
\put(960,-5299){\rotatebox{45.0}{\makebox(0,0)[lb]{\smash{\SetFigFont{22}{14.4}{\rmdefault}{\mddefault}{\updefault}{\color[rgb]{0,0,0}$d(\y_1,\ell(\y_1,x))=r\,d(\y_1,\xi_1(2,x))$}%
}}}}
\put(2026,-7036){\rotatebox{45.0}{\makebox(0,0)[lb]{\smash{\SetFigFont{22}{14.4}{\rmdefault}{\mddefault}{\updefault}{\color[rgb]{0,0,0}$d(\y_1,\xi_1(2,x))$}%
}}}}
\put(3867,-6839){\rotatebox{320.0}{\makebox(0,0)[lb]{\smash{\SetFigFont{22}{14.4}{\rmdefault}{\mddefault}{\updefault}{\color[rgb]{0,0,0}$\xi_1(2,x)$}%
}}}}
\end{picture}

%% file: two_extrema.pstex_t
\begin{picture}(0,0)%
\includegraphics{two_extrema.pstex}%
\end{picture}%
\setlength{\unitlength}{3947sp}%
\begingroup\makeatletter\ifx\SetFigFont\undefined%
\gdef\SetFigFont#1#2#3#4#5{%
  \reset@font\fontsize{#1}{#2pt}%
  \fontfamily{#3}\fontseries{#4}\fontshape{#5}%
  \selectfont}%
\fi\endgroup%
\begin{picture}(11349,8124)(-311,-8023)
\put(9301,-6586){\makebox(0,0)[lb]{\smash{{\SetFigFont{25}{16.8}{\rmdefault}{\mddefault}{\updefault}{\color[rgb]{0,0,0}$\y_2=(1,0)$}%
}}}}
\put(826,-6586){\makebox(0,0)[lb]{\smash{{\SetFigFont{25}{16.8}{\rmdefault}{\mddefault}{\updefault}{\color[rgb]{0,0,0}$\y_1=(0,0)$}%
}}}}
\put(3976,-511){\makebox(0,0)[lb]{\smash{{\SetFigFont{25}{16.8}{\rmdefault}{\mddefault}{\updefault}{\color[rgb]{0,0,0}$\y_3=(c_1,c_2)$}%
}}}}
\put(4276,-2761){\makebox(0,0)[lb]{\smash{{\SetFigFont{25}{16.8}{\rmdefault}{\mddefault}{\updefault}{\color[rgb]{0,0,0}$(x,y)$}%
}}}}
\end{picture}%

%% file: two_ext_F.pstex_t
\begin{picture}(0,0)%
\includegraphics{two_ext_F.pstex}%
\end{picture}%
\setlength{\unitlength}{3947sp}%
\begingroup\makeatletter\ifx\SetFigFont\undefined%
\gdef\SetFigFont#1#2#3#4#5{%
  \reset@font\fontsize{#1}{#2pt}%
  \fontfamily{#3}\fontseries{#4}\fontshape{#5}%
  \selectfont}%
\fi\endgroup%
\begin{picture}(11349,8124)(-311,-8023)
\put(9301,-6586){\makebox(0,0)[lb]{\smash{{\SetFigFont{22}{16.8}{\rmdefault}{\mddefault}{\updefault}{\color[rgb]{0,0,0}$y_2=(1,0)$}%
}}}}
\put(826,-6586){\makebox(0,0)[lb]{\smash{{\SetFigFont{22}{16.8}{\rmdefault}{\mddefault}{\updefault}{\color[rgb]{0,0,0}$y_1=(0,0)$}%
}}}}
\put(3976,-511){\makebox(0,0)[lb]{\smash{{\SetFigFont{22}{16.8}{\rmdefault}{\mddefault}{\updefault}{\color[rgb]{0,0,0}$y_3=(c_1,c_2)$}%
}}}}
\put(2401,-5161){\makebox(0,0)[lb]{\smash{{\SetFigFont{22}{16.8}{\rmdefault}{\mddefault}{\updefault}{\color[rgb]{0,0,0}$(x,y)$}%
}}}}
\end{picture}%

%% file: N_CSexample2.pstex_t
\begin{picture}(0,0)%
\includegraphics{N_CSexample2.pstex}%
\end{picture}%
\setlength{\unitlength}{3947sp}%
\begingroup\makeatletter\ifx\SetFigFont\undefined%
\gdef\SetFigFont#1#2#3#4#5{%
  \reset@font\fontsize{#1}{#2pt}%
  \fontfamily{#3}\fontseries{#4}\fontshape{#5}%
  \selectfont}%
\fi\endgroup%
\begin{picture}(10244,7264)(1051,-6569)
\put(6226,-4711){\makebox(0,0)[lb]{\smash{\SetFigFont{22}{14.4}{\rmdefault}{\mddefault}{\updefault}{$x$}%
}}}
\put(1051,-6511){\makebox(0,0)[lb]{\smash{\SetFigFont{22}{14.4}{\rmdefault}{\mddefault}{\updefault}{$\y_1$}%
}}}
\put(4876,539){\makebox(0,0)[lb]{\smash{\SetFigFont{22}{14.4}{\rmdefault}{\mddefault}{\updefault}{$\y_3$}%
}}}
\put(2551,-2761){\makebox(0,0)[lb]{\smash{\SetFigFont{22}{14.4}{\rmdefault}{\mddefault}{\updefault}{$e_2$}%
}}}
\put(2251,-5911){\makebox(0,0)[lb]{\smash{\SetFigFont{22}{14.4}{\rmdefault}{\mddefault}{\updefault}{$d(x,e_3^{\tau}(x))=\tau\,d(x,e(x))$}%
}}}
\put(4726,-3211){\makebox(0,0)[lb]{\smash{\SetFigFont{22}{14.4}{\rmdefault}{\mddefault}{\updefault}{$M_C$}%
}}}
\put(5626,-5611){\makebox(0,0)[lb]{\smash{\SetFigFont{22}{14.4}{\rmdefault}{\mddefault}{\updefault}{$e_3^{\tau}(x)$}%
}}}
\put(6751,-3811){\makebox(0,0)[lb]{\smash{\SetFigFont{22}{14.4}{\rmdefault}{\mddefault}{\updefault}{$e_1^{\tau}(x)$}%
}}}
\put(8926,-5836){\makebox(0,0)[lb]{\smash{\SetFigFont{22}{14.4}{\rmdefault}{\mddefault}{\updefault}{$d(x,e(x))$}%
}}}
\put(5551,-6511){\makebox(0,0)[lb]{\smash{\SetFigFont{22}{14.4}{\rmdefault}{\mddefault}{\updefault}{$e_3=e(x)$}%
}}}
\put(8176,-2836){\makebox(0,0)[lb]{\smash{\SetFigFont{22}{14.4}{\rmdefault}{\mddefault}{\updefault}{$e_1$}%
}}}
\put(11251,-6136){\makebox(0,0)[lb]{\smash{\SetFigFont{22}{14.4}{\rmdefault}{\mddefault}{\updefault}{$\y_2$}%
}}}
\put(3451,-4036){\makebox(0,0)[lb]{\smash{\SetFigFont{22}{14.4}{\rmdefault}{\mddefault}{\updefault}{$e_2^{\tau}(x)$}%
}}}
\end{picture}

%% file: N_CSGam1.pstex_t
\begin{picture}(0,0)%
\includegraphics{N_CSGam1.pstex}%
\end{picture}%
\setlength{\unitlength}{3947sp}%
\begingroup\makeatletter\ifx\SetFigFont\undefined%
\gdef\SetFigFont#1#2#3#4#5{%
  \reset@font\fontsize{#1}{#2pt}%
  \fontfamily{#3}\fontseries{#4}\fontshape{#5}%
  \selectfont}%
\fi\endgroup%
\begin{picture}(10738,7297)(1051,-6578)
\put(5551,-6511){\makebox(0,0)[lb]{\smash{{\SetFigFont{25}{16.8}{\rmdefault}{\mddefault}{\updefault}{\color[rgb]{0,0,0}$e_3=e(x)$}%
}}}}
\put(5776,-5536){\makebox(0,0)[lb]{\smash{{\SetFigFont{23}{16.8}{\rmdefault}{\bfdefault}{\updefault}{\color[rgb]{0,0,0}$\zeta_2(\tau,x)$}%
}}}}
\put(5026,-5236){\makebox(0,0)[lb]{\smash{{\SetFigFont{23}{16.8}{\rmdefault}{\bfdefault}{\updefault}{\color[rgb]{0,0,0}$\zeta_7(\tau,x)$}%
}}}}
\put(8176,-2836){\makebox(0,0)[lb]{\smash{{\SetFigFont{23}{16.8}{\rmdefault}{\mddefault}{\updefault}{\color[rgb]{0,0,0}$e_1$}%
}}}}
\put(11251,-6136){\makebox(0,0)[lb]{\smash{{\SetFigFont{23}{16.8}{\rmdefault}{\bfdefault}{\updefault}{\color[rgb]{0,0,0}$\mathsf{y}_2$}%
}}}}
\put(4876,539){\makebox(0,0)[lb]{\smash{{\SetFigFont{23}{16.8}{\rmdefault}{\bfdefault}{\updefault}{\color[rgb]{0,0,0}$\mathsf{y}_3$}%
}}}}
\put(1051,-6511){\makebox(0,0)[lb]{\smash{{\SetFigFont{23}{16.8}{\rmdefault}{\bfdefault}{\updefault}{\color[rgb]{0,0,0}$\mathsf{y}_1$}%
}}}}
\put(2551,-2761){\makebox(0,0)[lb]{\smash{{\SetFigFont{23}{16.8}{\rmdefault}{\mddefault}{\updefault}{\color[rgb]{0,0,0}$e_2$}%
}}}}
\put(5476,-5761){\makebox(0,0)[lb]{\smash{{\SetFigFont{18}{14.4}{\rmdefault}{\bfdefault}{\updefault}{\color[rgb]{0,0,0}$$}%
}}}}
\put(4201,-5836){\makebox(0,0)[lb]{\smash{{\SetFigFont{23}{16.8}{\rmdefault}{\bfdefault}{\updefault}{\color[rgb]{0,0,0}$\zeta_1(\tau,x)$}%
}}}}
\end{picture}%

%% file: N_CSGam2.pstex_t
\begin{picture}(0,0)%
\includegraphics{N_CSGam2.pstex}%
\end{picture}%
\setlength{\unitlength}{3947sp}%
\begingroup\makeatletter\ifx\SetFigFont\undefined%
\gdef\SetFigFont#1#2#3#4#5{%
  \reset@font\fontsize{#1}{#2pt}%
  \fontfamily{#3}\fontseries{#4}\fontshape{#5}%
  \selectfont}%
\fi\endgroup%
\begin{picture}(10738,7297)(1051,-6578)
\put(5551,-6511){\makebox(0,0)[lb]{\smash{{\SetFigFont{25}{16.8}{\rmdefault}{\mddefault}{\updefault}{\color[rgb]{0,0,0}$e_3=e(x)$}%
}}}}
\put(8176,-2836){\makebox(0,0)[lb]{\smash{{\SetFigFont{23}{16.8}{\rmdefault}{\mddefault}{\updefault}{\color[rgb]{0,0,0}$e_1$}%
}}}}
\put(11251,-6136){\makebox(0,0)[lb]{\smash{{\SetFigFont{23}{16.8}{\rmdefault}{\bfdefault}{\updefault}{\color[rgb]{0,0,0}$\mathsf{y}_2$}%
}}}}
\put(4876,539){\makebox(0,0)[lb]{\smash{{\SetFigFont{23}{16.8}{\rmdefault}{\bfdefault}{\updefault}{\color[rgb]{0,0,0}$\mathsf{y}_3$}%
}}}}
\put(1051,-6511){\makebox(0,0)[lb]{\smash{{\SetFigFont{23}{16.8}{\rmdefault}{\bfdefault}{\updefault}{\color[rgb]{0,0,0}$\mathsf{y}_1$}%
}}}}
\put(2551,-2761){\makebox(0,0)[lb]{\smash{{\SetFigFont{23}{16.8}{\rmdefault}{\mddefault}{\updefault}{\color[rgb]{0,0,0}$e_2$}%
}}}}
\put(3376,-5611){\makebox(0,0)[lb]{\smash{{\SetFigFont{18}{14.4}{\rmdefault}{\bfdefault}{\updefault}{\color[rgb]{0,0,0}$$}%
}}}}
\put(2026,-5686){\makebox(0,0)[lb]{\smash{{\SetFigFont{23}{16.8}{\rmdefault}{\bfdefault}{\updefault}{\color[rgb]{0,0,0}$\zeta_1(\tau,x)$}%
}}}}
\put(4051,-5536){\makebox(0,0)[lb]{\smash{{\SetFigFont{23}{16.8}{\rmdefault}{\bfdefault}{\updefault}{\color[rgb]{0,0,0}$\zeta_2(\tau,x)$}%
}}}}
\put(4201,-4786){\makebox(0,0)[lb]{\smash{{\SetFigFont{23}{16.8}{\rmdefault}{\bfdefault}{\updefault}{\color[rgb]{0,0,0}$\zeta_7(\tau,x)$}%
}}}}
\put(2476,-4486){\makebox(0,0)[lb]{\smash{{\SetFigFont{23}{16.8}{\rmdefault}{\bfdefault}{\updefault}{\color[rgb]{0,0,0}$\zeta_5(\tau,x)$}%
}}}}
\put(1651,-5161){\makebox(0,0)[lb]{\smash{{\SetFigFont{23}{16.8}{\rmdefault}{\bfdefault}{\updefault}{\color[rgb]{0,0,0}$\zeta_6(\tau,x)$}%
}}}}
\end{picture}%

%% file: N_CSGam3.pstex_t
\begin{picture}(0,0)%
\includegraphics{N_CSGam3.pstex}%
\end{picture}%
\setlength{\unitlength}{3947sp}%
\begingroup\makeatletter\ifx\SetFigFont\undefined%
\gdef\SetFigFont#1#2#3#4#5{%
  \reset@font\fontsize{#1}{#2pt}%
  \fontfamily{#3}\fontseries{#4}\fontshape{#5}%
  \selectfont}%
\fi\endgroup%
\begin{picture}(10738,7297)(1051,-6578)
\put(5551,-6511){\makebox(0,0)[lb]{\smash{{\SetFigFont{25}{16.8}{\rmdefault}{\mddefault}{\updefault}{\color[rgb]{0,0,0}$e_3=e(x)$}%
}}}}
\put(8176,-2836){\makebox(0,0)[lb]{\smash{{\SetFigFont{23}{16.8}{\rmdefault}{\mddefault}{\updefault}{\color[rgb]{0,0,0}$e_1$}%
}}}}
\put(11251,-6136){\makebox(0,0)[lb]{\smash{{\SetFigFont{23}{16.8}{\rmdefault}{\bfdefault}{\updefault}{\color[rgb]{0,0,0}$\mathsf{y}_2$}%
}}}}
\put(4876,539){\makebox(0,0)[lb]{\smash{{\SetFigFont{23}{16.8}{\rmdefault}{\bfdefault}{\updefault}{\color[rgb]{0,0,0}$\mathsf{y}_3$}%
}}}}
\put(1051,-6511){\makebox(0,0)[lb]{\smash{{\SetFigFont{23}{16.8}{\rmdefault}{\bfdefault}{\updefault}{\color[rgb]{0,0,0}$\mathsf{y}_1$}%
}}}}
\put(2551,-2761){\makebox(0,0)[lb]{\smash{{\SetFigFont{23}{16.8}{\rmdefault}{\mddefault}{\updefault}{\color[rgb]{0,0,0}$e_2$}%
}}}}
\put(3601,-5311){\makebox(0,0)[lb]{\smash{{\SetFigFont{23}{16.8}{\rmdefault}{\bfdefault}{\updefault}{\color[rgb]{0,0,0}$\zeta_1(\tau,x)$}%
}}}}
\put(6001,-5311){\makebox(0,0)[lb]{\smash{{\SetFigFont{23}{16.8}{\rmdefault}{\bfdefault}{\updefault}{\color[rgb]{0,0,0}$\zeta_2(\tau,x)$}%
}}}}
\put(5326,-5311){\makebox(0,0)[lb]{\smash{{\SetFigFont{18}{14.4}{\rmdefault}{\bfdefault}{\updefault}{\color[rgb]{0,0,0}$$}%
}}}}
\put(3826,-3811){\makebox(0,0)[lb]{\smash{{\SetFigFont{23}{16.8}{\rmdefault}{\bfdefault}{\updefault}{\color[rgb]{0,0,0}$\zeta_5(\tau,x)$}%
}}}}
\put(2926,-4411){\makebox(0,0)[lb]{\smash{{\SetFigFont{23}{16.8}{\rmdefault}{\bfdefault}{\updefault}{\color[rgb]{0,0,0}$\zeta_6(\tau,x)$}%
}}}}
\put(5176,-3661){\makebox(0,0)[lb]{\smash{{\SetFigFont{23}{16.8}{\rmdefault}{\bfdefault}{\updefault}{\color[rgb]{0,0,0}$\zeta_7(\tau,x)$}%
}}}}
\put(7126,-4336){\makebox(0,0)[lb]{\smash{{\SetFigFont{23}{16.8}{\rmdefault}{\bfdefault}{\updefault}{\color[rgb]{0,0,0}$\zeta_3(\tau,x)$}%
}}}}
\put(6526,-3811){\makebox(0,0)[lb]{\smash{{\SetFigFont{23}{16.8}{\rmdefault}{\bfdefault}{\updefault}{\color[rgb]{0,0,0}$\zeta_4(\tau,x)$}%
}}}}
\end{picture}%

%% file: N_CSGam4.pstex_t
\begin{picture}(0,0)%
\includegraphics{N_CSGam4.pstex}%
\end{picture}%
\setlength{\unitlength}{3947sp}%
\begingroup\makeatletter\ifx\SetFigFont\undefined%
\gdef\SetFigFont#1#2#3#4#5{%
  \reset@font\fontsize{#1}{#2pt}%
  \fontfamily{#3}\fontseries{#4}\fontshape{#5}%
  \selectfont}%
\fi\endgroup%
\begin{picture}(10738,7297)(1051,-6578)
\put(5551,-6511){\makebox(0,0)[lb]{\smash{{\SetFigFont{25}{16.8}{\rmdefault}{\mddefault}{\updefault}{\color[rgb]{0,0,0}$e_3=e(x)$}%
}}}}
\put(8176,-2836){\makebox(0,0)[lb]{\smash{{\SetFigFont{23}{16.8}{\rmdefault}{\mddefault}{\updefault}{\color[rgb]{0,0,0}$e_1$}%
}}}}
\put(11261,-6136){\makebox(0,0)[lb]{\smash{{\SetFigFont{23}{16.8}{\rmdefault}{\bfdefault}{\updefault}{\color[rgb]{0,0,0}$\mathsf{y}_2$}%
}}}}
\put(4876,539){\makebox(0,0)[lb]{\smash{{\SetFigFont{23}{16.8}{\rmdefault}{\bfdefault}{\updefault}{\color[rgb]{0,0,0}$\mathsf{y}_3$}%
}}}}
\put(1051,-6511){\makebox(0,0)[lb]{\smash{{\SetFigFont{23}{16.8}{\rmdefault}{\bfdefault}{\updefault}{\color[rgb]{0,0,0}$\mathsf{y}_1$}%
}}}}
\put(2551,-2761){\makebox(0,0)[lb]{\smash{{\SetFigFont{23}{16.8}{\rmdefault}{\mddefault}{\updefault}{\color[rgb]{0,0,0}$e_2$}%
}}}}
\put(3826,-5011){\makebox(0,0)[lb]{\smash{{\SetFigFont{23}{16.8}{\rmdefault}{\bfdefault}{\updefault}{\color[rgb]{0,0,0}$\zeta_1(\tau,x)$}%
}}}}
\put(5926,-4786){\makebox(0,0)[lb]{\smash{{\SetFigFont{23}{16.8}{\rmdefault}{\bfdefault}{\updefault}{\color[rgb]{0,0,0}$\zeta_2(\tau,x)$}%
}}}}
\put(5776,-3286){\makebox(0,0)[lb]{\smash{{\SetFigFont{23}{16.8}{\rmdefault}{\bfdefault}{\updefault}{\color[rgb]{0,0,0}$\zeta_4(\tau,x)$}%
}}}}
\put(6301,-3961){\makebox(0,0)[lb]{\smash{{\SetFigFont{23}{16.8}{\rmdefault}{\bfdefault}{\updefault}{\color[rgb]{0,0,0}$\zeta_3(\tau,x)$}%
}}}}
\put(4051,-3361){\makebox(0,0)[lb]{\smash{{\SetFigFont{23}{16.8}{\rmdefault}{\bfdefault}{\updefault}{\color[rgb]{0,0,0}$\zeta_5(\tau,x)$}%
}}}}
\put(3226,-4036){\makebox(0,0)[lb]{\smash{{\SetFigFont{23}{16.8}{\rmdefault}{\bfdefault}{\updefault}{\color[rgb]{0,0,0}$\zeta_6(\tau,x)$}%
}}}}
\put(5401,-4561){\makebox(0,0)[lb]{\smash{{\SetFigFont{18}{14.4}{\rmdefault}{\bfdefault}{\updefault}{\color[rgb]{0,0,0}$$}%
}}}}
\end{picture}%

%% file: Tri_Gamma.pstex_t
\begin{picture}(0,0)%
\includegraphics{Tri_Gamma.pstex}%
\end{picture}%
\setlength{\unitlength}{3947sp}%
\begingroup\makeatletter\ifx\SetFigFont\undefined%
\gdef\SetFigFont#1#2#3#4#5{%
  \reset@font\fontsize{#1}{#2pt}%
  \fontfamily{#3}\fontseries{#4}\fontshape{#5}%
  \selectfont}%
\fi\endgroup%
\begin{picture}(8669,7189)(1951,-6569)
\put(7727,-1397){\rotatebox{300.0}{\makebox(0,0)[lb]{\smash{{\SetFigFont{22}{14.4}{\rmdefault}{\mddefault}{\updefault}{\color[rgb]{0,0,0}$e_1(x)$}%
}}}}}
\put(6301,-3586){\makebox(0,0)[lb]{\smash{{\SetFigFont{22}{14.4}{\rmdefault}{\mddefault}{\updefault}{\color[rgb]{0,0,0}$M_C$}%
}}}}
\put(1951,-6511){\makebox(0,0)[lb]{\smash{{\SetFigFont{22}{14.4}{\rmdefault}{\mddefault}{\updefault}{\color[rgb]{0,0,0}$\y_1=(0,0)$}%
}}}}
\put(5701,464){\makebox(0,0)[lb]{\smash{{\SetFigFont{22}{14.4}{\rmdefault}{\mddefault}{\updefault}{\color[rgb]{0,0,0}$\y_3=\left(1/2,\sqrt{3}/2 \right)$}%
}}}}
\put(9076,-6511){\makebox(0,0)[lb]{\smash{{\SetFigFont{22}{14.4}{\rmdefault}{\mddefault}{\updefault}{\color[rgb]{0,0,0}$\y_2=(1,0)$}%
}}}}
\put(4501,-4786){\makebox(0,0)[lb]{\smash{{\SetFigFont{22}{14.4}{\rmdefault}{\mddefault}{\updefault}{\color[rgb]{0,0,0}$x_{e_3}=(x_3,y_3)$}%
}}}}
\put(7501,-2611){\makebox(0,0)[lb]{\smash{{\SetFigFont{22}{14.4}{\rmdefault}{\mddefault}{\updefault}{\color[rgb]{0,0,0}$x_{e_2}=(x_2,y_2)$}%
}}}}
\put(7276,-5611){\makebox(0,0)[lb]{\smash{{\SetFigFont{22}{14.4}{\rmdefault}{\mddefault}{\updefault}{\color[rgb]{0,0,0}$x_{e_1}=(x_1,y_1)$}%
}}}}
\end{picture}%

%% file: jntpdf1.pstex_t
\begin{picture}(0,0)%
\includegraphics{jntpdf1.pstex}%
\end{picture}%
\setlength{\unitlength}{3947sp}%
\begingroup\makeatletter\ifx\SetFigFont\undefined%
\gdef\SetFigFont#1#2#3#4#5{%
  \reset@font\fontsize{#1}{#2pt}%
  \fontfamily{#3}\fontseries{#4}\fontshape{#5}%
  \selectfont}%
\fi\endgroup%
\begin{picture}(8369,7189)(2251,-6569)
\put(6301,-3586){\makebox(0,0)[lb]{\smash{{\SetFigFont{22}{14.4}{\rmdefault}{\mddefault}{\updefault}{\color[rgb]{0,0,0}$M_C$}%
}}}}
\put(7727,-1397){\rotatebox{300.0}{\makebox(0,0)[lb]{\smash{{\SetFigFont{22}{14.4}{\rmdefault}{\mddefault}{\updefault}{\color[rgb]{0,0,0}$e_1(x)$}%
}}}}}
\put(5926,464){\makebox(0,0)[lb]{\smash{{\SetFigFont{22}{14.4}{\rmdefault}{\mddefault}{\updefault}{\color[rgb]{0,0,0}$\y_3=\left(1/2,\sqrt{3}/2\right)$}%
}}}}
\put(9676,-6511){\makebox(0,0)[lb]{\smash{{\SetFigFont{22}{14.4}{\rmdefault}{\mddefault}{\updefault}{\color[rgb]{0,0,0}$\y_2=(1,0)$}%
}}}}
\put(2251,-6436){\makebox(0,0)[lb]{\smash{{\SetFigFont{22}{14.4}{\rmdefault}{\mddefault}{\updefault}{\color[rgb]{0,0,0}$\y_1=(0,0)$}%
}}}}
\put(7876,-2386){\makebox(0,0)[lb]{\smash{{\SetFigFont{22}{14.4}{\rmdefault}{\mddefault}{\updefault}{\color[rgb]{0,0,0}$x_{e_2}=(x_2,y_2)$}%
}}}}
\put(6976,-6136){\makebox(0,0)[lb]{\smash{{\SetFigFont{22}{14.4}{\rmdefault}{\mddefault}{\updefault}{\color[rgb]{0,0,0}$x_{e_1}=(x_1,y_1)$}%
}}}}
\put(2476,-3511){\makebox(0,0)[lb]{\smash{{\SetFigFont{22}{14.4}{\rmdefault}{\mddefault}{\updefault}{\color[rgb]{0,0,0}$x_{e_3}=(x_3,y_3)$}%
}}}}
\end{picture}%